\documentclass[symmetry,article,submit,moreauthors,pdftex,12pt,a4paper]{mdpi}

\lastpage{x}
\doinum{10.3390/------}
\pubvolume{xx}
\pubyear{2015}
\history{Received: xx / Accepted: xx / Published: xx}

\usepackage{amsfonts}
\usepackage{amssymb}
\usepackage{amsthm}
\usepackage{graphicx}

\theoremstyle{mdpi}
\newcounter{thm}
\newcounter{ex}
\newcounter{re}
\newcommand{\R}{\mathbb{R}}

\Title{Some Elementary Aspects of $4$-dimensional Geometry}

\Author{J. Scott Carter $^{1,}$* and David A. Mullens $^{1}$}

\address{%
$^{1}$ University of South Alabama, Department of Mathematics and Statistics, Mobile, AL 36688}

\corres{carter@southalabama.edu, 251-460-6264, Fax 251-460-7969}

\abstract{We indicate that Heron's formula (which relates the square of the area of a triangle to a quartic function of its edge lengths) can be interpreted as a scissors congruence  in $4$-dimensional space. In the process of demonstrating this, we examine a number of decompositions of hypercubes, hyper-parallelograms,
and other elementary $4$-dimensional solids. }


\begin{document}

\section{Introduction}

Heron's formula gives the explicit relationship among the area, $A$,  of a triangle and its edge lengths. We assume, now and throughout, that a triangle with edge lengths, $a,b,c$ where $0< a \le b \le c$ is given. The {\it semi-perimeter} is the quantity: $s=(a+b+c)/2$. One representation of the formula is 
$$A=\sqrt{s(s-a)(s-b)(s-c)}.$$
The formula was written in Heron's (Hero's) book Metrica circa $60$AD. According to wikipedia \cite{wiki}, the formula may have been known to Archimedes more than two centuries earlier. But Heron's book was a compendium of a number of known results. Here, we square both sides of the equation and multiply by the denominator to get the following expression
$$16 A^2 =(a+b+c)(-a+b+c)(a-b+c)(a+b-c).$$
This can be considered an expression that relates the $4$-dimensional volumes of a number of hyper-solids. As the text develops, we will further manipulate this expression all the while maintaining a $4$-dimensional awareness of the algebraic manipulations that are performed.

Our primary purpose here is to indicate a $4$-dimensional scissors congruence among the hyper-solids that are involved in Heron's formula. But in the process, we will also develop several secondary goals along the way. These are our own interpretations of conversations with several people. Most importantly, we hope to develop the reader's intuition and interest in elementary $4$-dimensional geometry.

In particular, many algebraic identities, such as the binomial and multi-nomial theorems, have elegant and useful interpretations in terms of relating the higher dimensional volumes of various $n$-cubes. The symmetries in the expressions are manifest in the symmetries among the pieces of these decompositions. 

\begin{figure}
\begin{center}
\includegraphics[scale=.4]{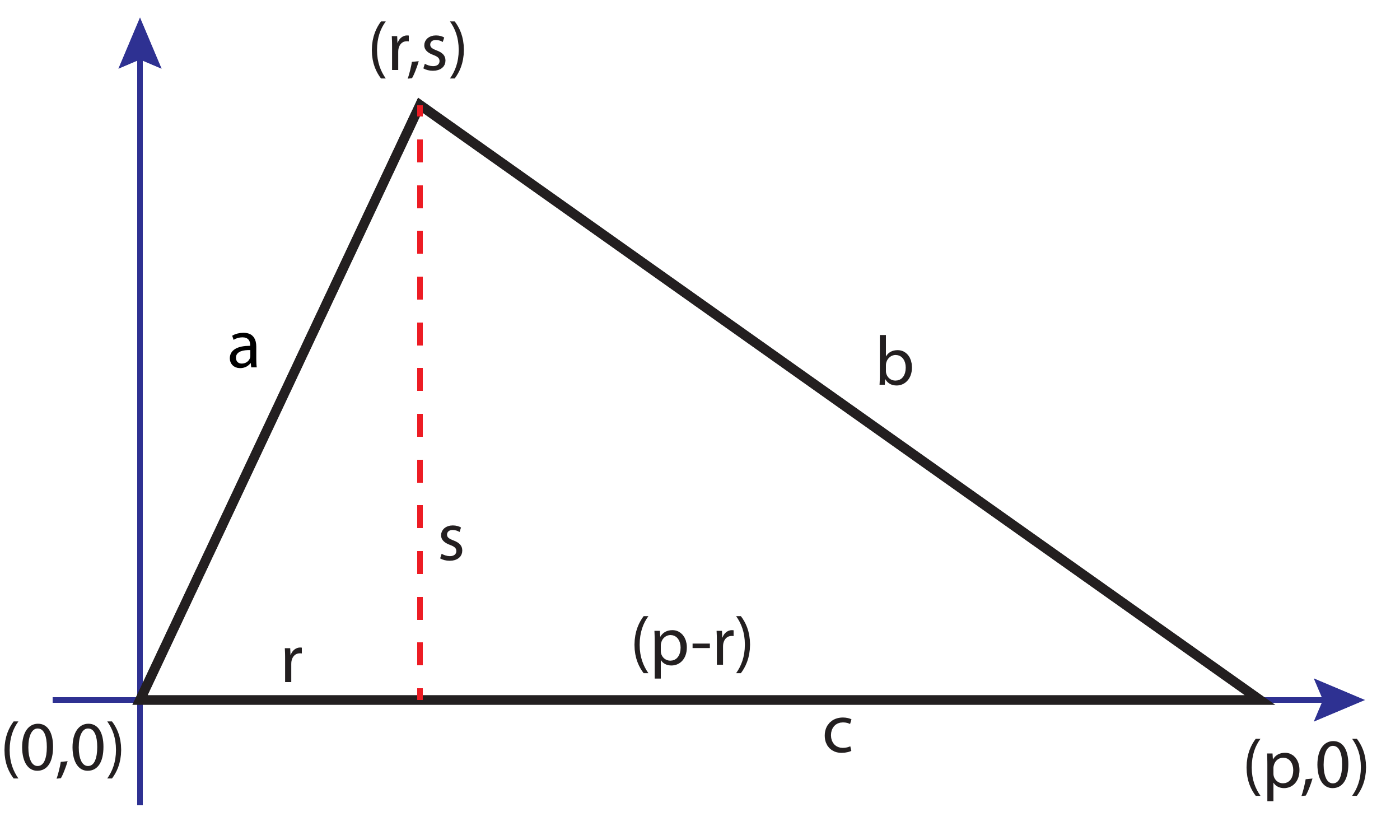}
\end{center}
\caption{A standard triangle of side lengths $a\le b \le c$}
\label{righttriangle}
\end{figure}

We remark here that the algebraic proof of Heron's formula given in this form goes as follows. 
Suppose that we arrange the triangle in the coordinate plane as indicated in Fig.~\ref{righttriangle}.
In this way, one vertex is at the origin, the longest side (the side of length $c$) lies along the $x$-axis with a vertex at the point $(p,0)$ (Here a different coordinate is chosen so that  in the new variables are of the same flavor). The remaining vertex lies at the point $(r,s)$. In this way, $a^2=r^2+s^2$, and $b^2=s^2+(p-r)^2$. The left-hand-side of Heron's formula can be shown to be $4 p^2s^2.$ The right-hand-side reduces to 
$${\mbox{\rm RHS}}=2a^2b^2+2a^2c^2+2b^2c^2 -a^4-b^4-c^4.$$
 We expand the terms on the right-hand-side as follows:
 \begin{eqnarray*} 
 2a^2b^2 & =& 2(r^2+s^2)((p-r)^2+s^2) \\ 
 & =& 2r^2(p-r)^2 +2s^2(p-r)^2 +2r^2s^2 +2s^4. \\ & & \\
  2a^2c^2 & =& 2(r^2+s^2)p^2 \\ 
 & =& 2r^2p^2+2s^2p^2. \\ & & \\
  2b^2c^2 & =& 2((p-r)^2+s^2)p^2 \\ 
 & =& 2(p-r)^2p^2 +2s^2p^2. \\ & & \\
 -a^4 & =& -(r^2+s^2)(r^2+s^2) \\ 
 & =& -r^4-2r^2s^2 -s^4. \\ & & \\
 -b^4 & =& -((p-r)^2+s^2)((p-r)^2+s^2) \\ 
 & =& -(p-r)^4-2s^2(p-r)^2 -s^4. \\ & & \\
-c^4 & =& -p^4. 
\end{eqnarray*}
We cancel like terms:
 \begin{eqnarray*} 
RHS & =& 2r^2(p-r)^2    \\
 & +& 2r^2p^2+2s^2p^2 \\ 
 & +& 2(p-r)^2p^2 +2s^2p^2 \\
 &+&   (-r^4 ) \\
 & +& (-(p-r)^4 )  \\
 & +& (-p^4). 
 \end{eqnarray*}
After regrouping:
 \begin{eqnarray*} 
RHS & =&4s^2p^2      \\
 & +& 2r^2(p-r)^2+2(p-r)^2p^2 -(p-r)^4 \\  
&+& (-p^4) + 2r^2p^2 - r^4.
 \end{eqnarray*}
We factor the last line
 \begin{eqnarray*} 
RHS & =&4s^2p^2      \\
 & +& 2r^2(p-r)^2+2(p-r)^2p^2 -(p-r)^4 \\  
&+& (-(p^2-r^2)^2) \\ 
&=& 4s^2p^2      \\
 & +& 2r^2(p-r)^2+2(p-r)^2p^2 -(p-r)^4 \\  
&+& -(p-r)^2(p+r)^2.
 \end{eqnarray*}
 And combine like terms again:
  \begin{eqnarray*} 
RHS 
&=& 4s^2p^2      \\
 & +&(p-r)^2[ 2r^2+2p^2 -(p-r)^2- (p+r)^2] \\
 &=& 4s^2p^2.    
  \end{eqnarray*}
  
  As such the formula has been proven. However, we point out that even in our most quiet moments, the simplification $$(a+b+c)(a+b-c)(a-b+c)(-a+b+c)= 2a^2b^2+a^2c^2+b^2c^2-a^4-b^4-c^4$$ is quite tedious. After some preliminary work on developing higher dimensional models and the reader's intuition on $4$-dimensional geometry, our first task will be to show how to use $4$-dimensional geometry to encode the simplification. The second task that we endeavor upon is to capitalize upon a scissor's congruence proof of the Pythagorean Theorem, to expand the $4$-dimensional solids of size $a^4$, $b^4$, and $a^2b^2.$ 
Each of $a^2$ and $b^2$ is the square of a hypotenuse of a right triangle, and  can each be written as a sum of squares. To do so, one decomposes the square figure into five pieces and reassembles them. In $4$-dimensions, the quartics $a^4$, $b^4$, and $a^2b^2$ are decomposed into twenty-five pieces and reassembled as the sum of four hypercubes (or more precisely, hyper-rectangles). The grouping and factoring that is indicated above can also be seen in terms of $4$-dimensional scissor's congruences. Meanwhile, there is a little work to rearrange the factorization $4s^2p^2$ from four copies of a $4$-dimensional figure that is a parallelogram times itself into a hyper-rectangle.

In Section~\ref{ncube}, we discuss the $n$-cube from several points of view including its representation as a configuration space. Then we decompose it into right $n$-dimensional simplicies. This decomposition gives Nicomachus's Theorem \cite{Nico}. Our interpretation yields that the cube $[0,a]^4$ can be decomposed as the union of four congruent figures, each of the form a right isosceles triangle times itself. In Section~\ref{multi}, we discuss the multinomial theorem as another decomposition of the $n$-cube, and we use this idea to organize the computation of the right-hand-side of Heron's formula. In Section~\ref{pythag}, we demonstrate that the proof of the Pythagorean Theorem can be adjusted to a proof of a $4$-dimensional result: if $z^2=x^2+y^2$, and $u^2+v^2=w^2$, then the hyper-cube $[0,z]\times[0,z]\times[0,w]\times[0,w]$ can be cut into $25$ pieces that can be reassembled into a figure that is the product of a pair of sums of squares. This is the geometric interpretation of the identity
$$z^2w^2=x^2u^2+x^2v^2+y^2u^2+y^2v^2.$$
In Section~\ref{sumdiff}, we tie up some loose ends in the algebraic proof above. Section~\ref{LHS} discusses the left-hand-side of the equation, and Section~\ref{con} is a conclusion.

\section{The $n$-dimensional cube}
\label{ncube}

From a mathematician's perspective, the $n$-dimensional unit cube is a very natural object. It is written as $[0,1]^n$ and can be considered as the set of functions from $\{1,2,\ldots, n\}$ to the unit interval. A given function is determined by its set of values, and these are expressed as coordinates $(x_1,x_2, \ldots, x_n)$ where $x_j$ is the value of a given function on $j \in \{1,2,\ldots, n\}$. For  $j=1,2,\ldots, n,$ each $x_j$ is bounded between $0$ and $1$. We occasionally write $0\le x_1,x_2,\ldots, x_n \le 1$ to indicate the collection of such relations.
In the rest of this section, we develop a number of metaphors so that different intuitions can be developed. 

\subsection{Fingers}

 Before we consider the unit $n$-cube in general, let us consider binary digits --- your fingers, each in one of two possible states. Hold up your right hand and make a fist. Now sequentially start from your thumb and proceed through each of the remaining digits. Hold up your thumb as if hitch-hiking or simply giving a ``thumbs-up;'' retract it. Now hold up your index finger and retract it. Do so with each of your middle, annular, and your baby finger (pinky). Hold up some pair of fingers. In how many ways can two fingers be held up (10)? How many ways can you hold up three fingers (again, 10)? How about four fingers (5)? 
 
 These are digits since they are fingers. They are binary since they are either up or down. Assign the number $2^0=1$ to your thumb; assign $2^1=2$ to your index finger; assign  $2^2=4$ to your middle finger, assign $8$ to your ring finger, and assign $16$ to your pinky. When you consider all possible finger configurations on one hand, you can see that you can represent each integer from $0$ to $31$ as the sum of the values of the upheld fingers in each finger configuration.
 
\begin{figure}
\begin{center}
\includegraphics[scale=.04]{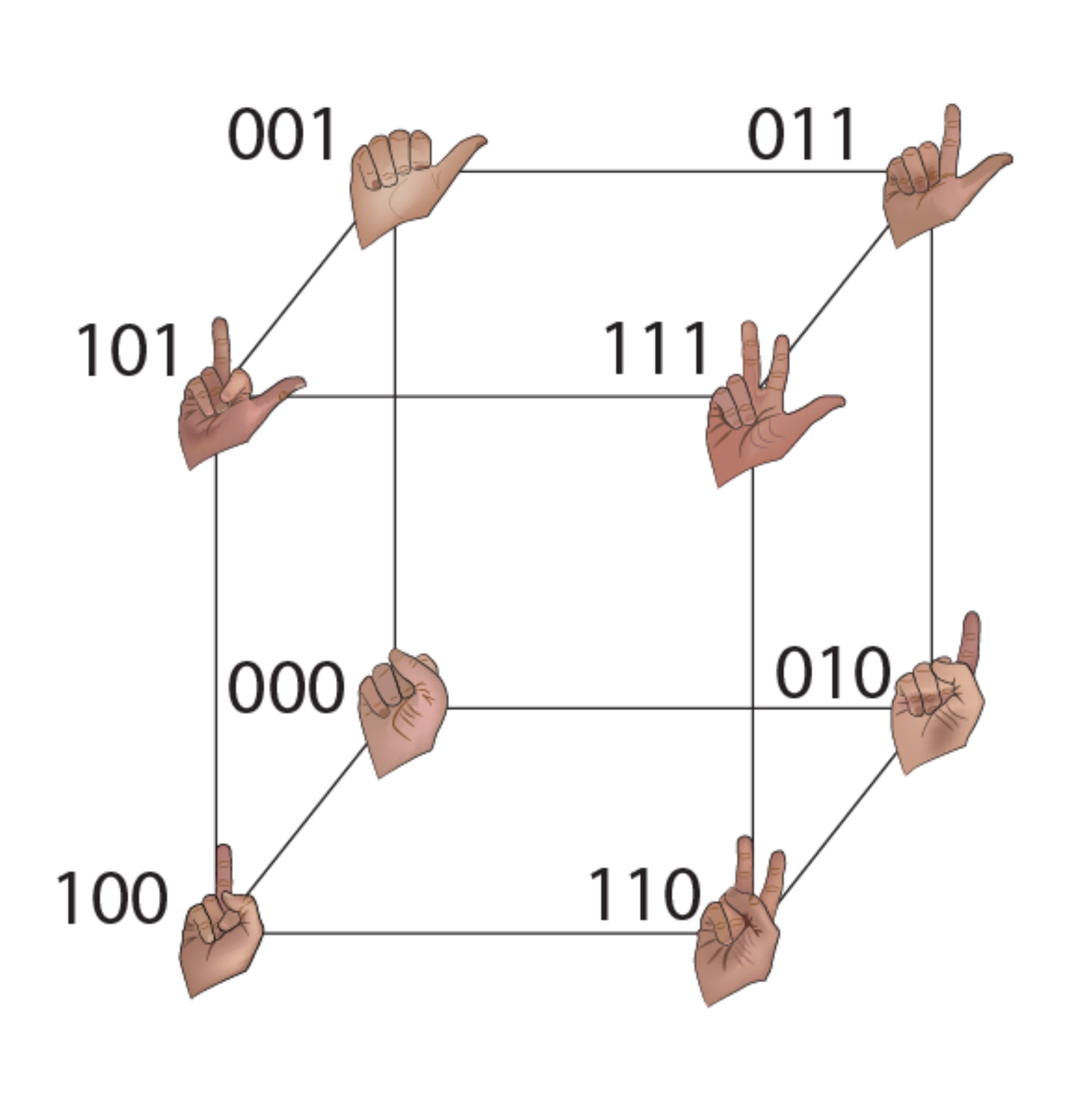}
\end{center}
\caption{The vertices of a cube correspond to binary digits}
\label{digitalcube}
\end{figure}

 In Fig.~\ref{digitalcube}, we have arranged the first eight of these configurations at the vertices of a cube in $3$-dimensional space. 
 
 In returning to the description of the unit $n$-dimensional cube as $[0,1]^n$ --- the set of functions from $\{1, \ldots , n\}$ into the unit interval, we observe that the binary sequences are the extreme values of such functions: $x_j=0$ or $x_j=1$ for each $j=1, \ldots, n.$
 
 \subsection{Kitchen drawers}
 
 Imagine a small kitchen that has exactly four drawers. Each drawer can be open or closed. As such, we can describe the state of the kitchen cabinet as a binary sequence. But let us suppose that each drawer can be open any distance from $0$ cubits up to fully open at a cubit.{\footnote {By convention a cubit is roughly half a meter. The units are chosen to appear to be realistic from the point of view of a kitchen.}} Let us label the drawers with integers $1$ through $4$. Then the configurations of the set of drawers corresponds to a function from $\{1,2,3,4\}$ into $[0,1]$. So the possible states for this kitchen cabinet corresponds to the points of a $4$-dimensional cube. We call this cube the {\it hypercube}. 
 
 \subsection{An Archive}
 
 It is easy to generalize from the kitchen cabinet to a large archive in a bureaucracy. This archive has a large number of cabinets, say $n$ of them. Each cabinet consists of a single drawer. The cabinets lie along the south wall of a building and are numbered east-to-west from $1$ through $n$. The drawers can open up to a full meter. And so the set of possible configurations of the cabinets corresponds to the unit $n$-dimensional cube $[0,1]^n$ since if $x_j$ is the distance (in meters) that the $j\/$th drawer is open, then $0\le x_j \le 1$. The configuration of any one drawer is independent of the configuration of all the other drawers. 
 
 \subsection{Triangles, tetrahedra, and so forth}
 
   \begin{figure}[htb]
\begin{center}
\includegraphics[scale=.4]{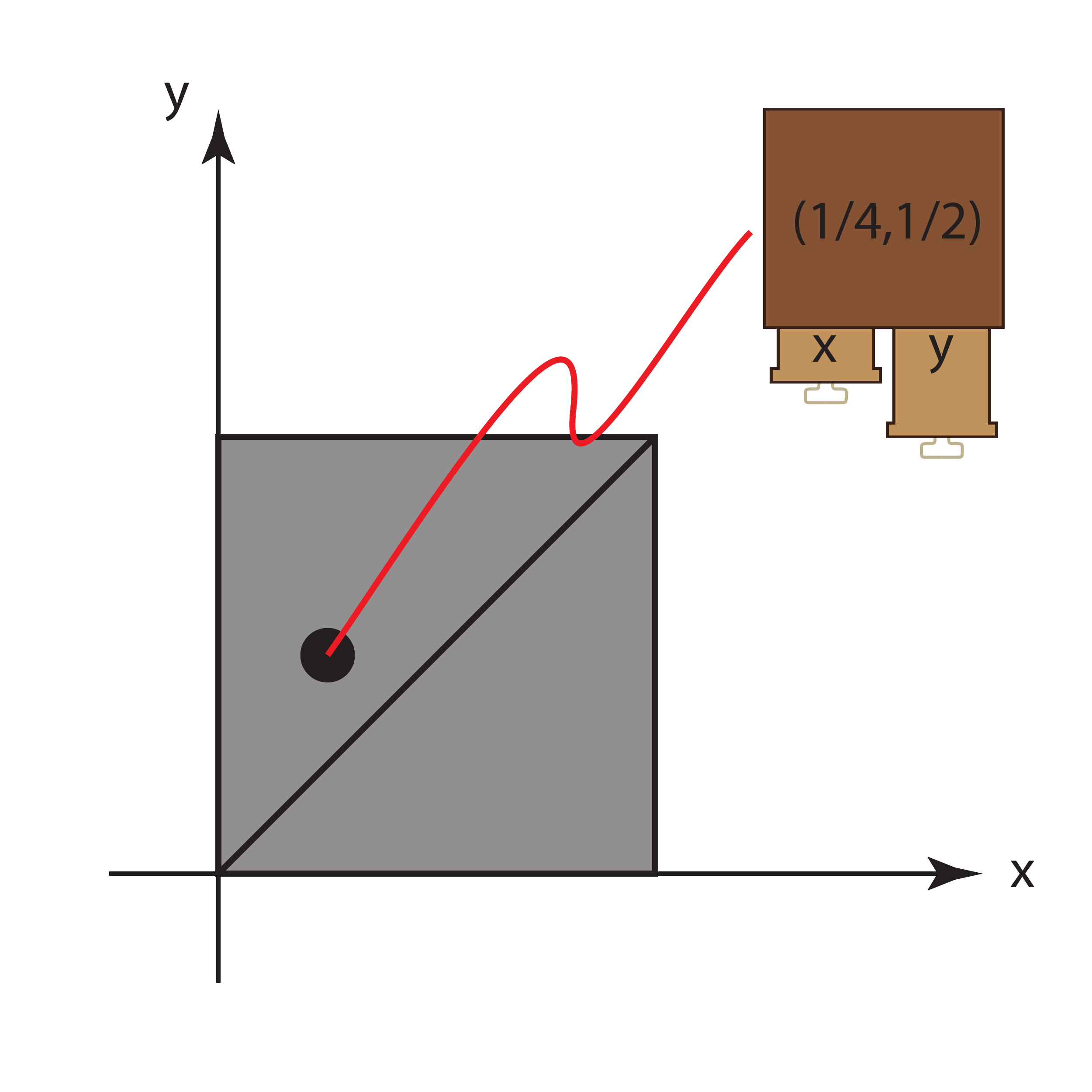}
\end{center}
\caption{The square as a pair of open drawers}
\label{drawersin2D}

\end{figure}
 
 We begin this subsection with a simple case. Consider the possible configurations of a pair of drawers. Let $x$ denote the distance that the first drawer is opened, and let $y$ denote the distance that the second drawer is opened. Of course $0\le x,y  \le 1$ --- that is each drawer can be opened any amount from $0$ meters to $1$ meter. This situation can be subdivided into two cases. Either the drawer on the left is no more open than the drawer on the right is, or vice versa. Thus either $0\le x \le y \le 1$ or $0\le y \le x \le 1.$ These cases are not disjoint; their coincidence occurs when both drawers are equally open: $x=y.$ The set of all configurations of the drawers correspond to a unit square. We have decomposed the square as the union of two isosceles right triangles whose legs are a unit length. They intersect along their diagonals. 
 
 For example, the point $(1/4,1/2)$ in the unit square corresponds to the configuration in which the drawer on the left is open half as much as the drawer on the right which is half-open. See Fig.~\ref{drawersin2D}.

   \begin{figure}[htb]
\begin{center}
\includegraphics[scale=.5]{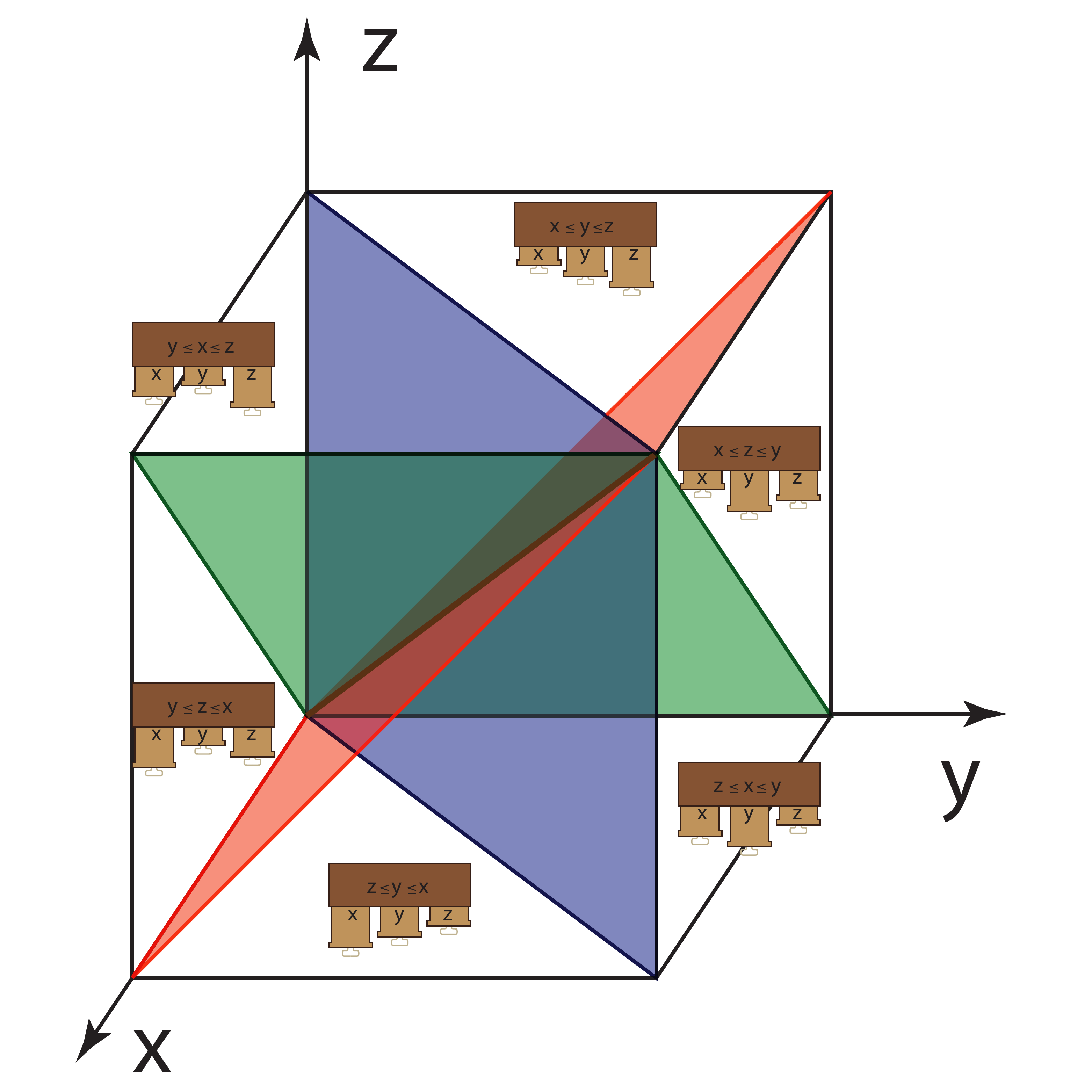}
\end{center}
\caption{Decomposing the 3-cube as the union of six tetrahedra}
\label{drawercube}
\end{figure}
 
 Now consider three drawers. Let $x$, $y$, and $z$ denote the amount that each is pulled out. As before $0\le x,y,z\le 1$, but there is no relation among the openness of the drawers in general. We decompose this situation into six configurations: 
 $$0 \le x \le y \le z \le 1, \quad 0 \le y \le x \le z \le 1,
  \quad 0 \le x \le z \le y \le 1,$$ $$ \quad 0 \le z \le x \le y \le 1, \quad
  0 \le y \le z \le x \le 1, \quad \& \quad 0 \le z \le y \le x \le 1.$$
as indicated in Fig.~\ref{drawercube}. Using the coordinates of the $3$ dimensional cube that are indicated, the cube has been decomposed into six congruent tetrahedra. Each tetrahedron has three edges with length 1, two edges of length $\sqrt{2}$ and one edge of length $\sqrt{3}$. All the tetrahedra intersect along the edge of length $\sqrt{3}$ since this line corresponds to the configuration of drawers in which all are equally open. The interior triangular faces with edge lengths $1$, $\sqrt{2}$ and $\sqrt{3}$ correspond to configurations of drawers in which at least two are equally open and the other drawer is either open no more or no less than the two equal drawers. The exterior triangular faces that have edge lengths $1$, $1$ and $\sqrt{2}$ correspond to exactly one drawer being fully open or fully closed.

In full generality, we have $n$ drawers with drawer number $j$ open to a distance $x_j$. In this case $0\le x_1, x_2, \ldots, x_n \le 1$. There are $n!=n(n-1)!$ (where $0!=1$) possible orderings of the distances that the drawers can be open. To see this observe that the $j\/$th drawer may be pulled out a greater distance than all the others. Thus there are $n$ possible ways for one drawer to be more extended than all the rest. By induction, the remaining $(n-1)$ drawers are ordered in $(n-1)!$ configurations. 

A {\it right $n$-dimensional simplex} consists of the point set $$\Delta^n = \{ (x_1, x_2, \ldots, x_n) \in [0,1]^n: 0 \le x_1 \le x_2 \le \cdots \le x_n \le 1 \}.$$
The $n$-dimensional cube can be decomposed into $n!$ such simplices that are all congruent. So the $n$-dimensional volume of such a simplex is $\frac{1}{n!}$. 

There is a typographical trick called {\it the $L7$-trick}\/{\footnote{Cf. the second verse to {\it Woolly Bully} as sung by Sam the Sham and the Pharohs}}, to obtain the set of vertices for each simplex. Let $n$ be fixed. Consider a row of $n$ $0\/$s (separated, if you like, by commas): $00\ldots0$. Also construct a column of $n$ $1\/$s. Now the column of $1\/$s can be appended below any one of the $0\/$s. Thus there are $n$ possible ways of doing so. For each such shape, we can insert  the $(n-1)!$ shapes that exist by induction. 
The first step of the induction is to consider the array $\left[ \begin{array}{c} 0 \\ 1 \end{array} \right]$. In the next step, we have the arrays
$$\begin{array}{cc} \left[ \begin{array}{cc} 0& 0 \\ & 1 \\ & 1 \end{array}  \right] & \left[ \begin{array}{cc} 0&0 \\ 1& \\ 1& \end{array} \right] \end{array}$$ which can be filled as

$$\begin{array}{cc} \left[ \begin{array}{cc} 0& 0 \\ 0 & 1 \\  1& 1 \end{array}  \right] & \left[ \begin{array}{cc} 0&0 \\ 1& 0\\ 1&1 \end{array} \right] \end{array}.$$
The rows are the coordinates of the vertices of the two triangles $0\le x \le y \le 1$ and $0 \le y \le x \le 1$ into which the unit square was decomposed. 
 For the cube, we start from the three arrays:
 $$\begin{array}{ccc} 
 \left[\begin{array}{ccc} 0 & 0 & 0 \\ & & 1 \\ & & 1 \\ & & 1 \end{array} \right] &    \left[\begin{array}{ccc} 0 & 0 & 0 \\ & 1&  \\ &1 &  \\ & 1&  \end{array} \right] &   \left[\begin{array}{ccc} 0 & 0 & 0 \\ 1 & &  \\  1 & &  \\ 1  & &  \end{array} \right]  \end{array} $$
 and fill them respectively with the two arrays above.

$$\begin{array}{ccc} 
 \left[\begin{array}{ccc} 0 & 0 & 0 \\ 0 & 0 & 1 \\ 0& 1 & 1 \\ 1 &1 & 1 \end{array} \right] &    \left[\begin{array}{ccc} 0 & 0 & 0 \\ 0& 1&0  \\0 &1 & 1 \\ 1& 1&1  \end{array} \right] &   \left[\begin{array}{ccc} 0 & 0 & 0 \\ 1 &0 & 0 \\  1 & 0&1  \\ 1  &1 & 1 \end{array} \right]  \end{array} $$

 $$\begin{array}{ccc} 
  \left[\begin{array}{ccc} 0 & 0 & 0 \\ 0&0 & 1 \\1 &0 & 1 \\ 1&1 & 1 \end{array} \right] &    \left[\begin{array}{ccc} 0 & 0 & 0 \\ 0& 1&0  \\ 1&1 &0  \\ 1& 1&1  \end{array} \right] &   \left[\begin{array}{ccc} 0 & 0 & 0 \\ 1 & 0&0  \\  1 &1 &0  \\ 1  &1 &1  \end{array} \right]  \end{array}. $$
 In the next step, each of these six matrices are fed into the four patterns 
  $$\begin{array}{cccc}  
  \left[\begin{array}{cccc} 0 & 0 & 0 &0 \\& &  & 1 \\ & & & 1 \\ & &  & 1 \end{array} \right] &
 \left[\begin{array}{cccc} 0 & 0 & 0 &0 \\& & 1 &  \\ & &1 &  \\ & & 1 &  \end{array} \right] &
  \left[\begin{array}{cccc} 0 & 0 & 0 &0 \\& 1&  &  \\ & 1 & &  \\ & 1 &  &  \end{array} \right] &
   \left[\begin{array}{cccc} 0 & 0 & 0 &0 \\1& &  &  \\ 1& & &  \\ 1& &  &  \end{array} \right] 
    \end{array}. $$

\begin{figure}
\includegraphics[scale=.07]{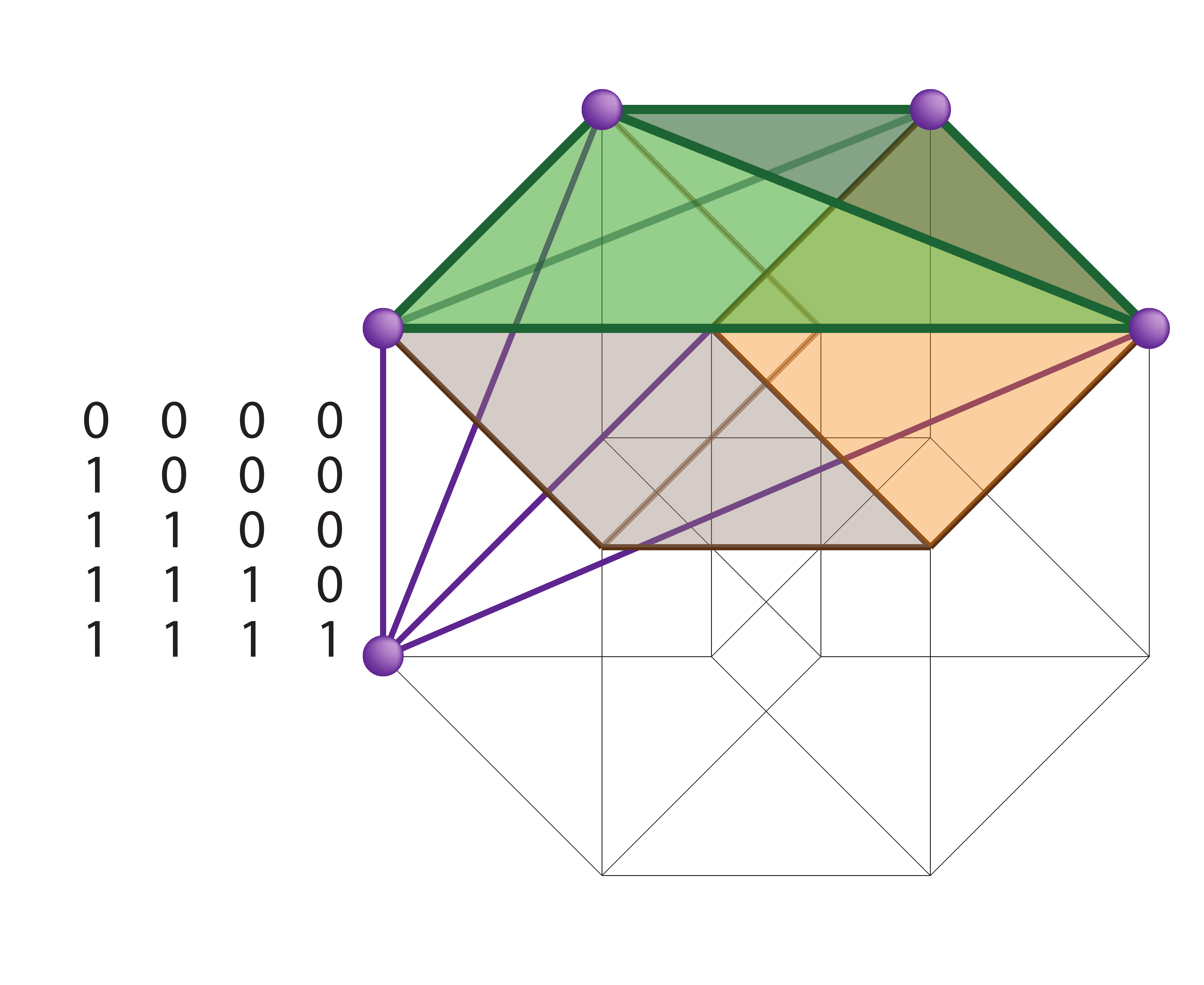}
\includegraphics[scale=.07]{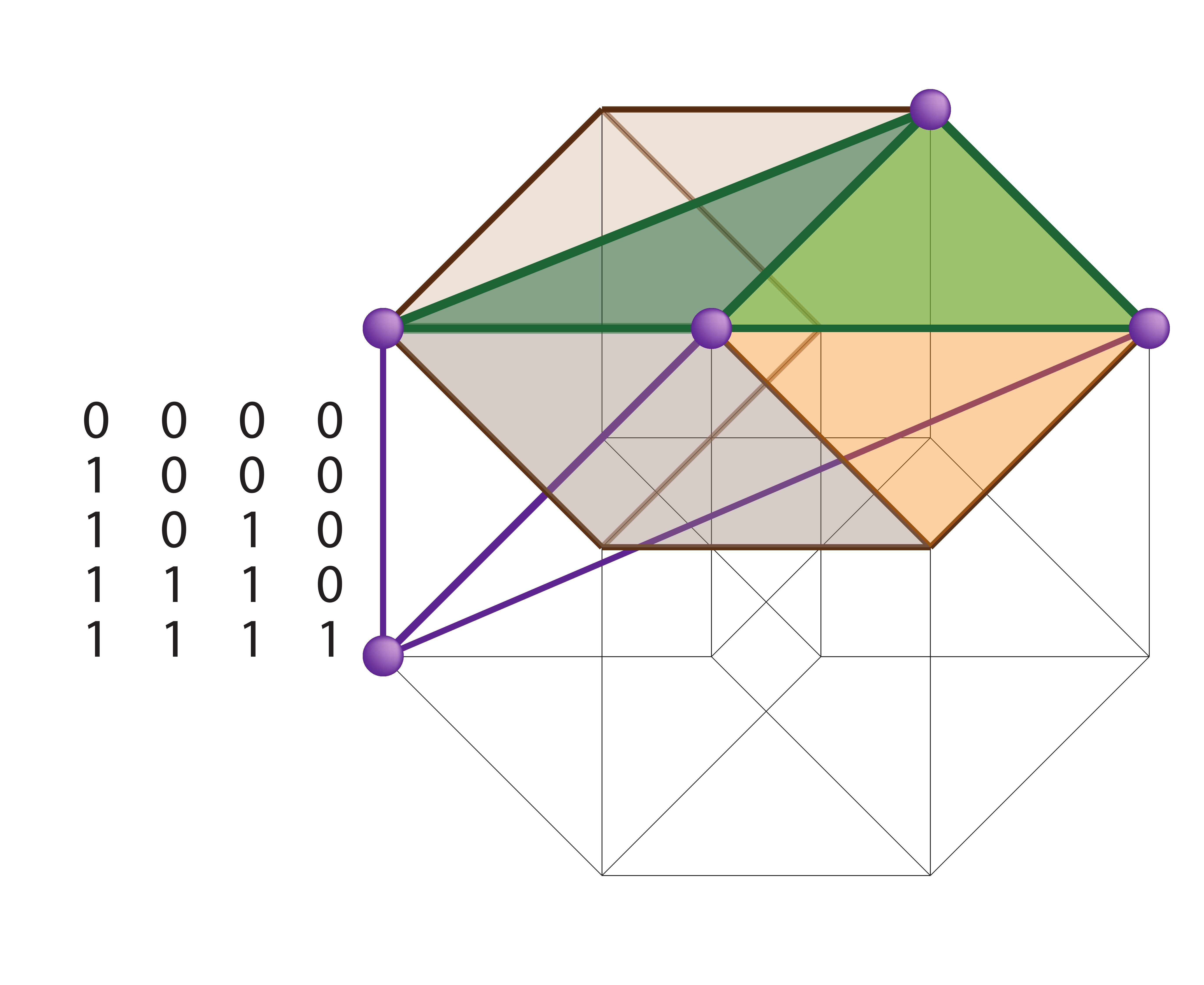}
\includegraphics[scale=.07]{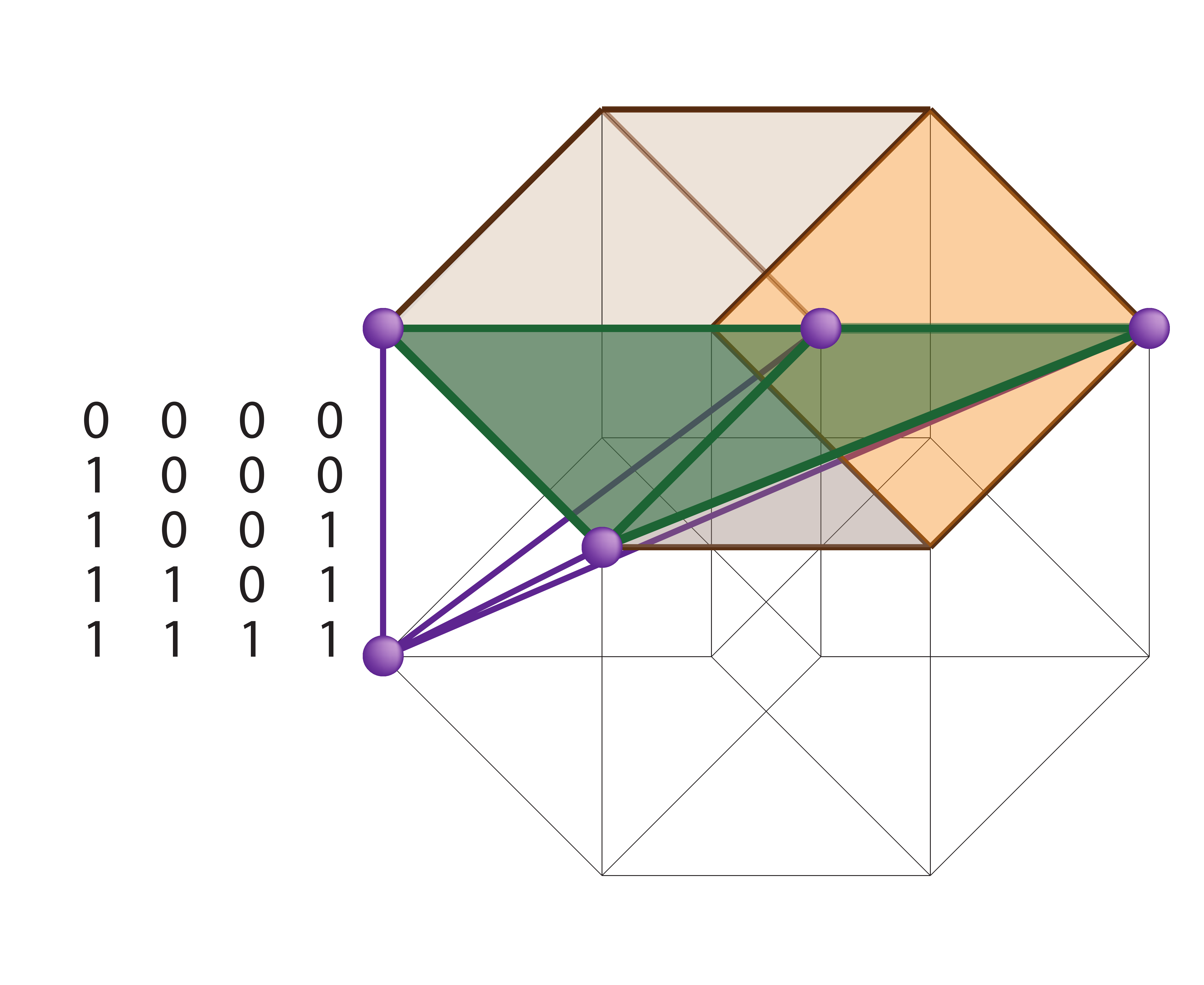}
\includegraphics[scale=.07]{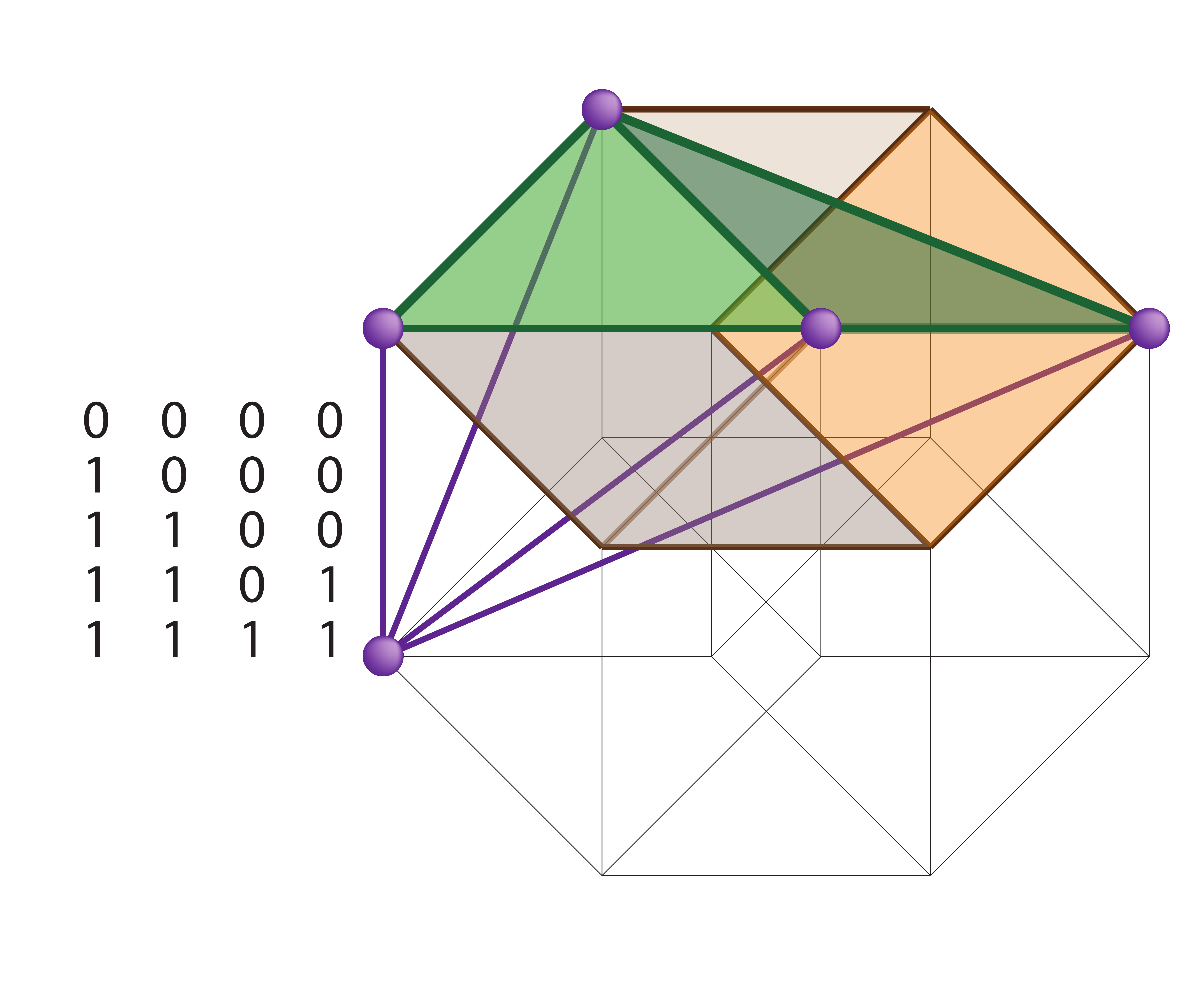}
\includegraphics[scale=.07]{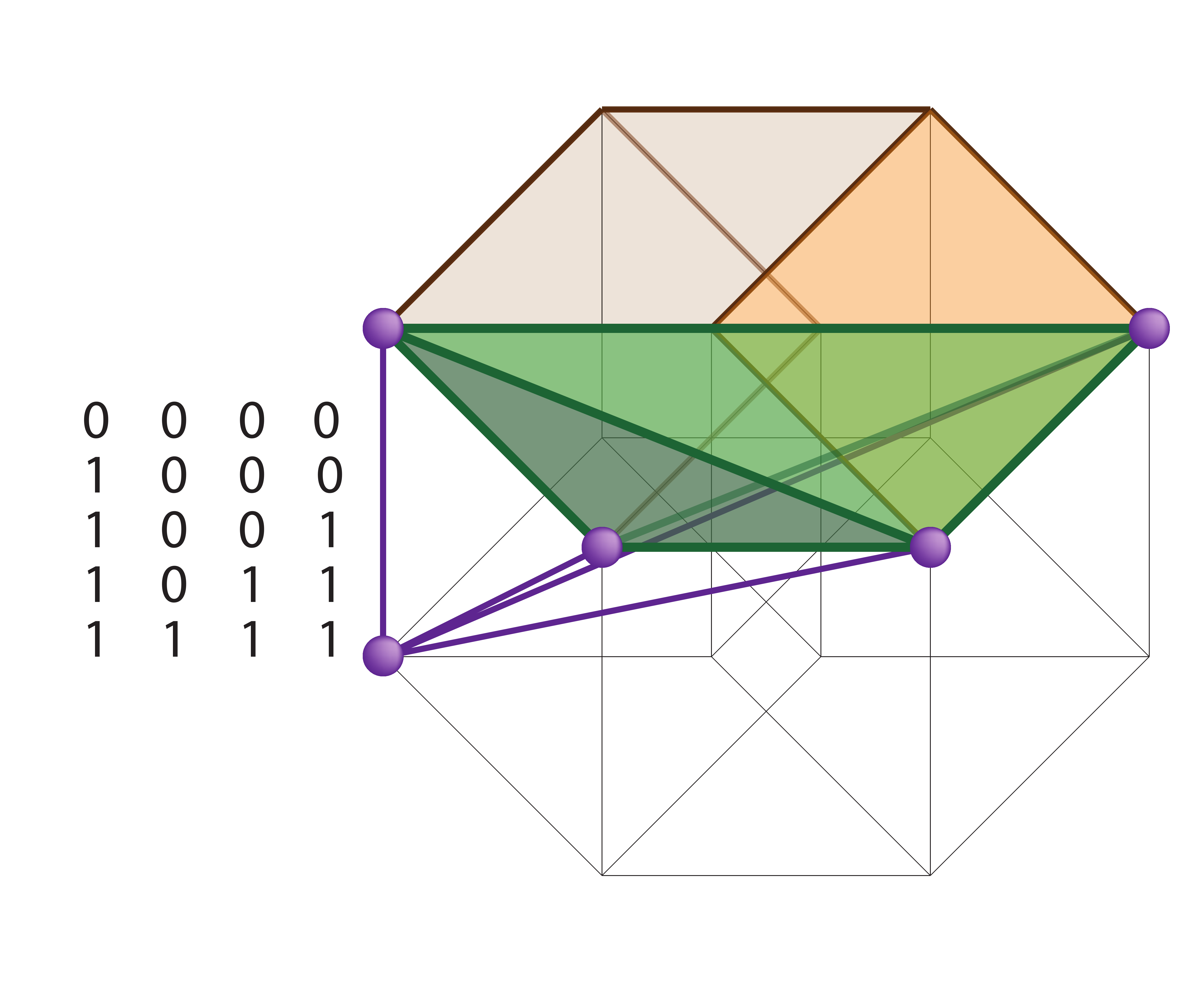}
\includegraphics[scale=.07]{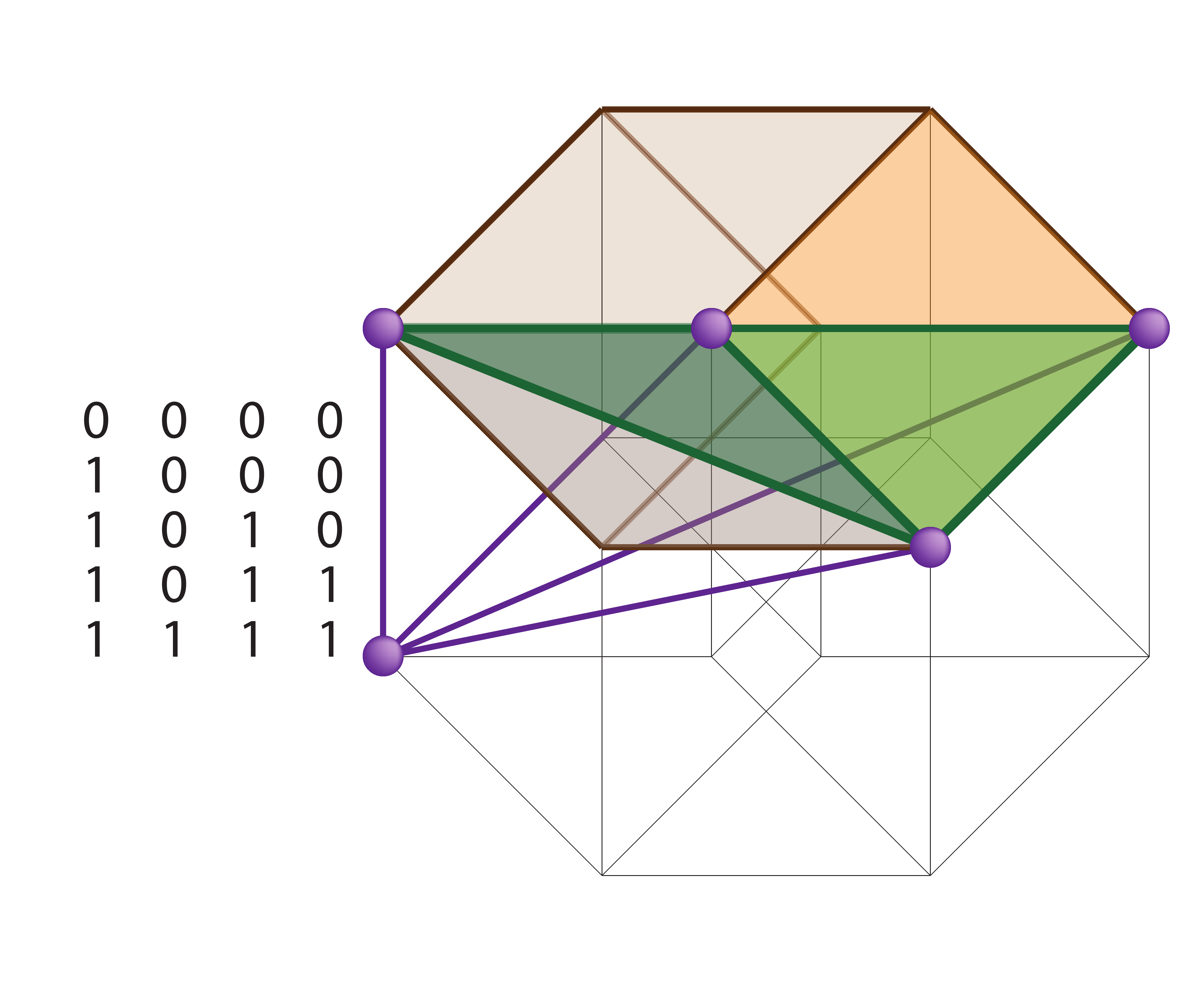}
\includegraphics[scale=.07]{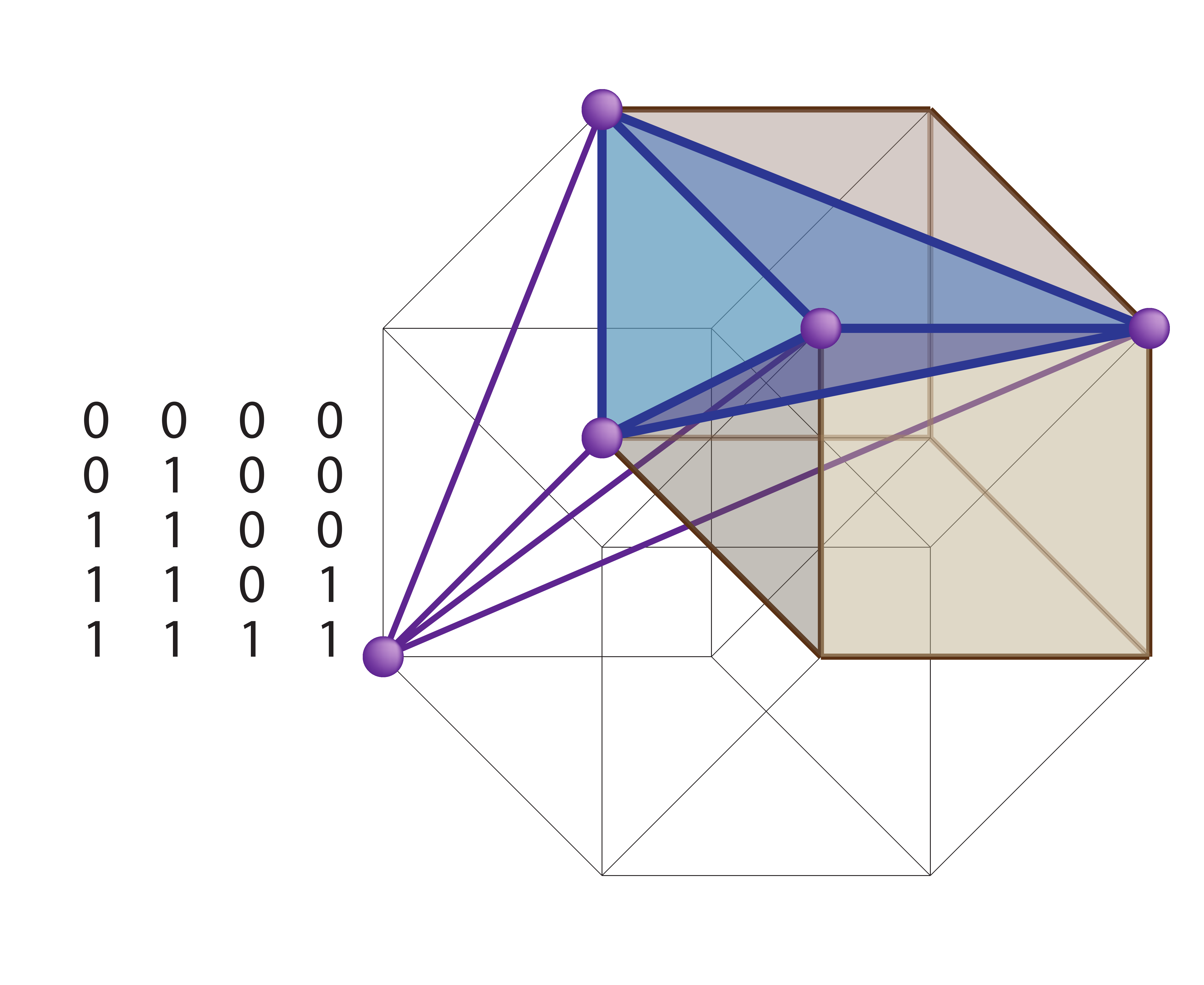}
\includegraphics[scale=.07]{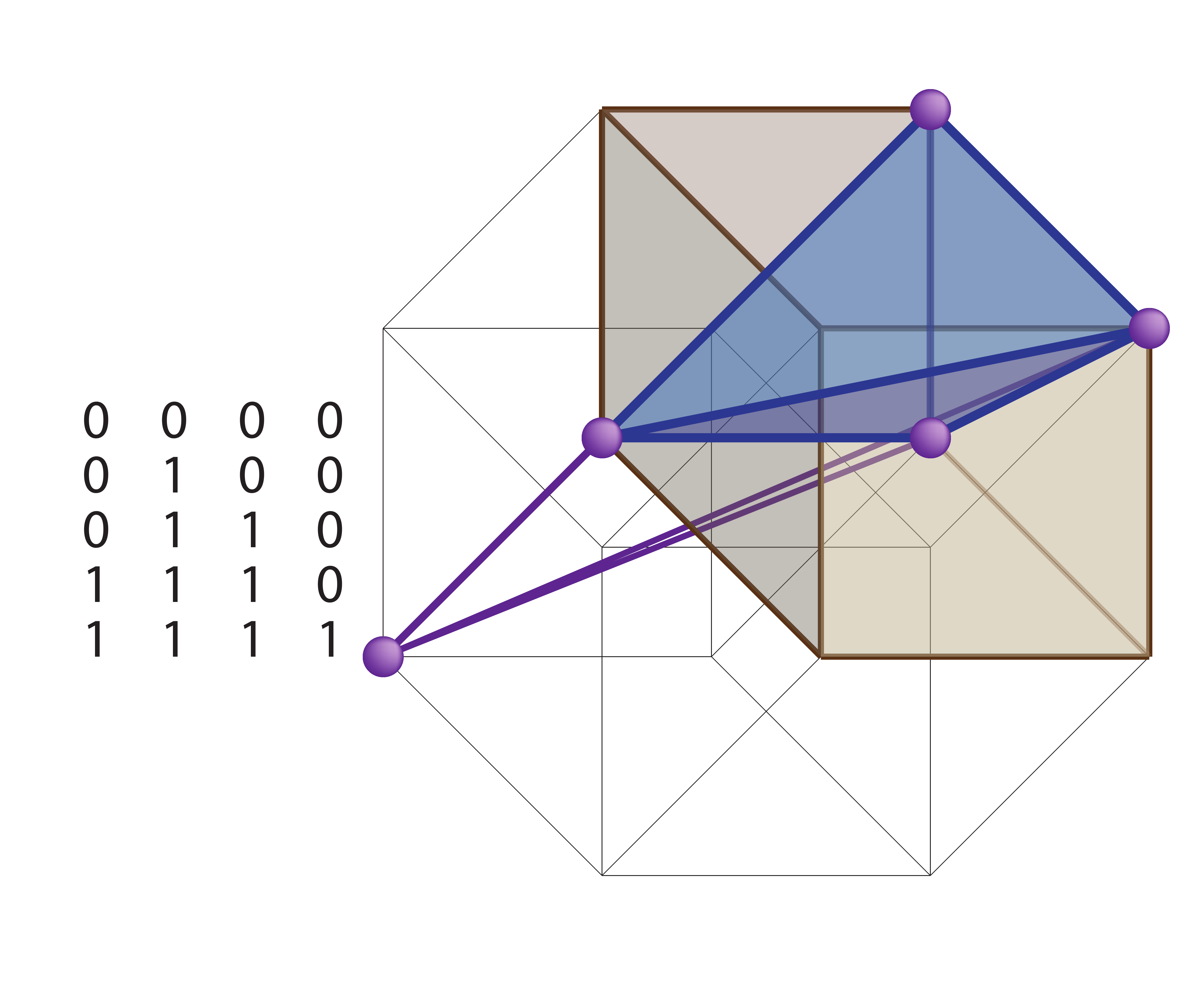}
\includegraphics[scale=.07]{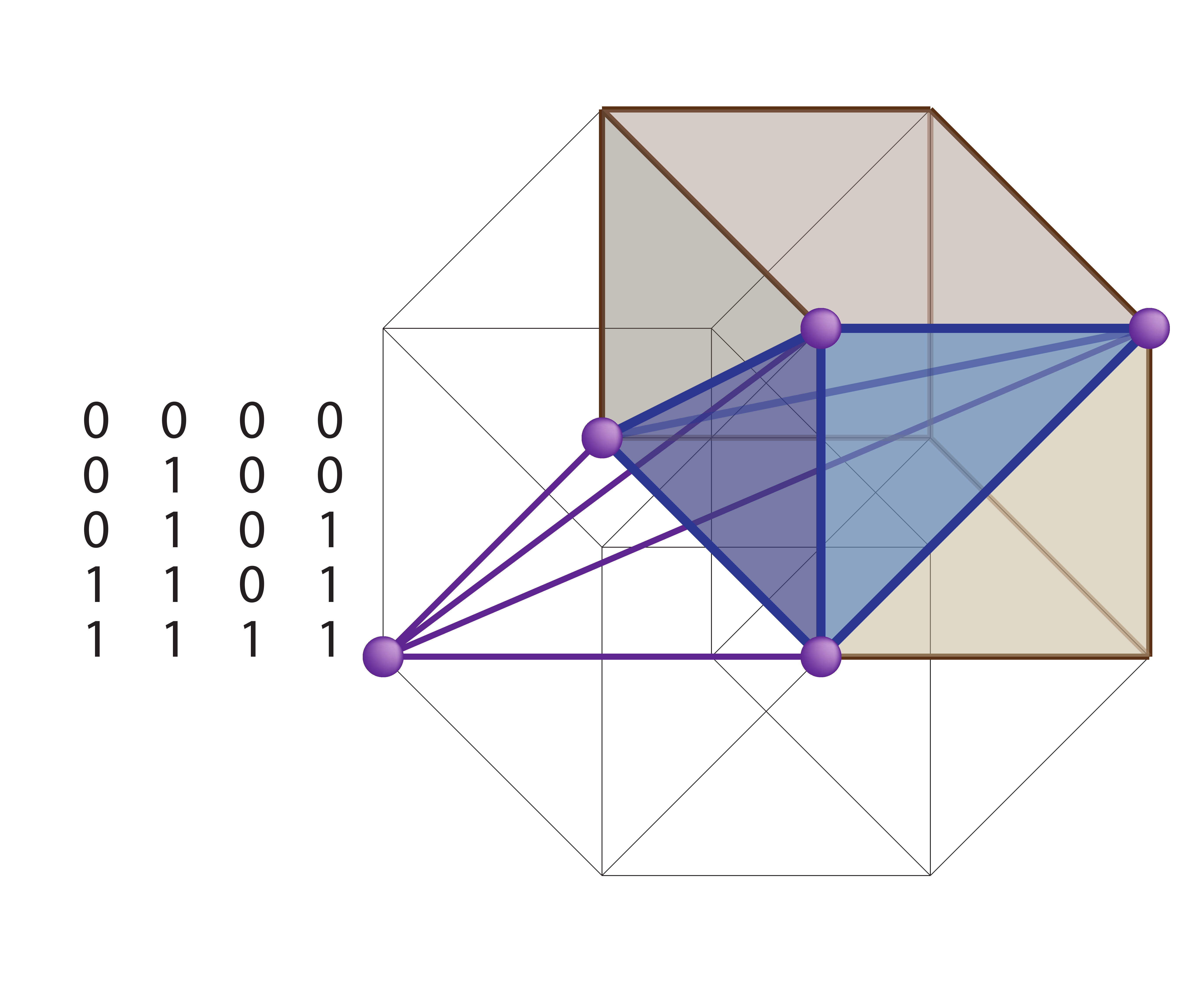}
\includegraphics[scale=.07]{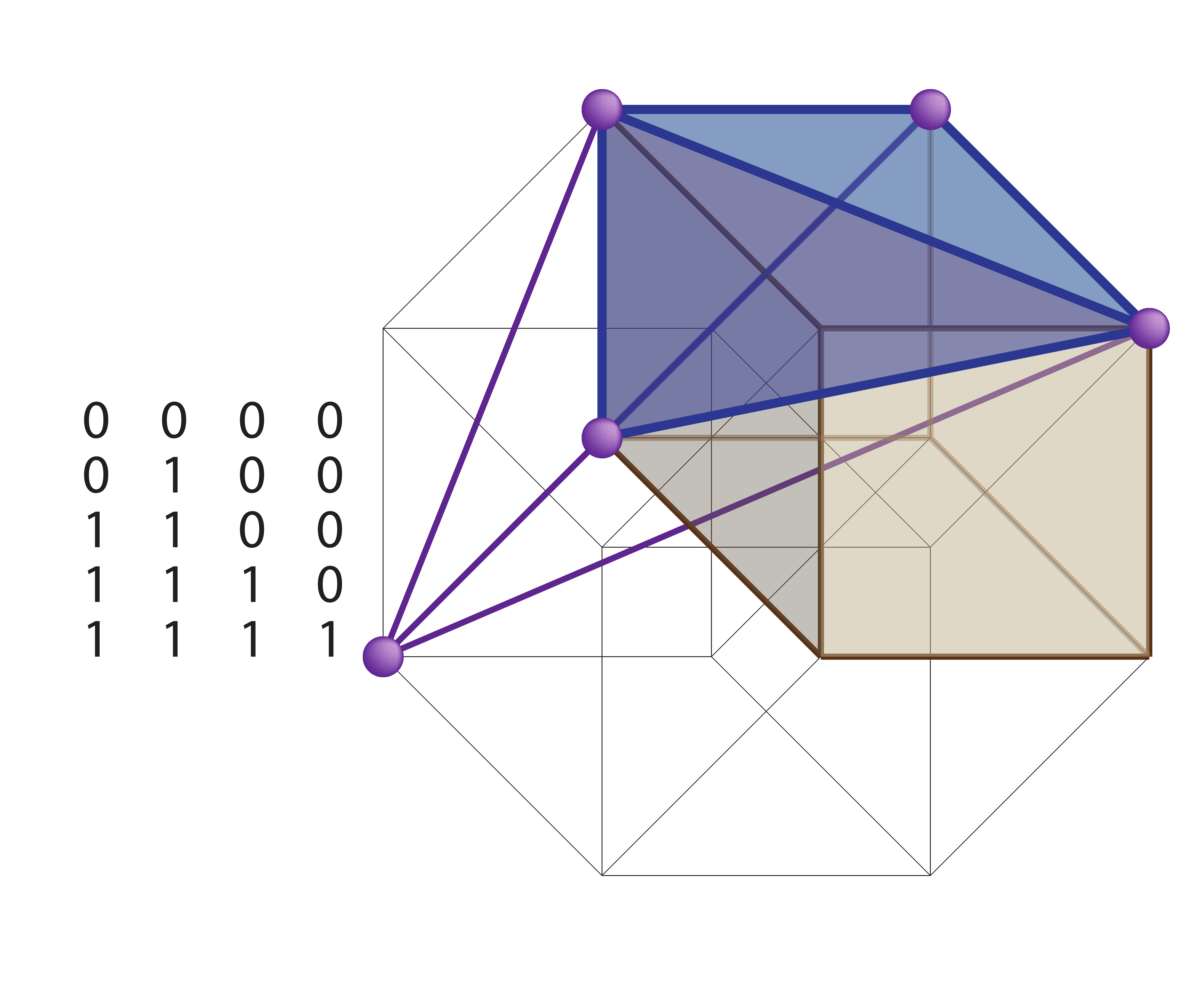}
\includegraphics[scale=.07]{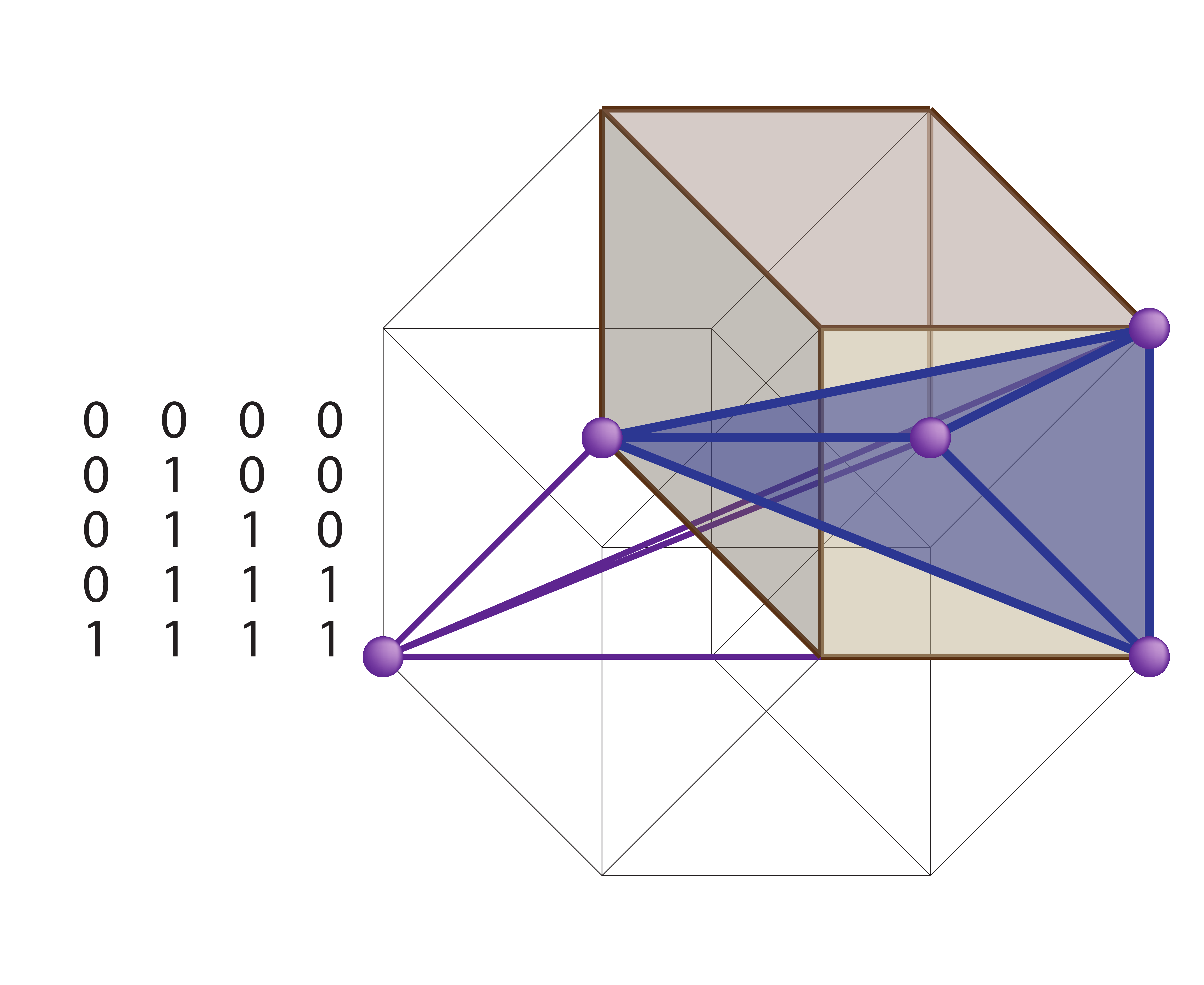}
\includegraphics[scale=.07]{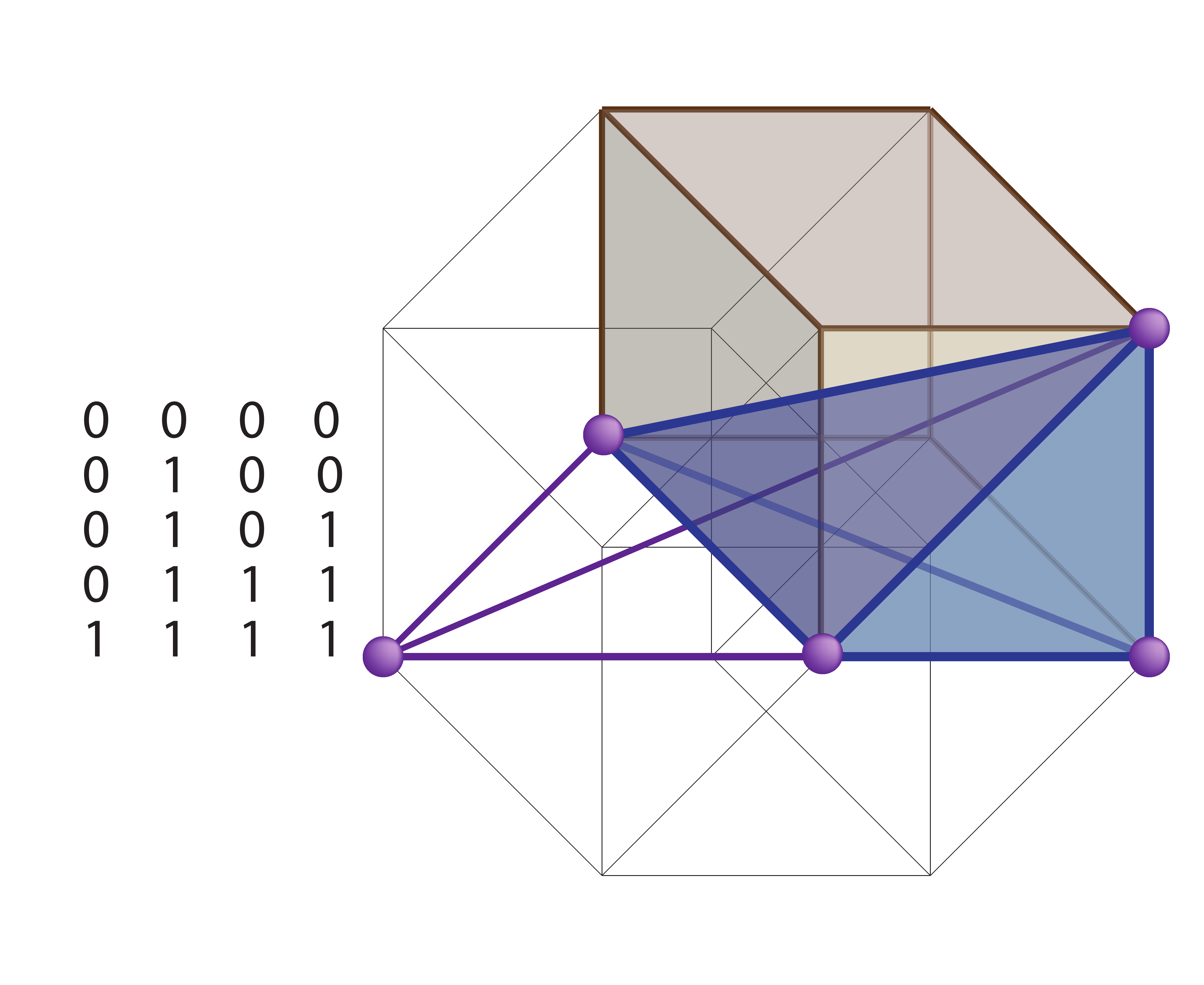}
\includegraphics[scale=.07]{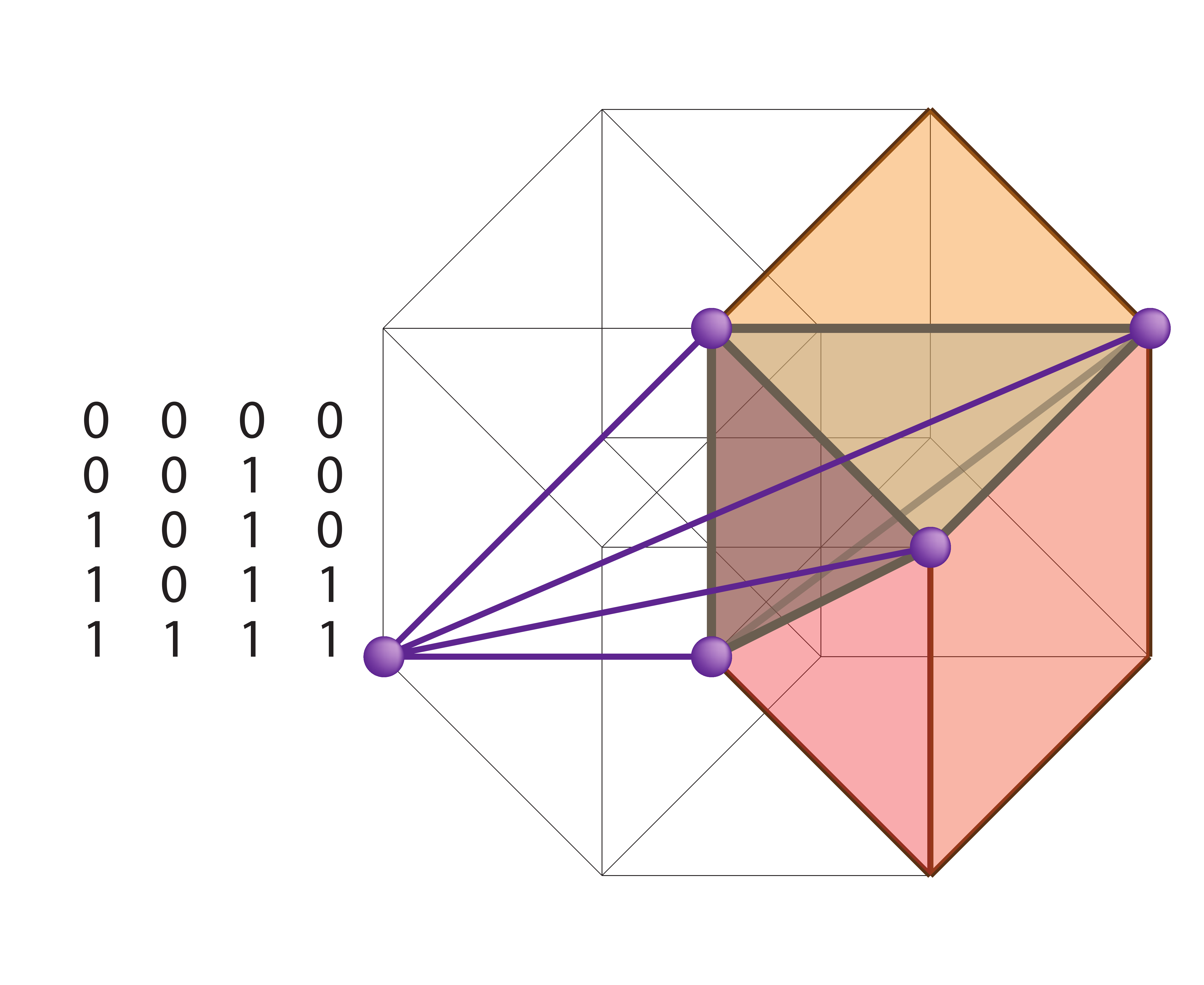}
\includegraphics[scale=.07]{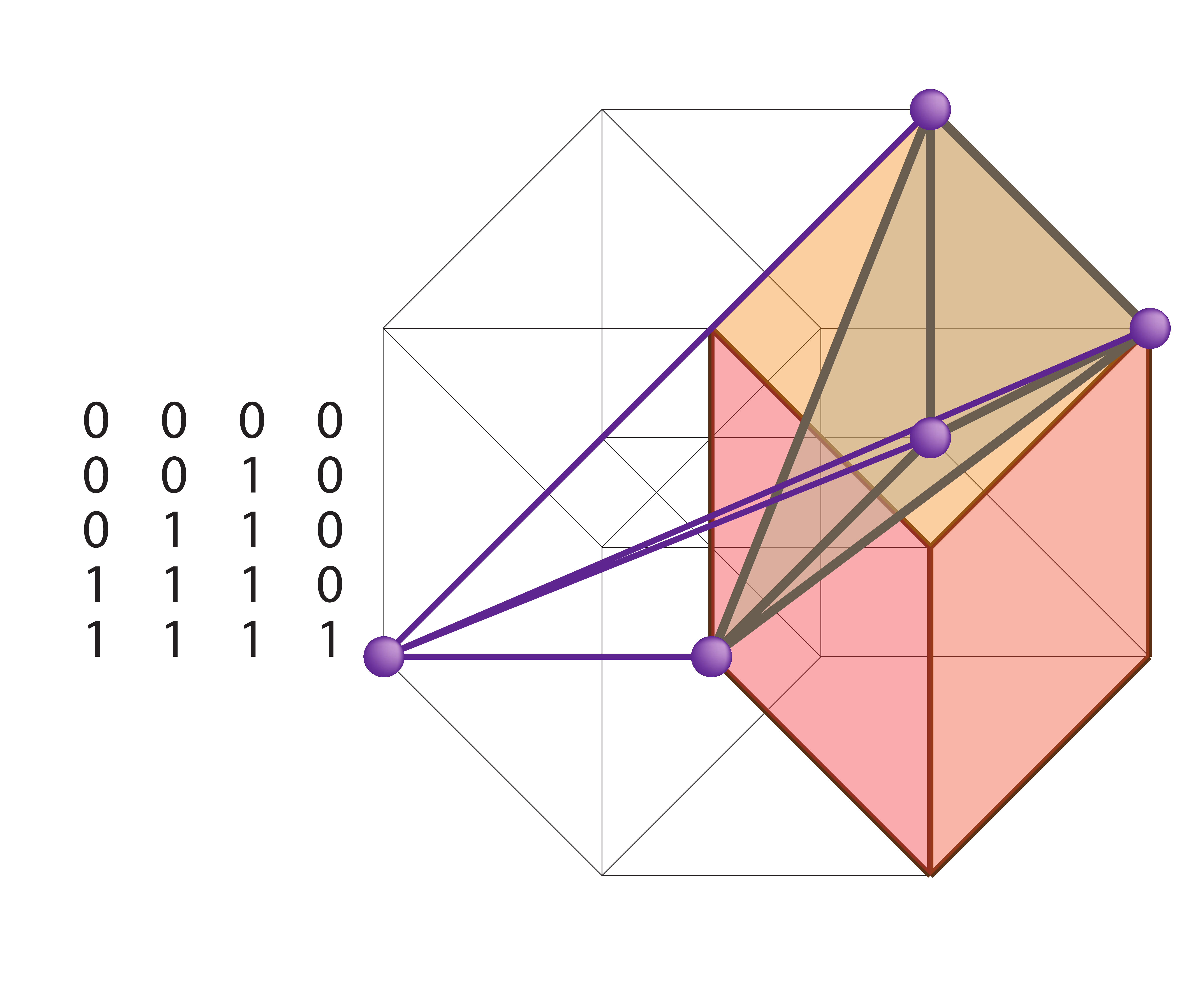}
\includegraphics[scale=.07]{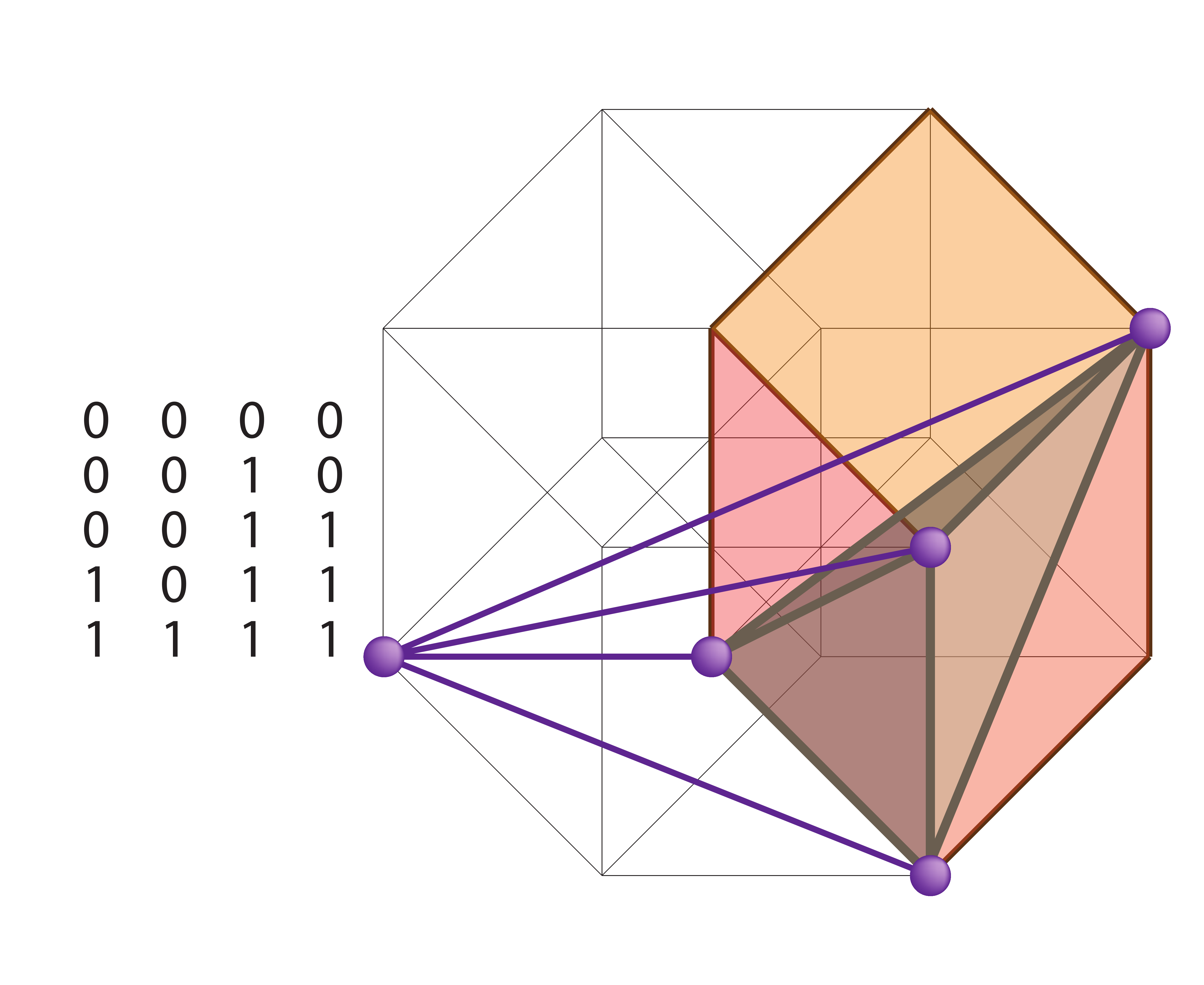}
\includegraphics[scale=.07]{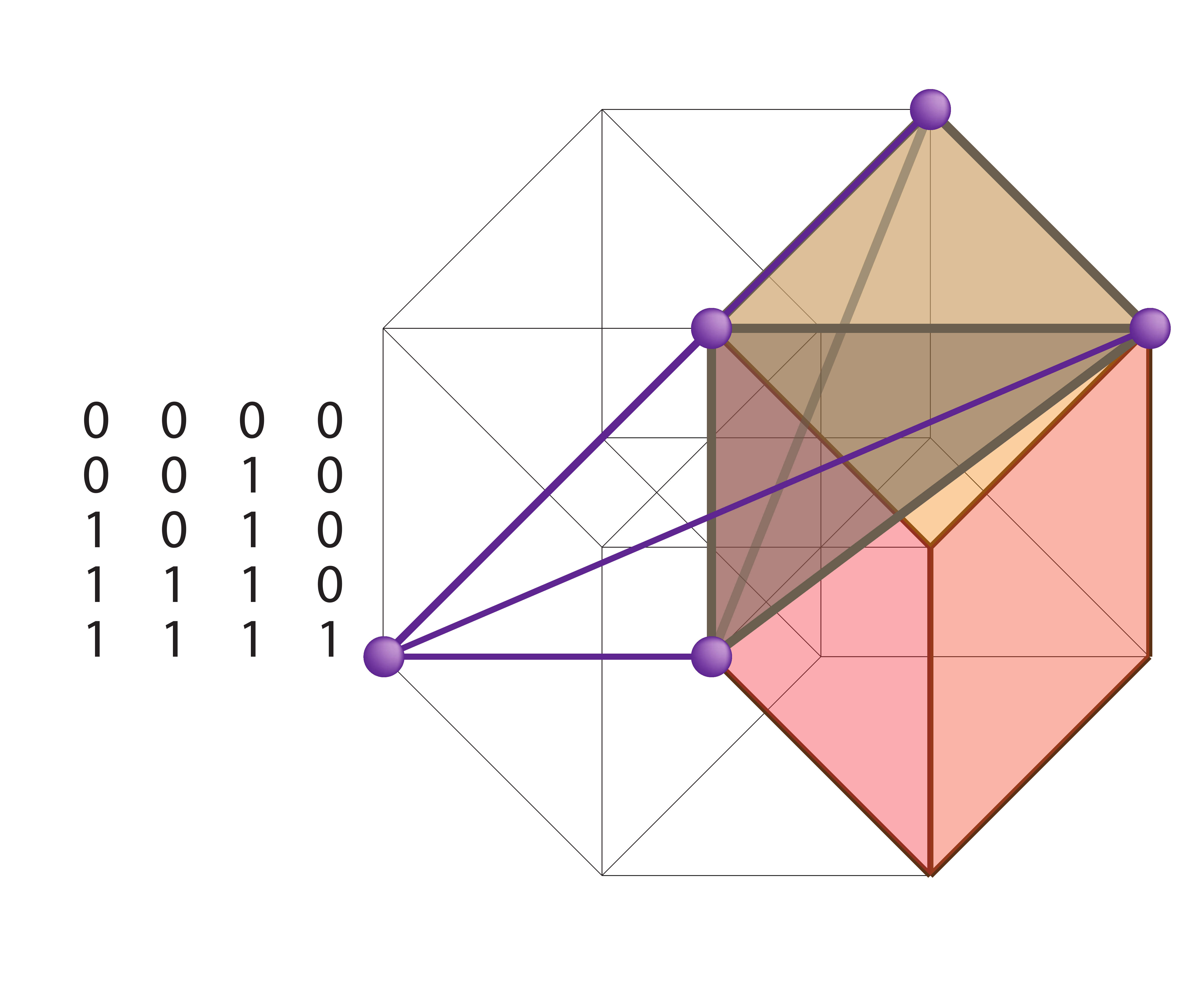}
\includegraphics[scale=.07]{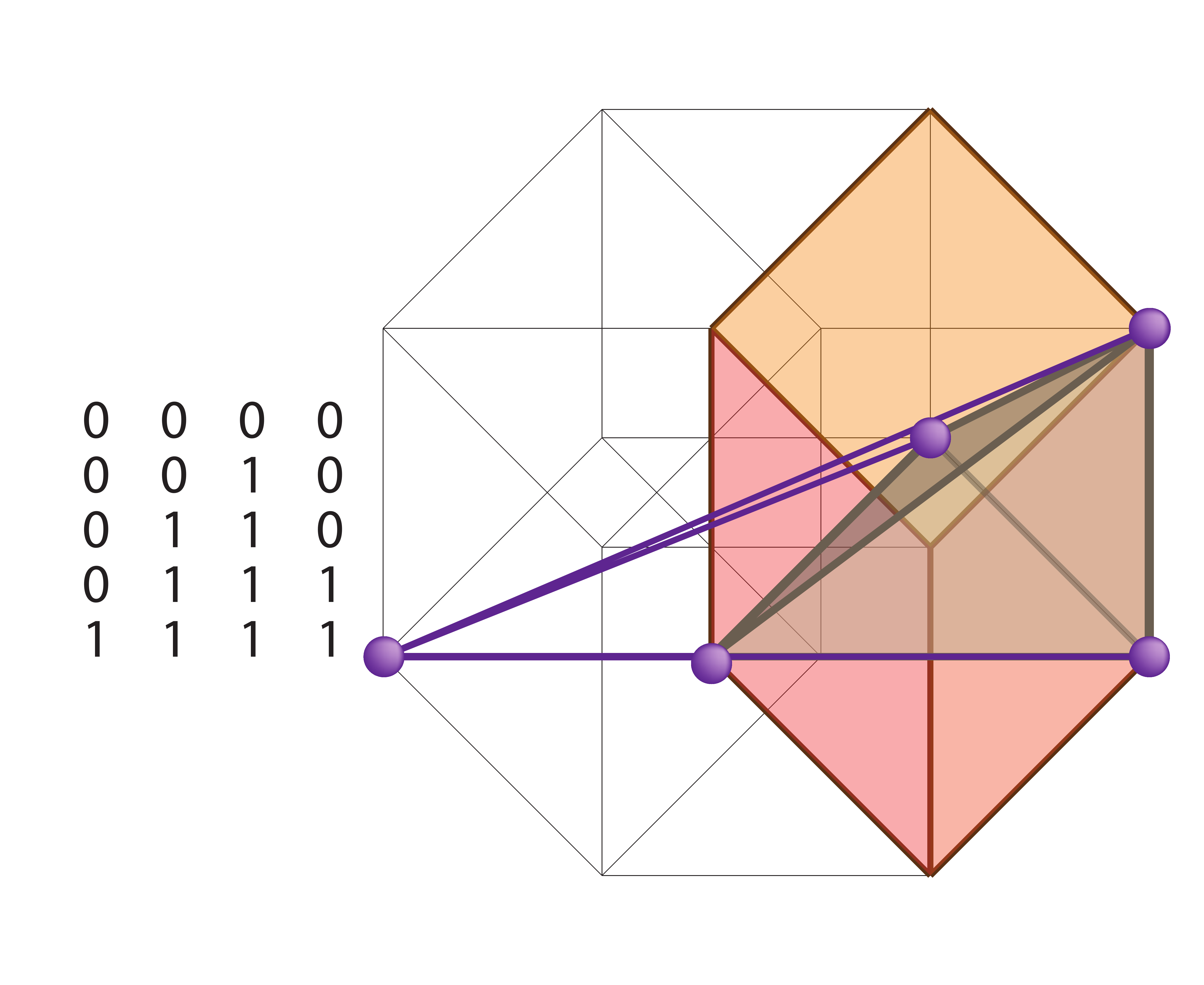}
\includegraphics[scale=.07]{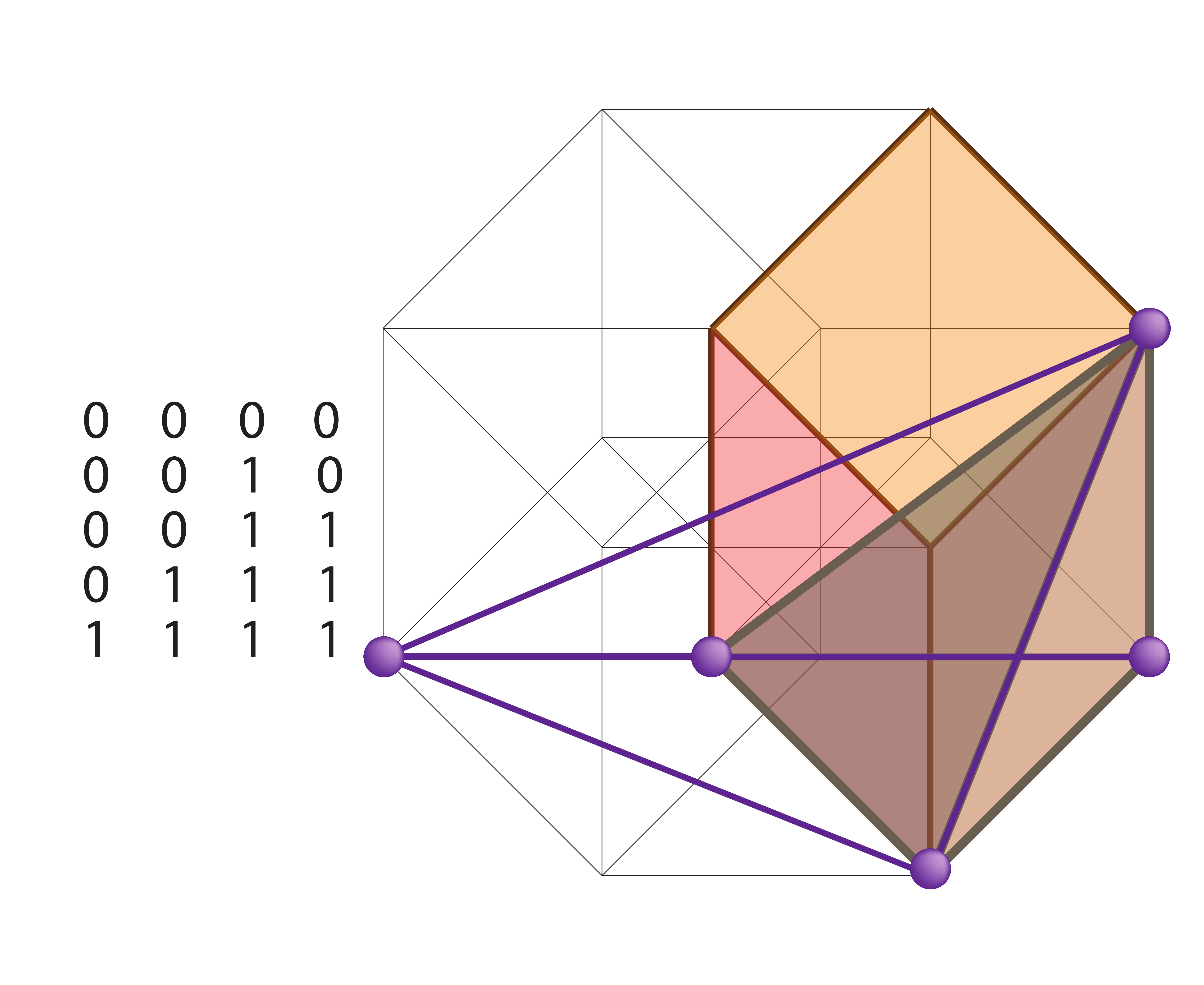}
\includegraphics[scale=.07]{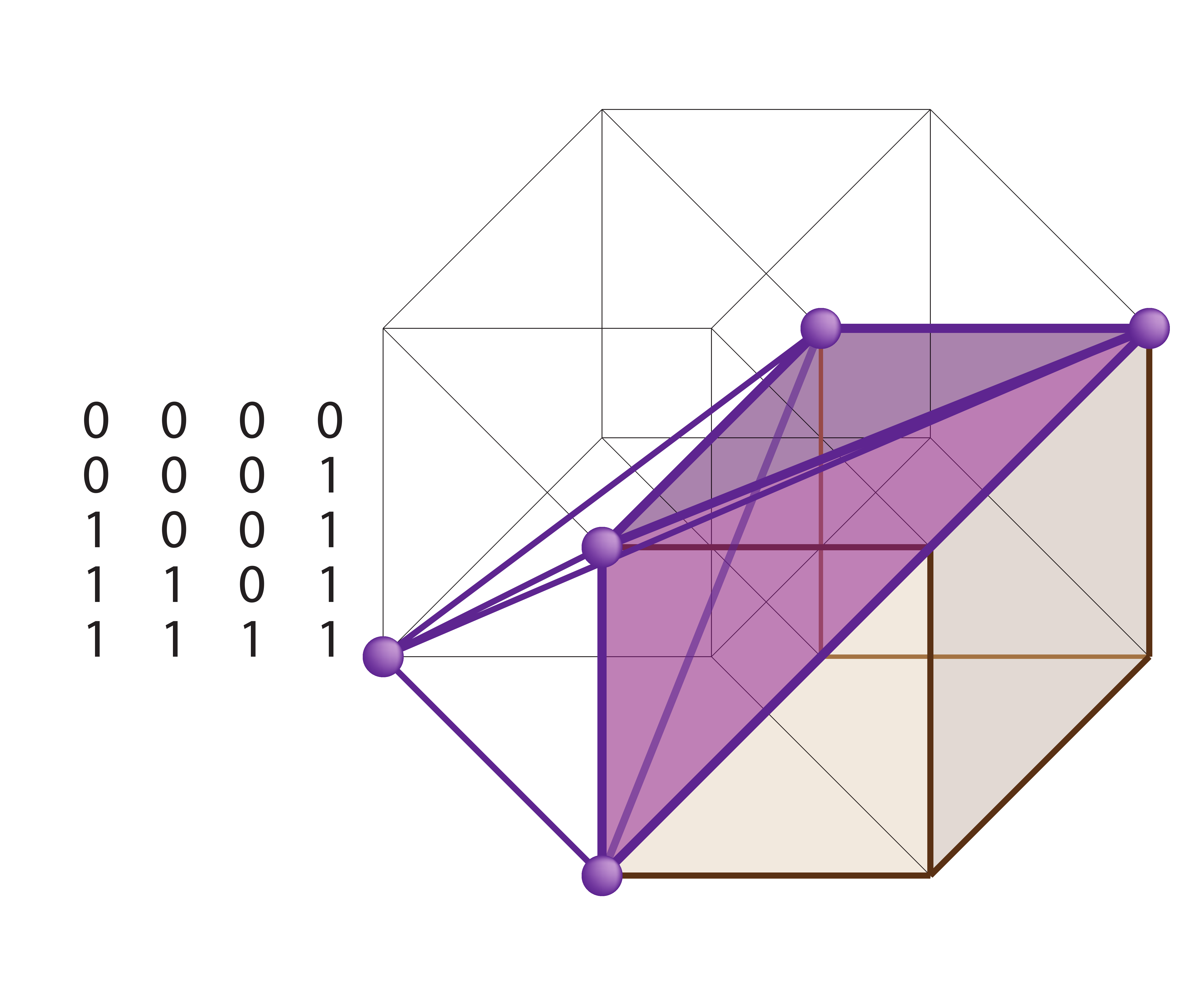} 
\includegraphics[scale=.07]{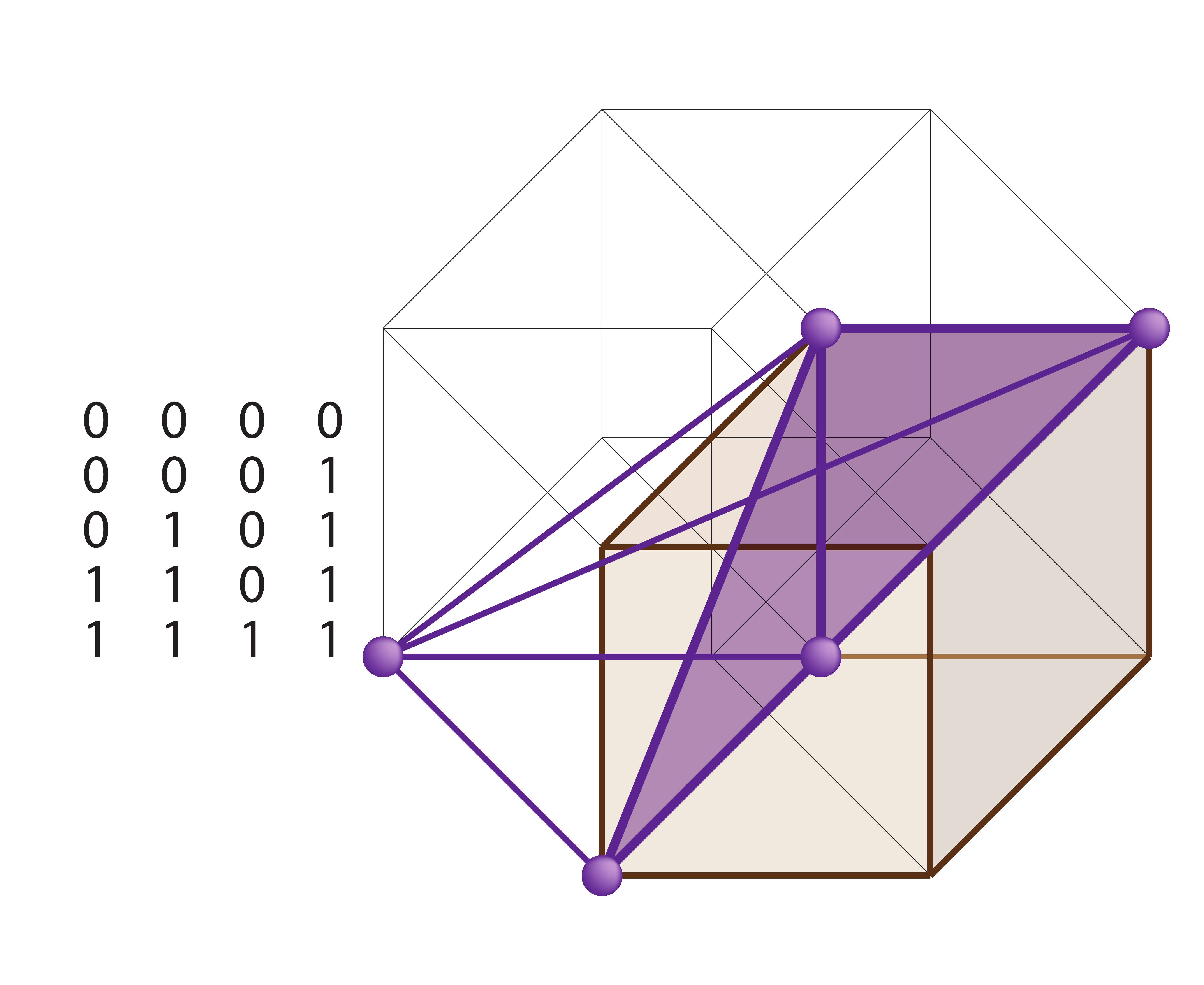}
\includegraphics[scale=.07]{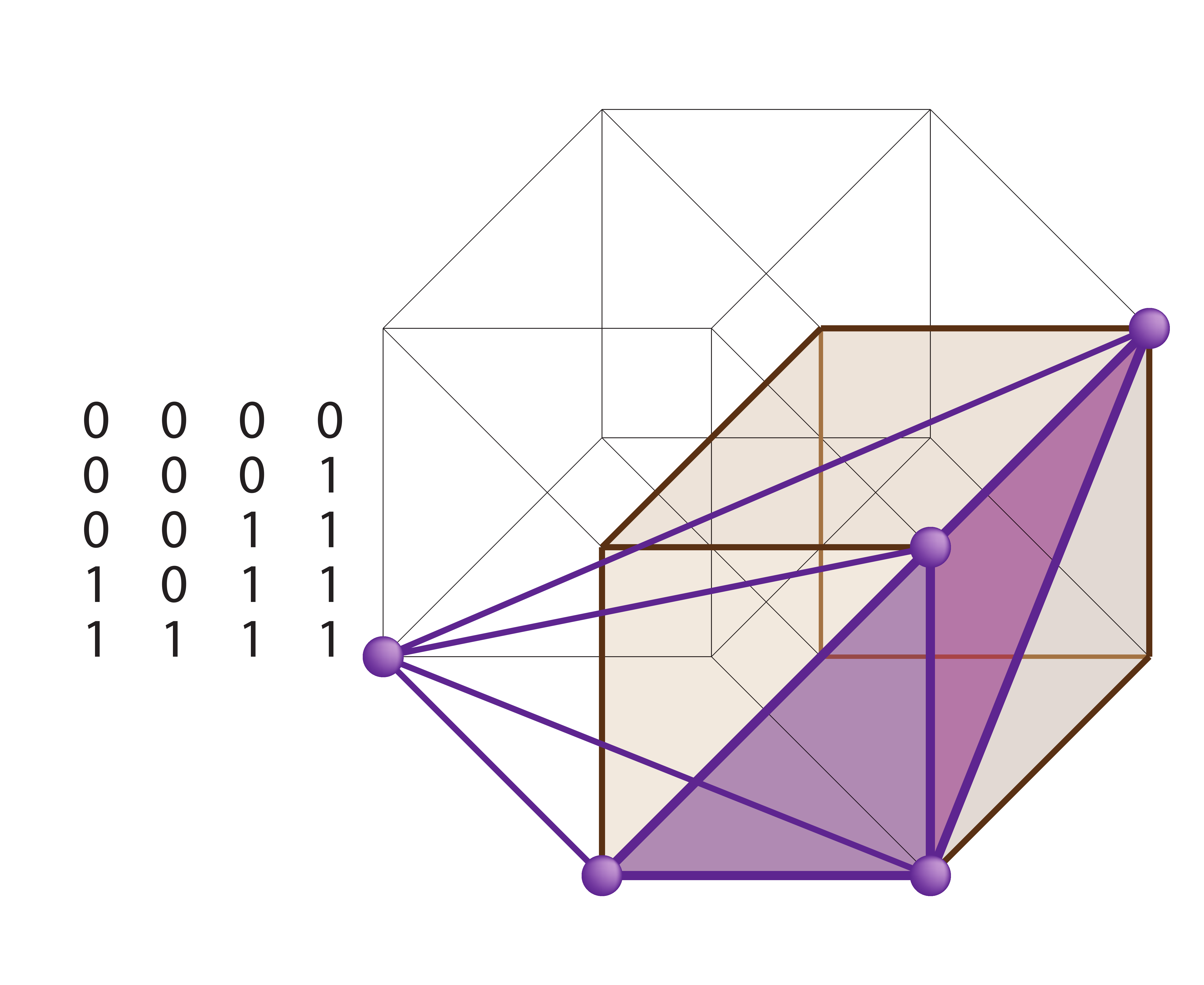}
\includegraphics[scale=.07]{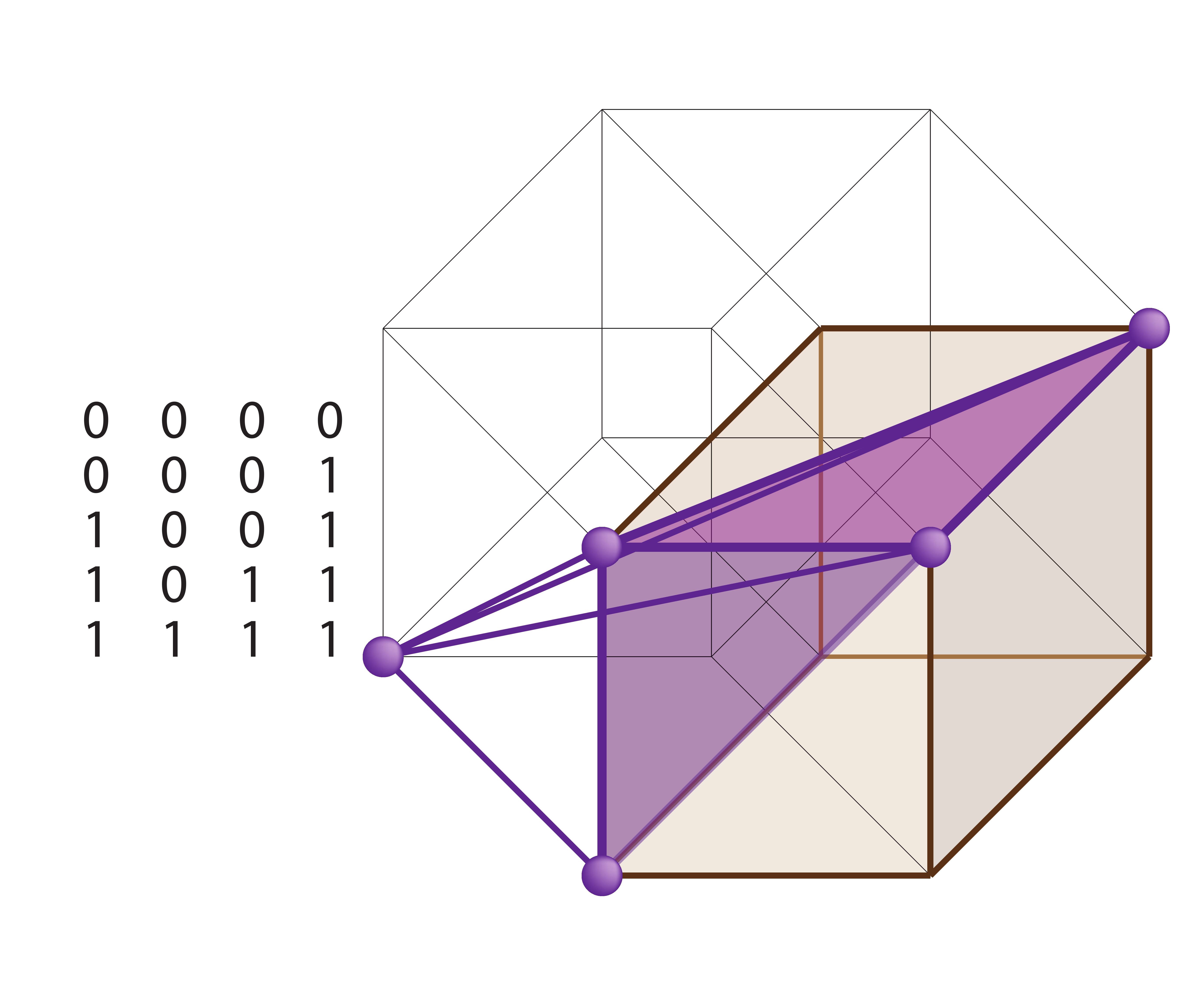}
\includegraphics[scale=.07]{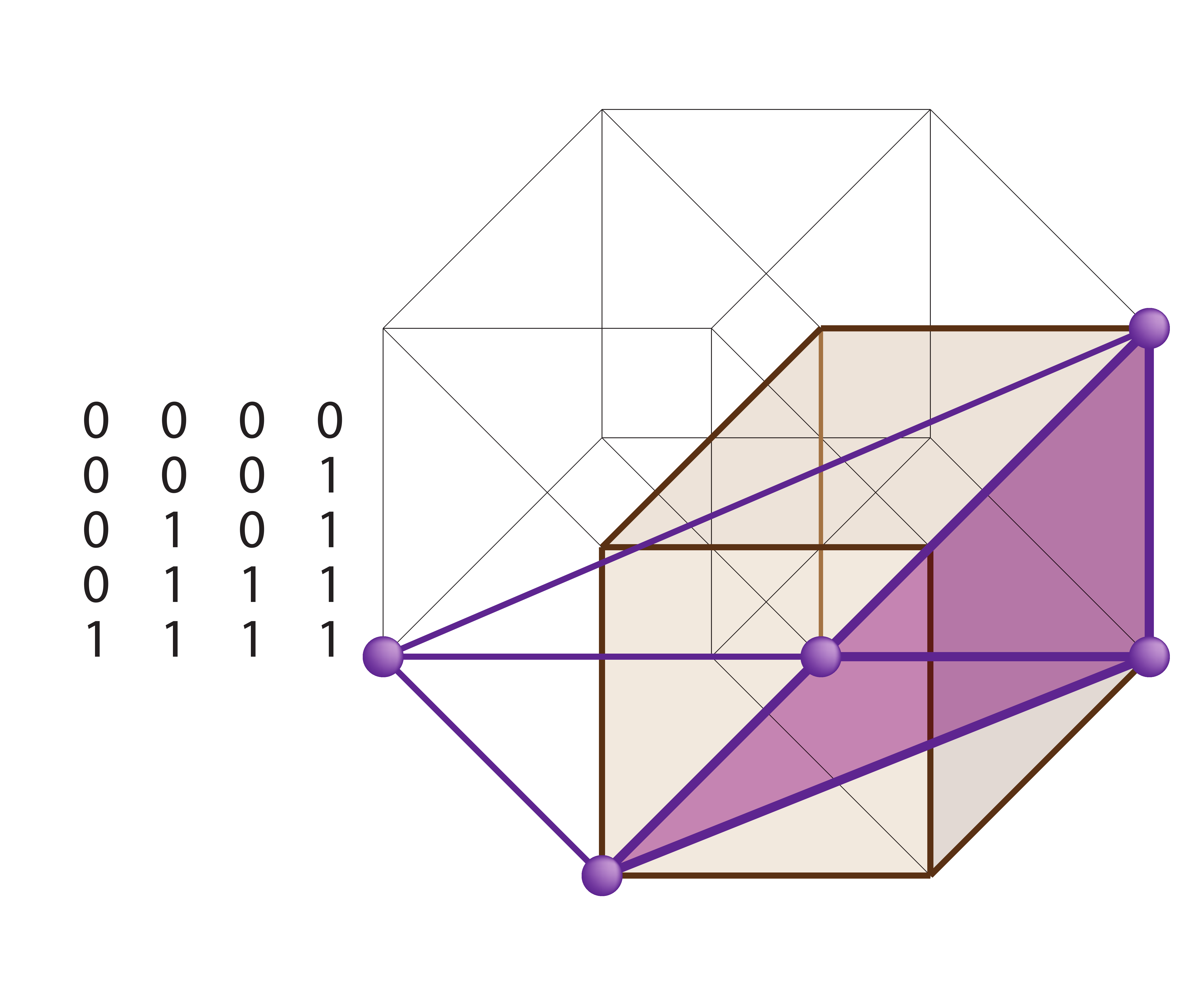}
\includegraphics[scale=.07]{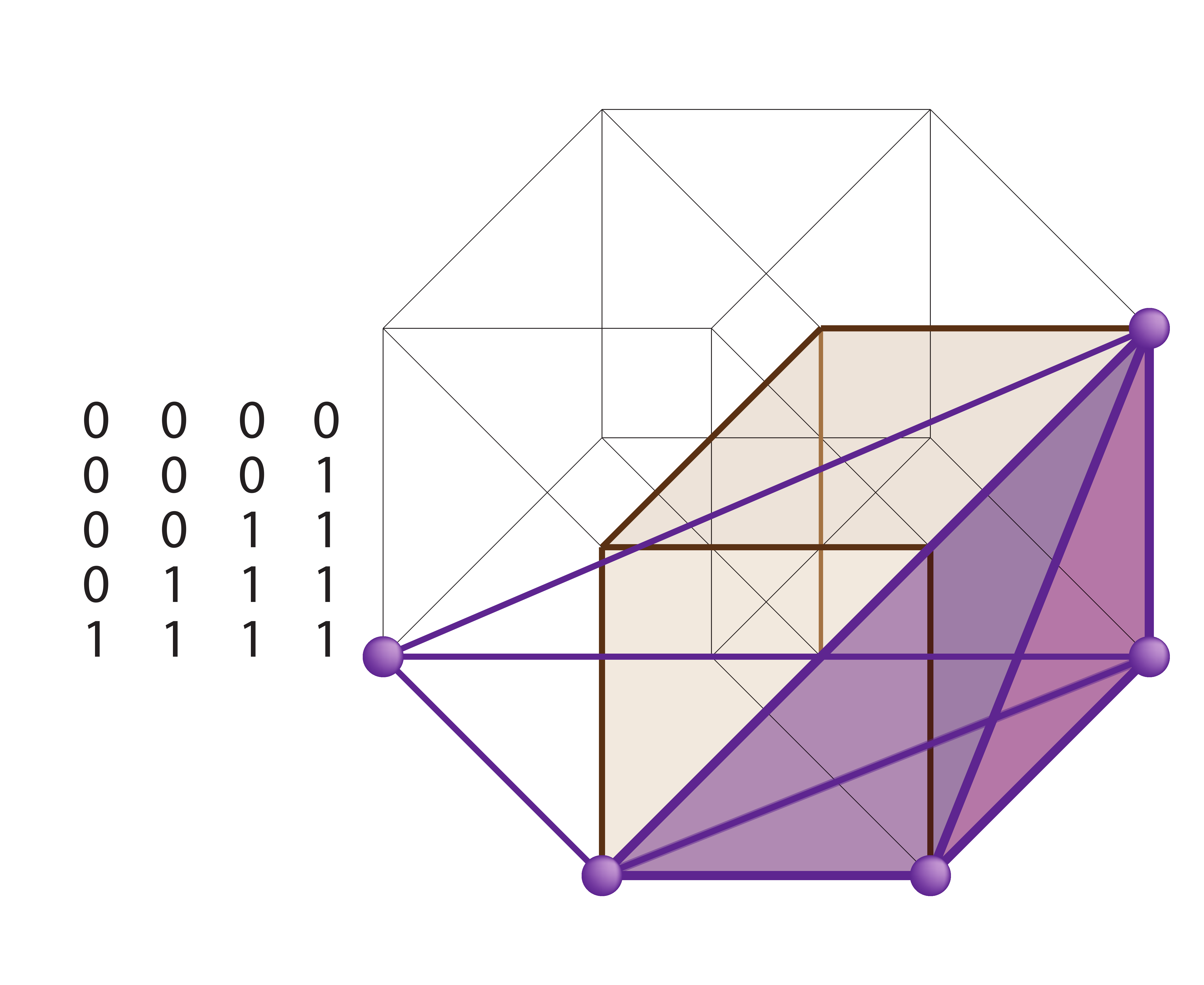}
\caption{The twenty-four $4$-simplices that fill the $4$-cube}
\label{24}
\end{figure}

The resulting rows of each matrix are the five vertices of each $4$-simplex. The union of these fills the $4$-cube. The situations are illustrated in Fig.~\ref{24}.
In this figure six copies of each $3$-dimensional cube $x=1$, $y=1$, $z=1$, and $w=1$ 
are illustrated. In each such cube a tetrahedron is drawn, and this is to be coned to the origin in the figure. The rows of the matrices that appear to the left of each hypercube indicate the coordinates of the vertices of the $4$-simplices which are drawn as small purple spheres.

\subsection{Pyramidal sets}

In the unit $n$-cube, $[0,1]^n$ consider the set $$P_n = \{ (x_1, \ldots , x_n): 0 \le x_j \le 1 \ ({\mbox{\rm for all}} \  j=1, \ldots, n) \ \& \ x_i\le x_n\le 1 \ ({\mbox{\rm for all}} \ i=1, \ldots, n-1) \}.$$ For any fixed value $x_n=C$, the cross-sectional set
$\{(x_1, \ldots , x_{n-1},C) : 0\le x_j\le C \}$ is an $(n-1)$-dimensional cube 
$[0,C]^{n-1}$ whose volume{\footnote{ We will use the term {\it volume} to mean $n$-dimensional volume}}, then is $C^{n-1}.$

Considering the situation of drawers, we can imagine a ``high priority drawer" which can be opened any distance between $0$ and $1$ meters. If it is opened to $C$ meters, then the remaining drawers can only be opened that far. Or we can imagine an alternative situation in which a king or queen (sovereign) is seated in one of $n$-chairs. All the subjects are also seated, and no-one is to stand taller than the sovereign.  As the sovereign stands, all the inferiors attempt to stand. The set of configurations of partially standing people corresponds to the configurations of partially opened drawers and of the points in the {\it pyramidal set $P_n$}. Here, note that the subscript $n$ on the expression $P_n$ indicates that the last coordinate is the position of the sovereign. 

More generally, let 
$$P_k = \{ (x_1, \ldots , x_n): 0 \le x_j \le 1 \ ({\mbox{\rm for all}} \  j=1, \ldots, n) \ \& \ x_i\le x_k\le 1 \ ({\mbox{\rm for all}} \ i \ne k ) \}.$$

Observe, that the $n$-cube is decomposed into $n$ congruent pyramids. One can then interpret this decomposition as a geometric manifestation of the fact that $$\int^1_0 x^{n-1} \ dx = \frac{1}{n}.$$
This interpretation was given in \cite{CC} and earlier in \cite{Barth}. 

\subsection{The cone on a cube is a square of a triangle}

A favorite formula from discrete mathematics is Nicomachus's Theorem
$$\sum_{k=1}^n k^3 = \left( {{n+1}\choose{2}} \right)^2.$$
One usually proves this by induction. Since 
$$n^2(n+1)^2/4 + (n+1)^3 = (n+1)^2 \left[ n^2/4 + (n+1) \right] = \frac{(n+1)^2}{4} (n^2+4n+4),$$
the proof follows easily. 
That formula is often used to compute the integral $\int^1_0 x^3 \ dx $ as the limit
${lim_{n\rightarrow \infty}} \left( \sum^N_{k=0} k^3 \right) \frac{1}{N^3}.$
The last expression comes about by subdividing the unit interval into $N$ equal subintervals and approximating the area under the curve $y=x^3$ by rectangles of base $1/N$ and height $(k/N^3)$.

   \begin{figure}[htb]
\begin{center}
\includegraphics[scale=.2]{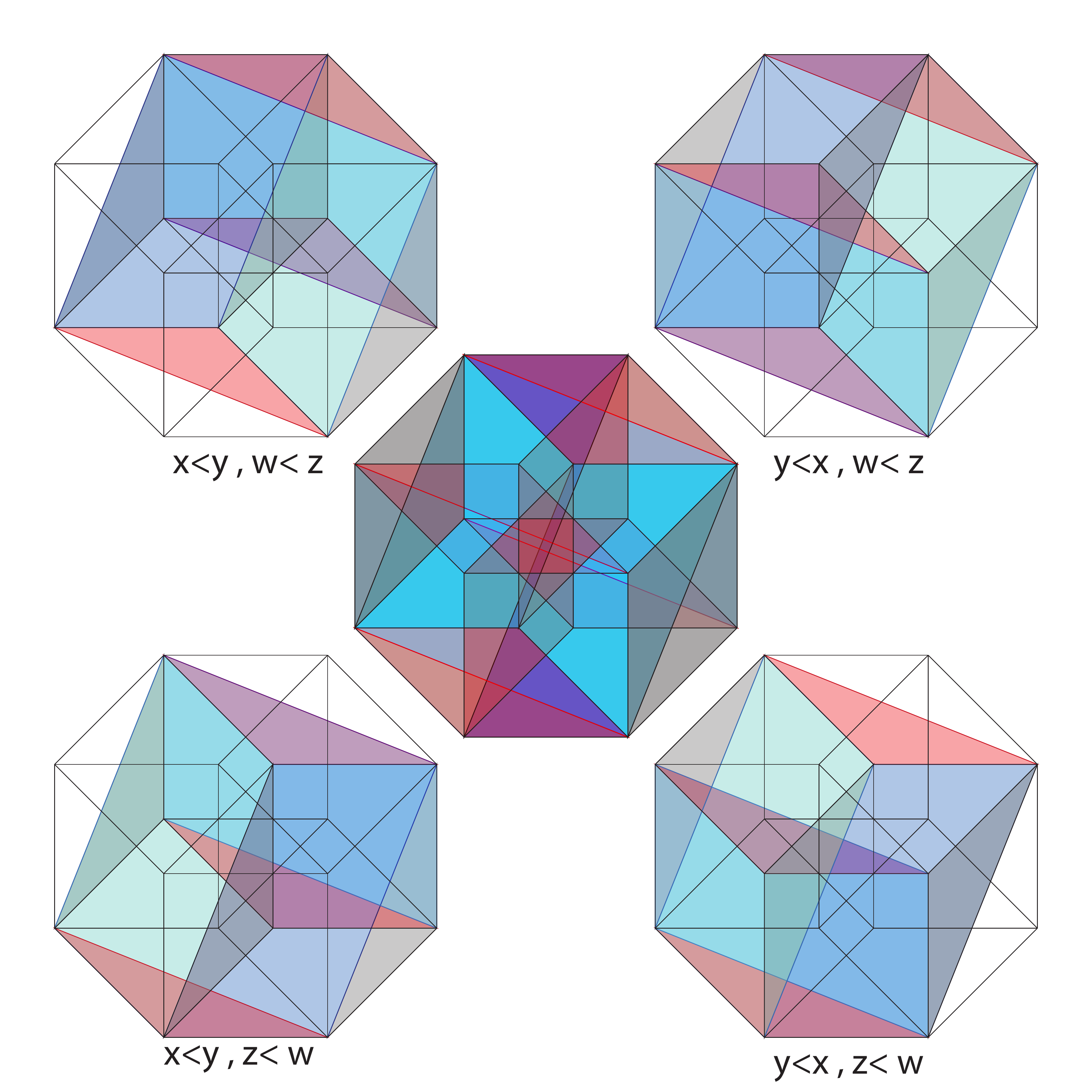}
\end{center}
\caption{Decomposing the $4$-cube as the union of four copies of $\Delta\times\Delta$}
\label{fourTxT}
\end{figure}

Dylan Thurston pointed out to JSC that one could consider the formula continuously. Thus the $4$-dimensional figure that is a cone on a cube has the same $4$-dimensional volume as the square of a triangle. Dylan's argument was developed over a dinner that he had with Dave Bayer and Walter Neumann. This section presents our modification of another dinner-time description that was given by Dror Bar-Natan in Hanoi on some evening in the period August 6-12, 2007. 

First, we consider the product of a right isosceles triangle with itself. This will be 
$$\{ (x,y,z,w): \ 0 \le x \le y \le 1 \ \& \ 0 \le z \le w \le 1 \}.$$
As a configuration of drawers, the first drawer is no more open than the second while the third drawer is no more open than the fourth. Alternatively, we consider a pair of married couples in which neither husband speaks more than his own wife, but either husband may speak more than the wife of the other. 

Now the hypercube $[0,1]^4$ can be decomposed into  four congruent figures:

 $$\Delta_{x\le y} \times \Delta_{z\le w}=  \{ (x,y,z,w): \ 0 \le x \le y \le 1 \ \& \ 0 \le z \le w \le 1 \}.$$
 $$\Delta_{y\le x} \times \Delta_{z\le w}= \{ (x,y,z,w): \ 0 \le y \le x \le 1 \ \& \ 0 \le z \le w \le 1 \}.$$
 $$\Delta_{x\le y} \times \Delta_{w\le z}= \{ (x,y,z,w): \ 0 \le x \le y \le 1 \ \& \ 0 \le w \le z \le 1 \}.$$
 and
 $$\Delta_{y\le x} \times \Delta_{w\le z}= \{ (x,y,z,w): \ 0 \le y \le x \le 1 \ \& \ 0 \le w \le z \le 1 \}.$$
 
 See Figure~\ref{fourTxT} for an illustration.
 
 And it also can be decomposed into the four pyramidal sets 
 $$P_4=\{ (x,y,z,w): \ 0 \le x , y, z \le w \le 1 \}, $$
 $$P_3 = \{ (x,y,z,w): \ 0 \le x , y, w \le z \le 1 \},$$
 $$P_2 = \{ (x,y,z,w): \ 0 \le x , z,w  \le y \le 1 \},$$
  and 
  $$P_1= \{ (x,y,z,w): \ 0 \le  y, z,w \le x \le 1 \}.$$ 
  
  Now we decompose, $\Delta_{x\le y} \times \Delta_{z\le w}$ into six $4$-simplices
  in which
  $$0\le x \le y \le z \le w\le 1, \quad 0\le x \le z \le w \le y\le 1, $$
$$0\le x \le z \le y \le w\le 1, \quad 0\le z \le x \le w \le y\le 1, $$
$$0\le z \le x \le y \le w\le 1, \quad 0\le z \le w \le x \le y\le 1. $$

and we decompose $P_4$ into six $4$-simplices:
  $$0\le x \le y \le z \le w\le 1, \quad 0\le y \le x\le z \le w\le 1, $$
$$0\le x \le z \le y \le w\le 1, \quad 0\le y \le z \le x  \le w\le 1,$$
$$0\le z \le x \le y \le w\le 1, \quad 0\le z \le y \le x \le w \le 1;$$
 the three of which listed to the left correspond to the  three that are listed to the left above.
So, half of the volume of the pyramid $P_4$ coincides with half of the volume of $\Delta_{x\le y} \times \Delta_{z\le w}.$ The other half of $P_4$ coincides with half the volume of $\Delta_{y\le x} \times \Delta_{z\le w }.$ Dylan describes this as cutting the cube $\{(x,y,z): 0 \le x,y,z \le 1\}$ into two pyramids ($(x\le y)$ and $(y \le x)$), considering the cone on each, twisting one and reassembling it to the square of the triangle.

Observe that a  $4$-cube  of the form $[0,a]^4$ also can be cut into the union of four figures each of which is the product of two isosceles right triangles of edge length $a$, $a$ and $\sqrt{2}a$. This decomposition is considered as a special case of Heron's formula, for in this case $a=b,$ and $c=a\sqrt{2}.$ The right-hand-side reduces to $4a^4$. So sixteen copies of a right isosceles triangle times itself are reassembled into four hyper-cubes.

\section{The Distributive Law, the Multinomial Theorem, and Scissors Congruences}
\label{multi}

   \begin{figure}[htb]
\begin{center}
\includegraphics[scale=.1]{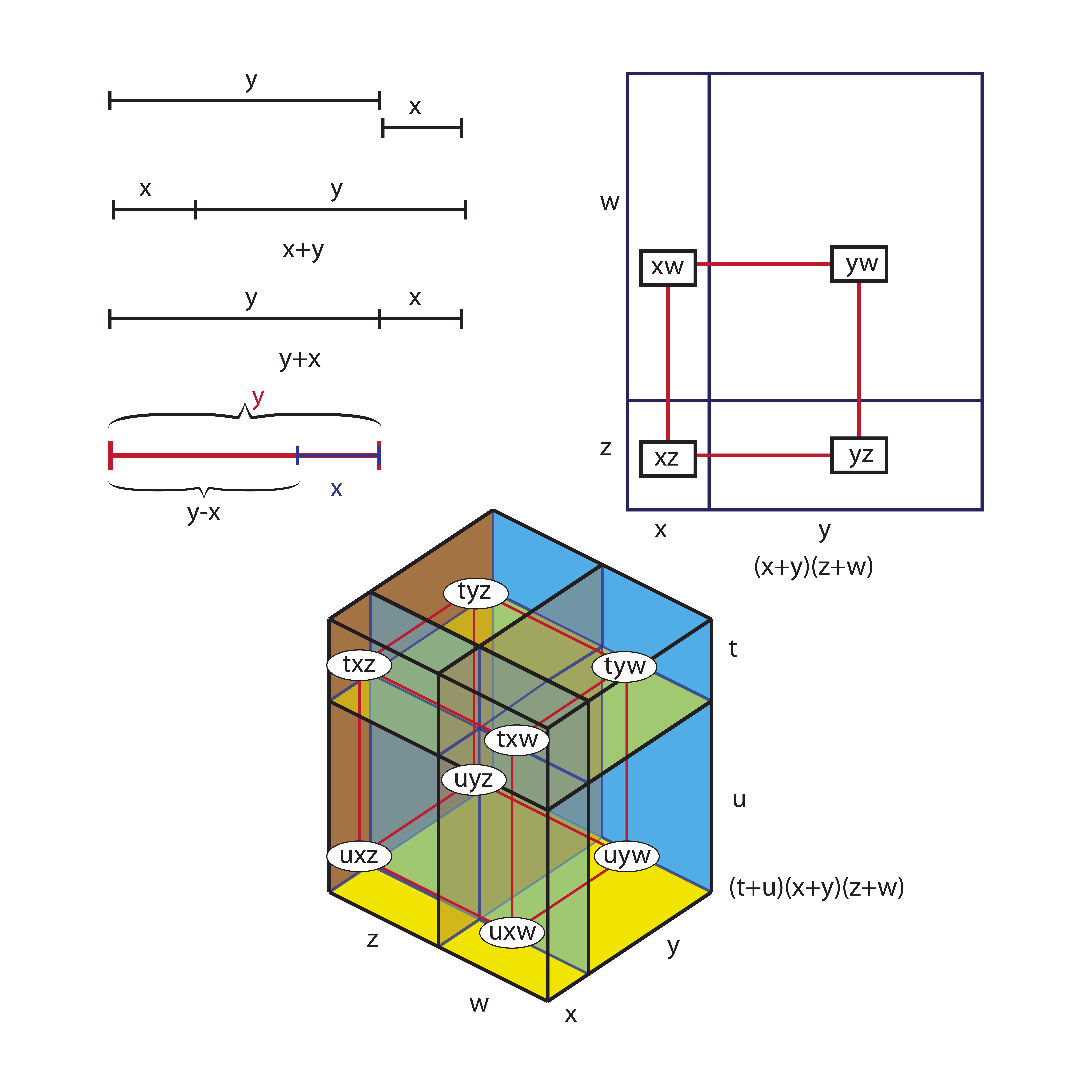}
\end{center}
\caption{The distributive law as a scissors congruence}
\label{distributive}
\end{figure}

Two $n$-dimensional polyhedra are said to be {\it scissors congruent} if and only if  they each can be cut into a union of  congruent polyhedra. For example above, we demonstrated that the $n$-cube is scissors congruent to a union of right $n$-simplices, and we demonstrated a scissors congruence between the product of a right isosceles triangle and itself with the cone on a cube. In this paper, we are primarily concerned with $4$-dimensional figures that are hypercubical, cones on three dimensional polyhedra, and the two-fold product of polygons. 

Famously, Max Dehn \cite{Dehn,Benko} demonstrated that a regular tetrahedron of unit volume is not scissors congruent with a unit cube. 

To begin this section, in which an $n$-dimensional cube with variable edge lengths is decomposed into several smaller sub-cubes, we first consider the case $n=1$.

A line segment whose length is the sum  $x+y$ of two positive quantities $x$ and $y$ can be cut into two segments with the segment of length $x$ on the left and the segment of length $y$ on the right (or vice versa). Similarly, suppose that $x<y$ and a segment of length $y-x$ is given. It could be produced by means of cutting a segment of length $x$ from the longer segment. Of course, we can, in general, discuss directed segments with the direction ``right" indicating positive length quantities, and the direction ``left" used to indicate negative quantities. In this way, each of $x+y$, $x-y$, and $y-x$ can be thought of as being built by composing directed segments of lengths $x$ and $y$. It is conceptually easier to think in this fashion rather than explicitly removing segments and hypercubes whose total volumes in algebraic expressions happen to be negative.

In case a given rectangular area of dimensions $(x+y)$-by-$z$ can be decomposed as a pair of rectangular areas $x$-by-$z$ and $y$-by-$z$. As above, we can indicate $(x-y)z$ as the difference between areas $xz$ and $yz$. We further exploit the distributive law as a scissors congruence when considering the case $(x+y)$-by-$(z+w)$. This rectangular area is decomposed as the union of  four rectangles of size $xz$,$xw$, $yz$ and $yw$. Some of these scissors congruences are indicated in Fig.~\ref{distributive}. We leave the reader to imagine others.

For convenience within cubes, hypercubes, and beyond, we dualize the pictures by indicating a vertex with a label that indicates its volume (and is drawn red if the net volume is negative) and connecting these by an edge if they share a codimension-$1$ face. If they share a
codimension $2$ face, then the dual  vertices will bound a polygon (always a rectangle in our cases), and so forth. The language of shared polygonal vertices, edges, faces, and so forth is initially quite awkward (thus the term ``codimension" is used), but in our cases the situation is a bit more straight-forward than one might presume. 

Our main interest here is to consider the $n$-dimensional cube that is decomposed as a product of edges in the form
$$RH=(a+b+c)(a+b-c)(a-b+c)(-a+b+c).$$
This expression $RH$ is called the {\it right-hand-side of Heron's formula}.
When we use the distributive law to expand the product as the union of $81$ hypercubes (more precisely hyper-rectangles), we want to identify which terms will cancel. It is not terribly difficult to determine that the product reduces to
$$2a^2b^2+2a^2c^2+2b^2c^2 -a^4-b^4-c^4.$$
``Not terribly difficult" is a relative phrase. An electronic computer can do so easily; a human has to develop some technique to organize and regroup the 81 terms that are initially involved in the expansion. 
Our purpose here will be to find within the hyper-solid the location of each piece. We will observe that the pieces that cancel can be found on specific affine $3$-dimensional planes.  

To help determine the addresses of the vertices in the dual complex, we will examine the terms in the multinomial expansion 
$$(x_1+x_2+\cdots + x_k)^n.$$
This situation is more general in that we are considering an $n$-dimensional quantity, but it is too specific for the case at hand since the terms in the expansion of $RH$ involve differing edge lengths. Still, once the addressing problem is solved for the multinomial expansion, it will be easy to apply here.

The quantity $(x_1+x_2+\cdots + x_k)^n$ corresponds to decomposing $[0, \sum_k x_k]^n$ into $k^n$ pieces. The volumes of each piece is a product of $x_i\/s$ where each $x_i$ is taken from one factor of the decomposition of the edges of the $n$-cube. Thus there is a specific $n$-cube{\footnote{The correct terminology might be $n$-rectangle since the edge lengths vary. Here we will use the more concise terminology, but elsewhere, we may speak of ``hyper-rectangles."}} of volume $x_{i_1} x_{i_2} \cdot \cdots \cdot x_{i_n}$ where each $i_\ell \in \{1,\ldots,k\}$ and this $\ell\/$th factor $x_{i_\ell}$ was selected from the $\ell\/$th factor in the product expansion:
 $$(x_1+ \cdots + x_k) \cdots \underbrace{(x_1+\cdots + \framebox{$x_{i_\ell}$}+ \cdots + x_k )}_{\ell\/{\mbox{\rm th factor}}} \cdots (x_1+ \cdots + x_k)$$

In order to keep track of these terms, we consider the coordinates of the points in the integral lattice that lie within the cube $[0,(k-1)]^n$. These coordinates range between $(0,0,\ldots,0)$ and $(k-1,k-1,\ldots, k-1)$.  We then label the point $(i_1-1,i_2-1, \ldots, i_n-1)$ in the integral lattice with the term $x_{i_1}x_{i_2} \cdot \cdots \cdot x_{i_n}$ in the expansion above.

The next step in understanding the multinomial theorem is to count all of the terms that have the same volume ($x_{i_1}x_{i_2} \cdot \cdots \cdot x_{i_n}$) as the cube $[0,x_{i_1}]\times[0,x_{i_2}]\times \cdots \times [0, x_{i_n}]$ does.
So we can think of $\{ \{ i_1-1,i_2-1, \ldots, i_n-1 \} \}$ as a multi-set --- that is a set possibly with repetitions. The products $x_{i_1}x_{i_2} \cdot \cdots \cdot x_{i_n}$ and $x_{j_1}x_{j_2} \cdot \cdots \cdot x_{j_n}$ have the same volumes if and only if the corresponding multisets are the same or that the sequences $(i_1-1,i_2-2, \ldots, i_n-1)$ and $(j_1-1,j_2-1, \ldots, j_n-1)$ are permutations of each other. 
If so, then these points lie along the same $(n-1)$-dimensional affine plane 
$\sum_{j=1}^n z_j = K$, where $K=\left( \sum_{\ell=1}^n x_{i_\ell} \right) -n.$
The value $K$ ranges between $0$ and $n(k-1).$

   \begin{figure}[htb]
\begin{center}
\includegraphics[scale=.1]{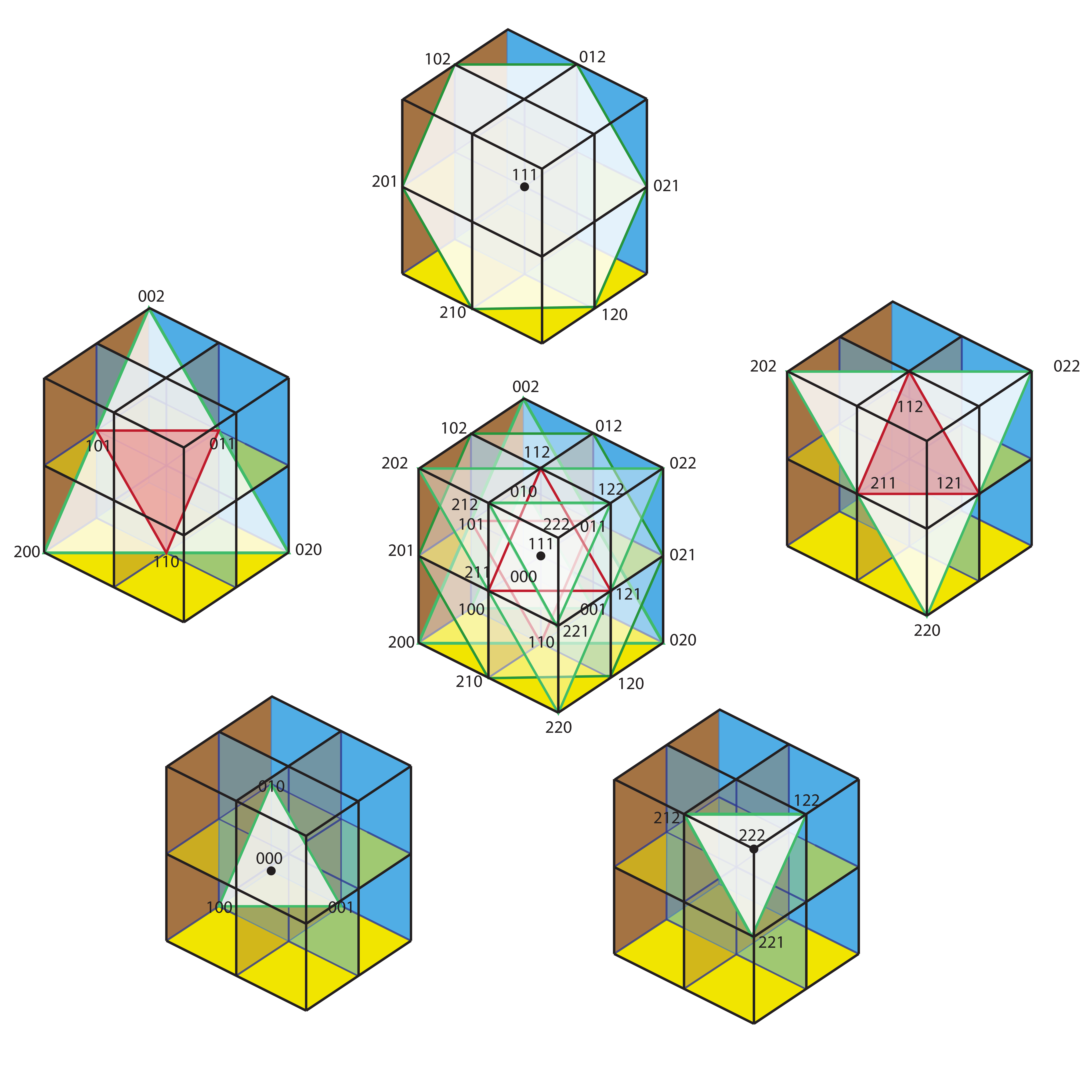}
\end{center}
\caption{The trinomial theorem $n=3$. Note $ijk$ corresponds to $x_{i+1}x_{j+1}x_{k+1}$.}
\label{trinomial}
\end{figure}

This last condition is necessary, but not sufficient. For example, in $\R^3$ the plane $x+y+z=2$ intersects the lattice in the points $(2,0,0),(0,2,0),(0,0,2)$ (which correspond respectively to the products $x_3 x_1 x_1$, $x_1 x_3 x_1$, and $x_1x_1x_3$) and the points $(1,1,0)$ $(1,0,1)$, $(0,1,1)$ (corresponding  respectively to the products $x_2x_2x_1$, $x_2 x_1 x_2$ and $x_1x_2x_2$). 

In the illustration Fig.~\ref{trinomial}, the $27$ terms in the expansion of $(x_1+x_2+x_3)^3$ are grouped according to the intersection with the plane $x+y+z=K$
as also indicated in the table below.

\begin{center}
\begin{tabular}{||c|c||} \hline \hline
$(0,0,0)$ & $x_1^3$ \\ \hline
$(1,0,0)$, $(0,1,0)$,$(0,0,1)$ & $3x_1^2 x_2$ \\ \hline
$(2,0,0)$, $(0,2,0)$,$(0,0,2)$ & $3x_1^2 x_3$ \\ \hline
$(1,1,0)$, $(1,0,1)$,$(0,1,1)$ & $3x_1x_2^2 $ \\ \hline
$(1,1,1)$ & $ x_2^3$ \\ \hline
$(2,2,0)$, $(2,0,2)$,$(0,2,2)$ & $3x_1x_3^2 $ \\ \hline
$(1,1,2)$, $(1,2,1)$,$(2,1,1)$ & $3x_2^2x_3 $ \\ \hline
$(2,2,1)$, $(2,1,2)$,$(1,2,2)$ & $3x_2x_3^2 $ \\ \hline
$(2,2,2)$ & $x_3^3$ \\ \hline \hline
\end{tabular}
\end{center}

The full statement of the multinomial theorem is 
$$(x_1+x_2 + \cdots + x_k)^n = \sum_{i_1+i_2+ \cdots + i_\ell=n} \left( \frac{n!}{i_1! i_2! \cdots i_\ell!} \right) x_{1}^{i_1}\cdots x_{2}^{i_2} \cdots x_{\ell}^{i_\ell}.$$
We envision the multinomial coefficient as an expression of the symmetry of the intersection of the affine $(n-1)$-plane $\sum_{j=1}^n z_j = K$ or equivalently the symmetry of the multi-set $$\{\{ \underbrace{x_1, \cdots,x_1}_{i_1}, \underbrace{x_2, \cdots,x_2}_{i_2},\cdots \underbrace{x_\ell, \cdots,x_\ell}_{i_\ell} \} \}.$$

   \begin{figure}[htb]
   \begin{center}
\includegraphics[scale=.07]{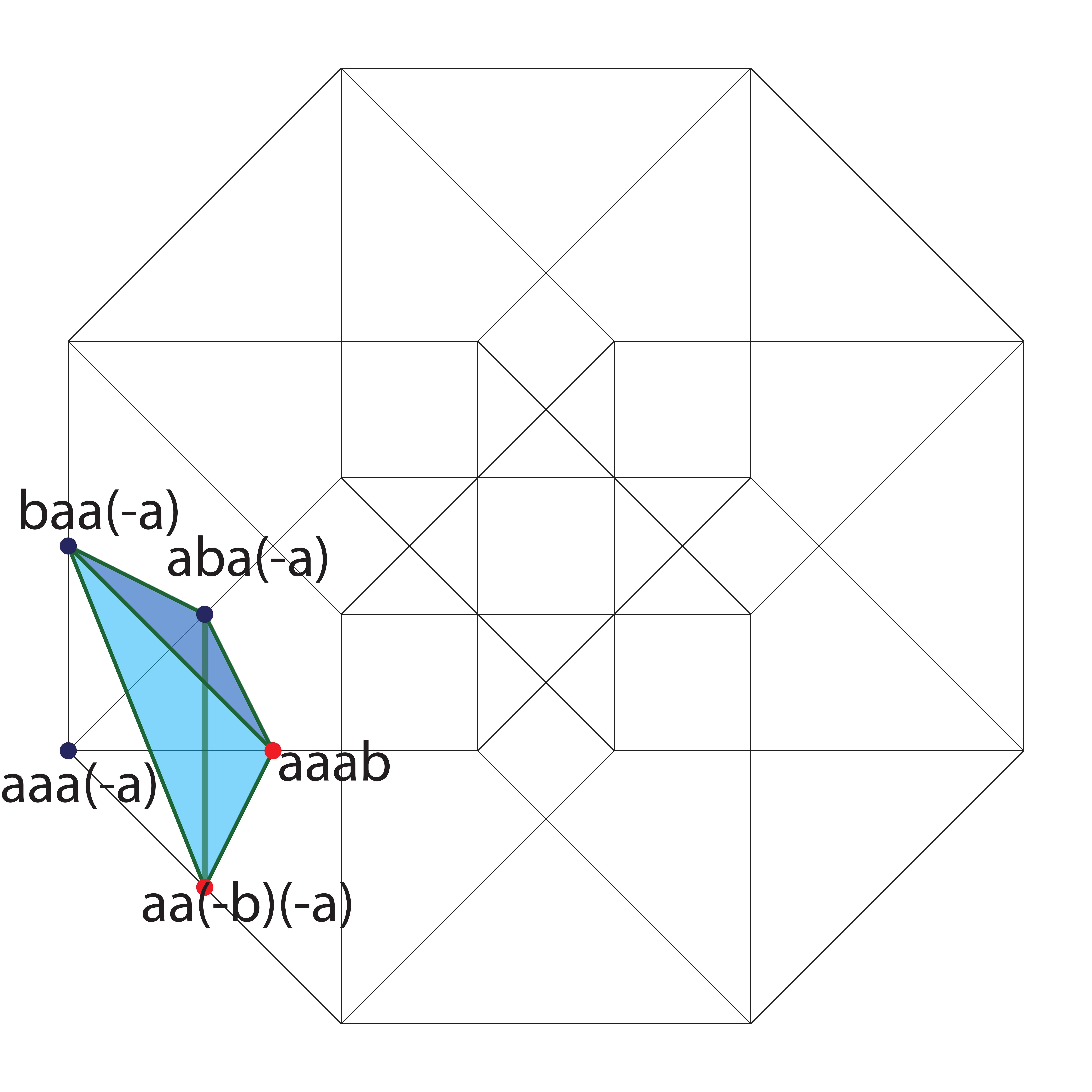}
\includegraphics[scale=.07]{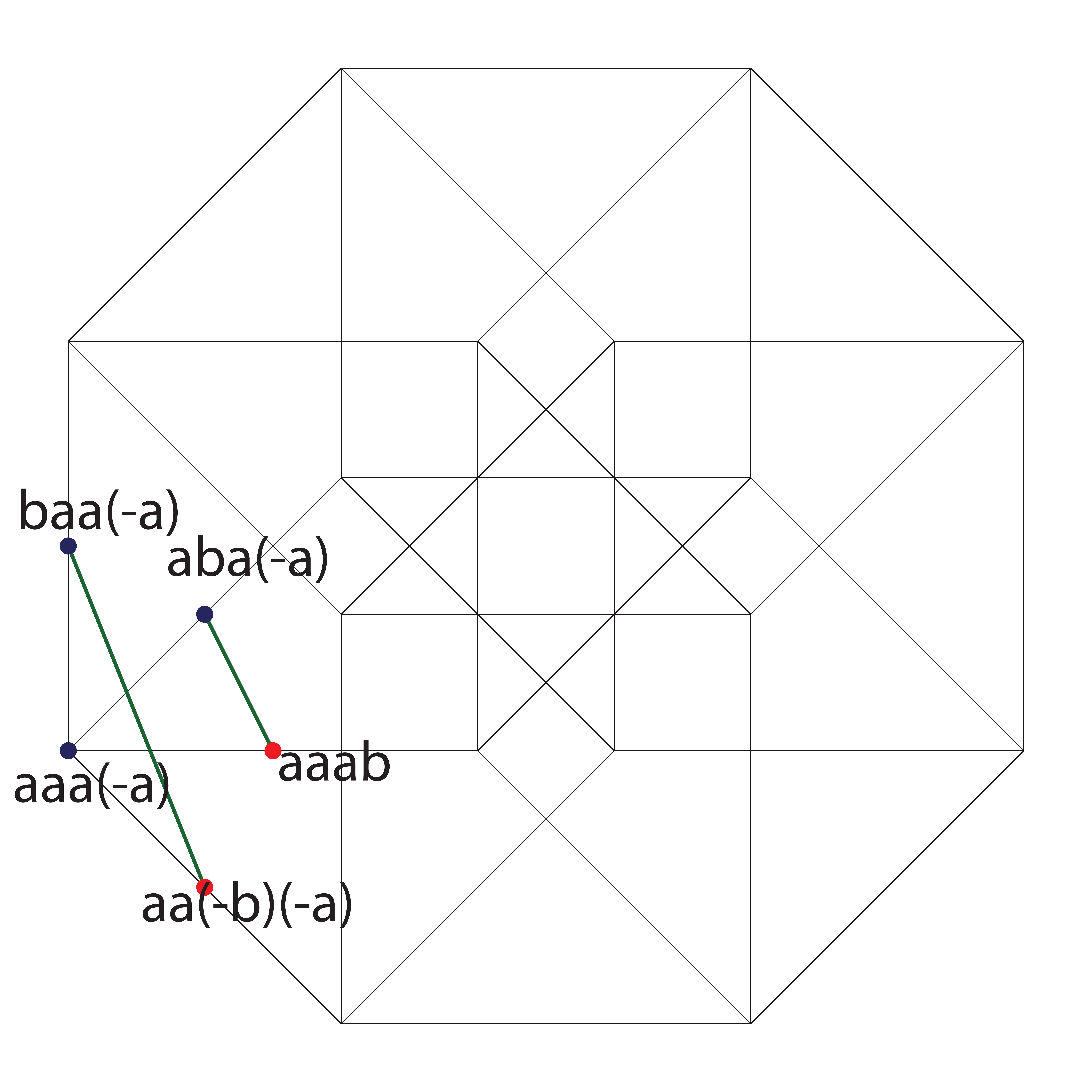}
\end{center}
\caption{The terms of the form $a^4$ and $a^3b$ in expansion of RHS }
\label{sudoko1}
\end{figure}

   \begin{figure}[htb]
      \begin{center}
\includegraphics[scale=.065]{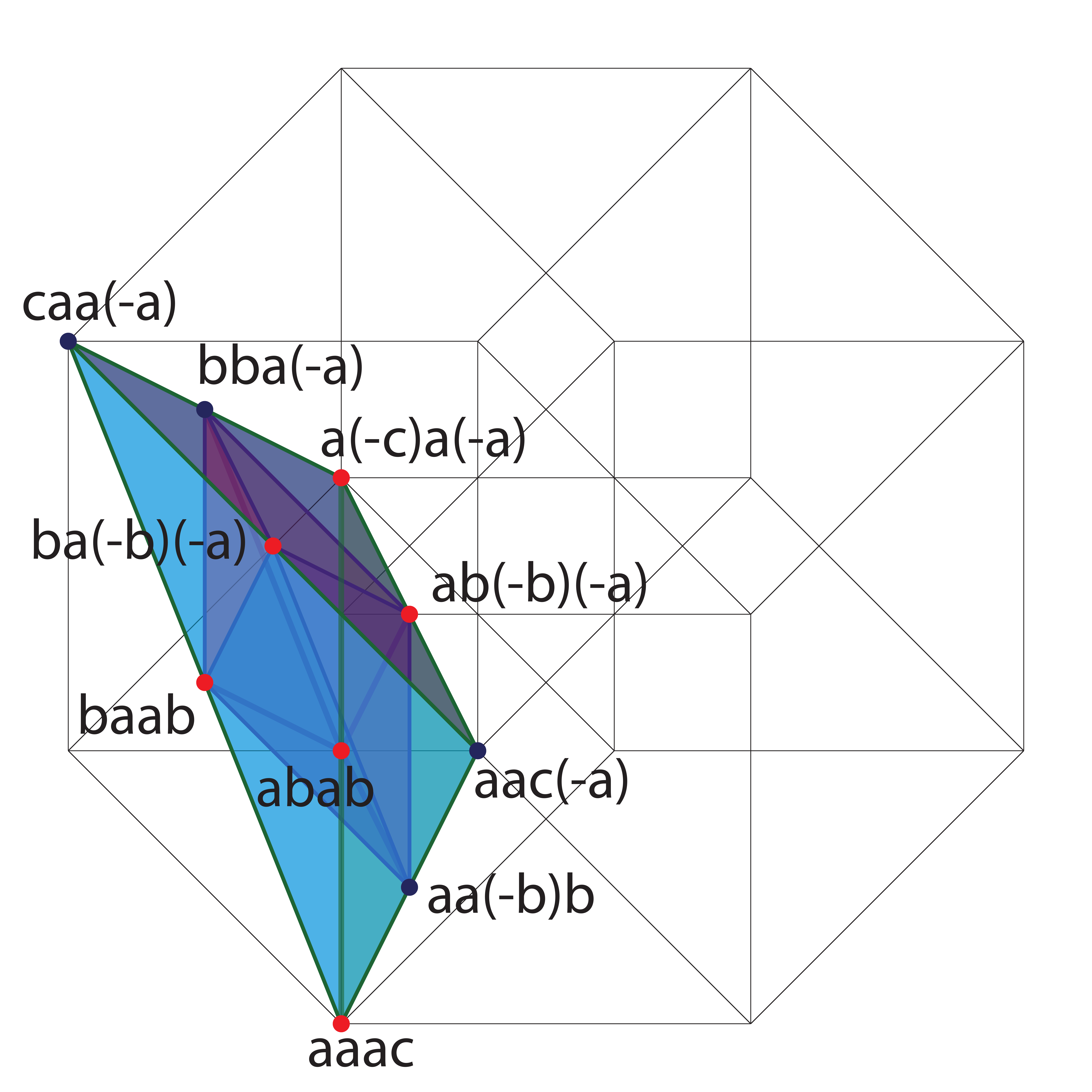}
\includegraphics[scale=.065]{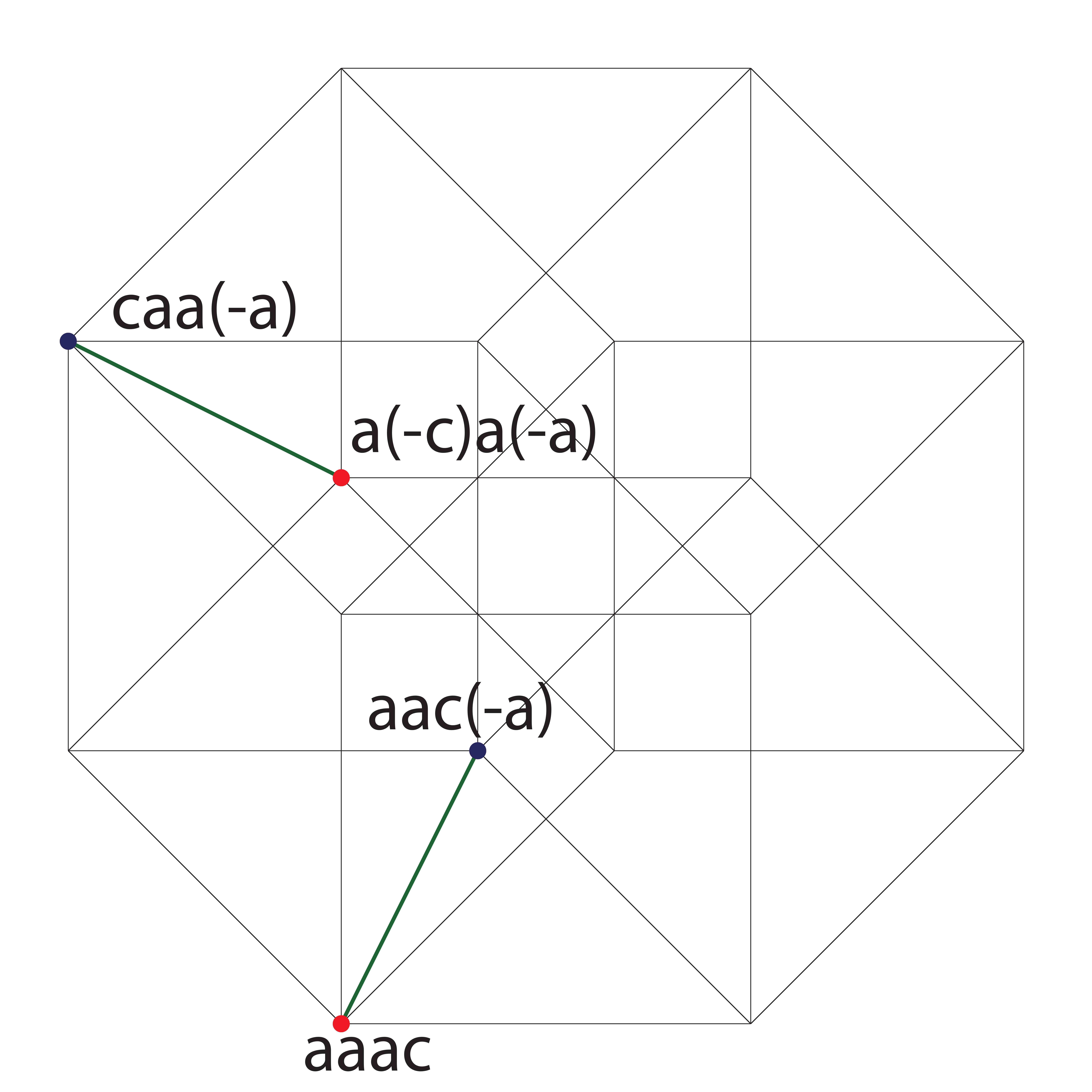}
\includegraphics[scale=.065]{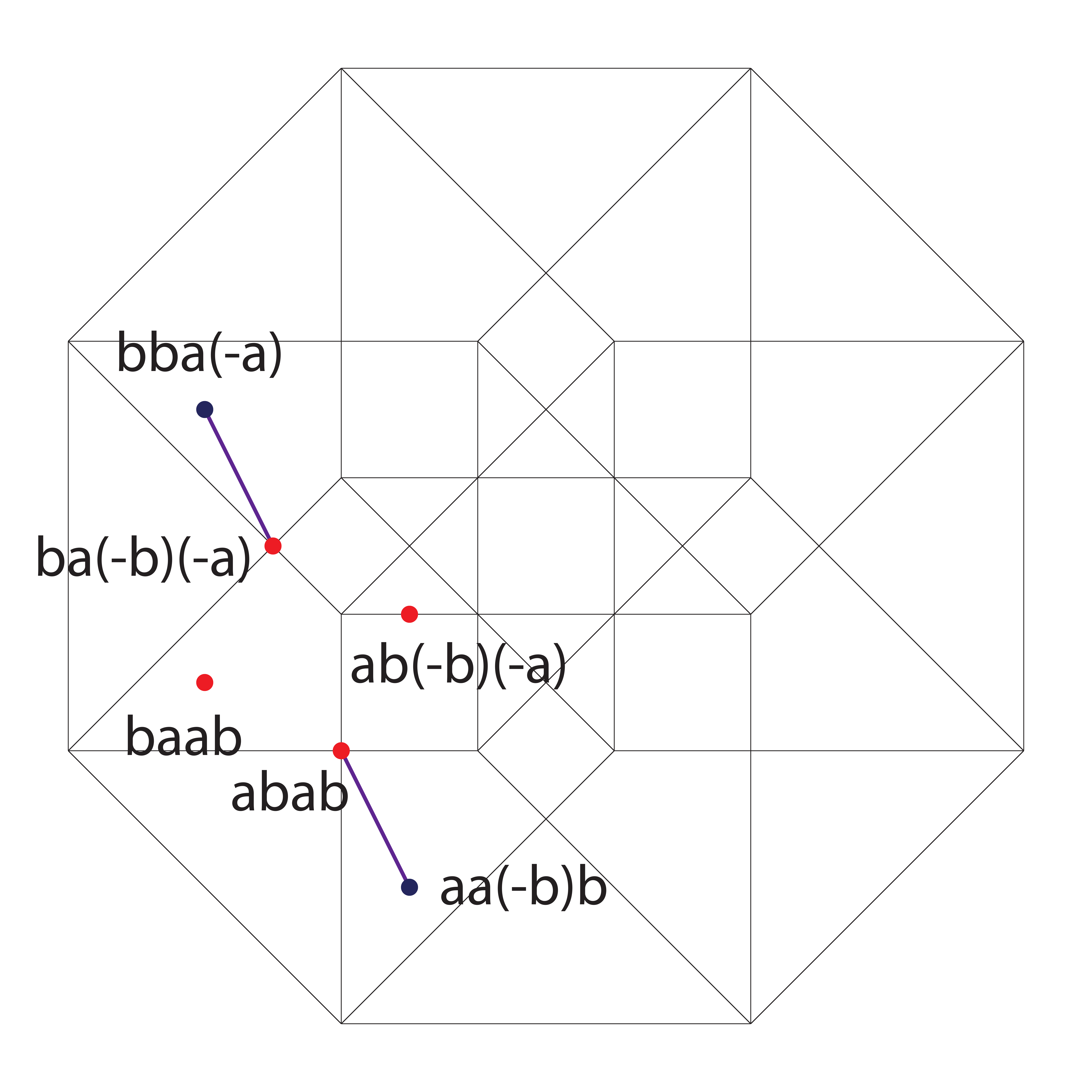}
\end{center}
\caption{The terms of the form $a^3c$ and $a^2b^2$ in the expansion of RHS}
\label{sudoko2}
\end{figure}

   \begin{figure}[htb]
   \begin{center}
\includegraphics[scale=.07]{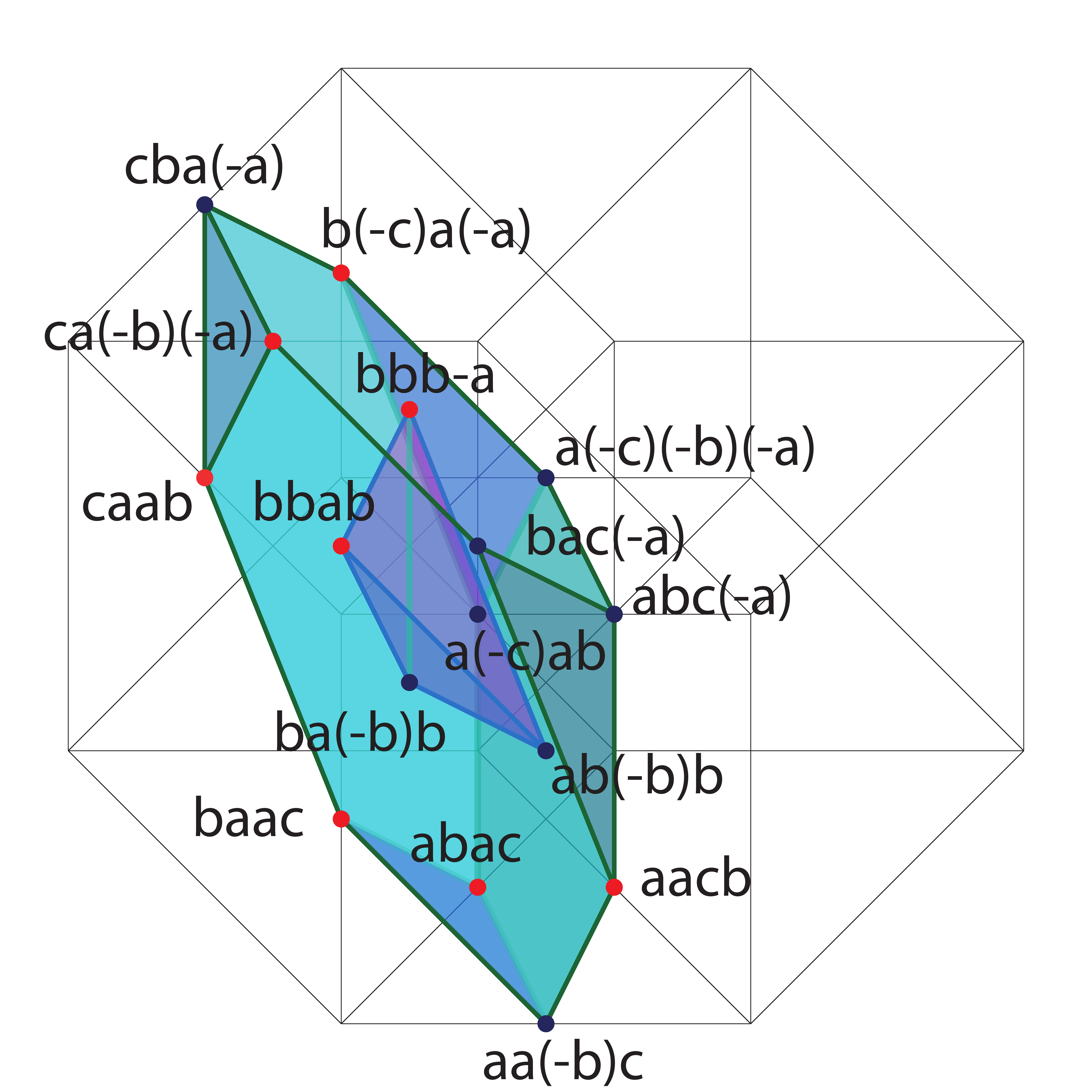}
\includegraphics[scale=.07]{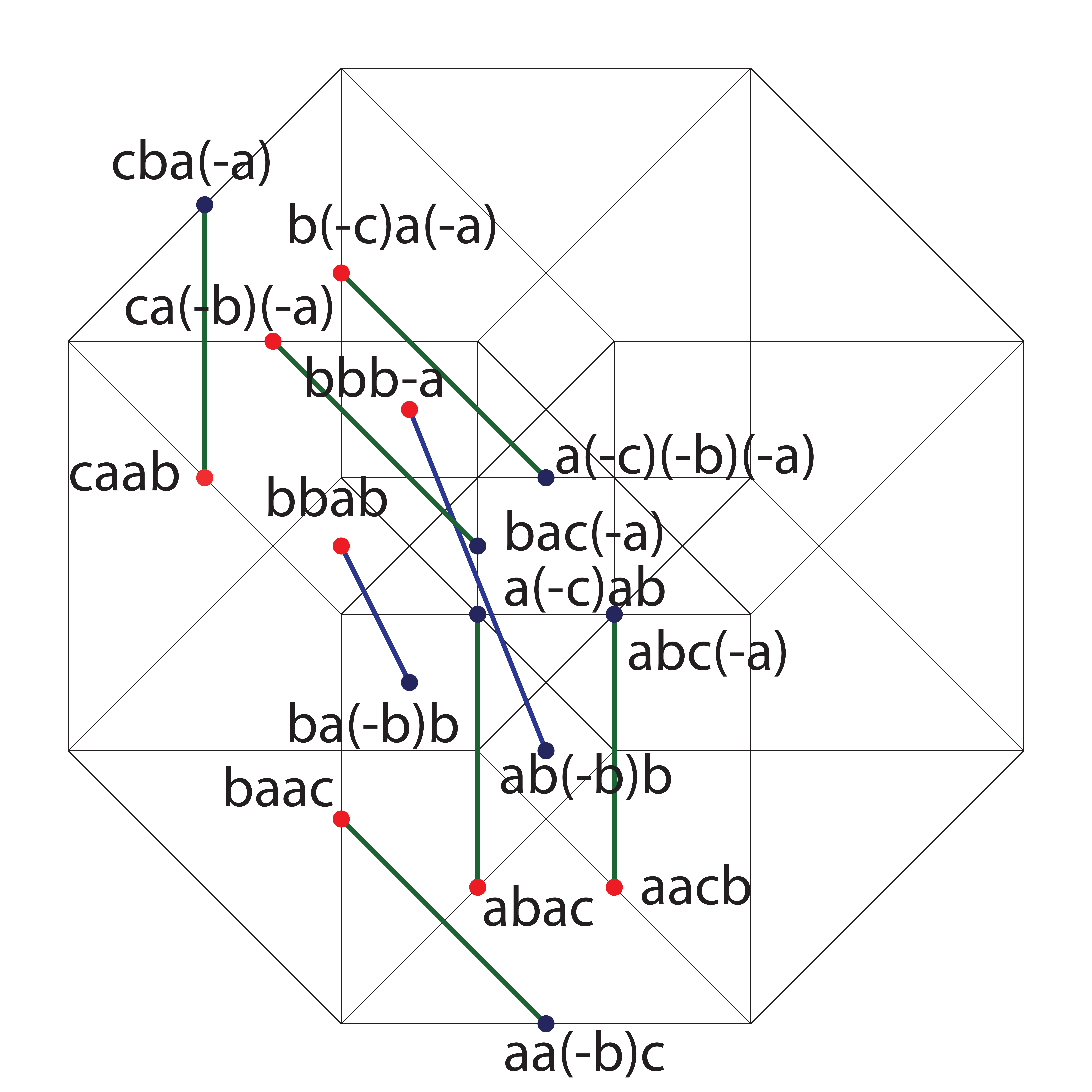}
\end{center}
\caption{The terms of the form $a^2bc$ and $ab^3$ in the expansion of RHS} 
\label{sudoko3}
\end{figure}


   \begin{figure}[htb]
   \begin{center}
\includegraphics[scale=.065]{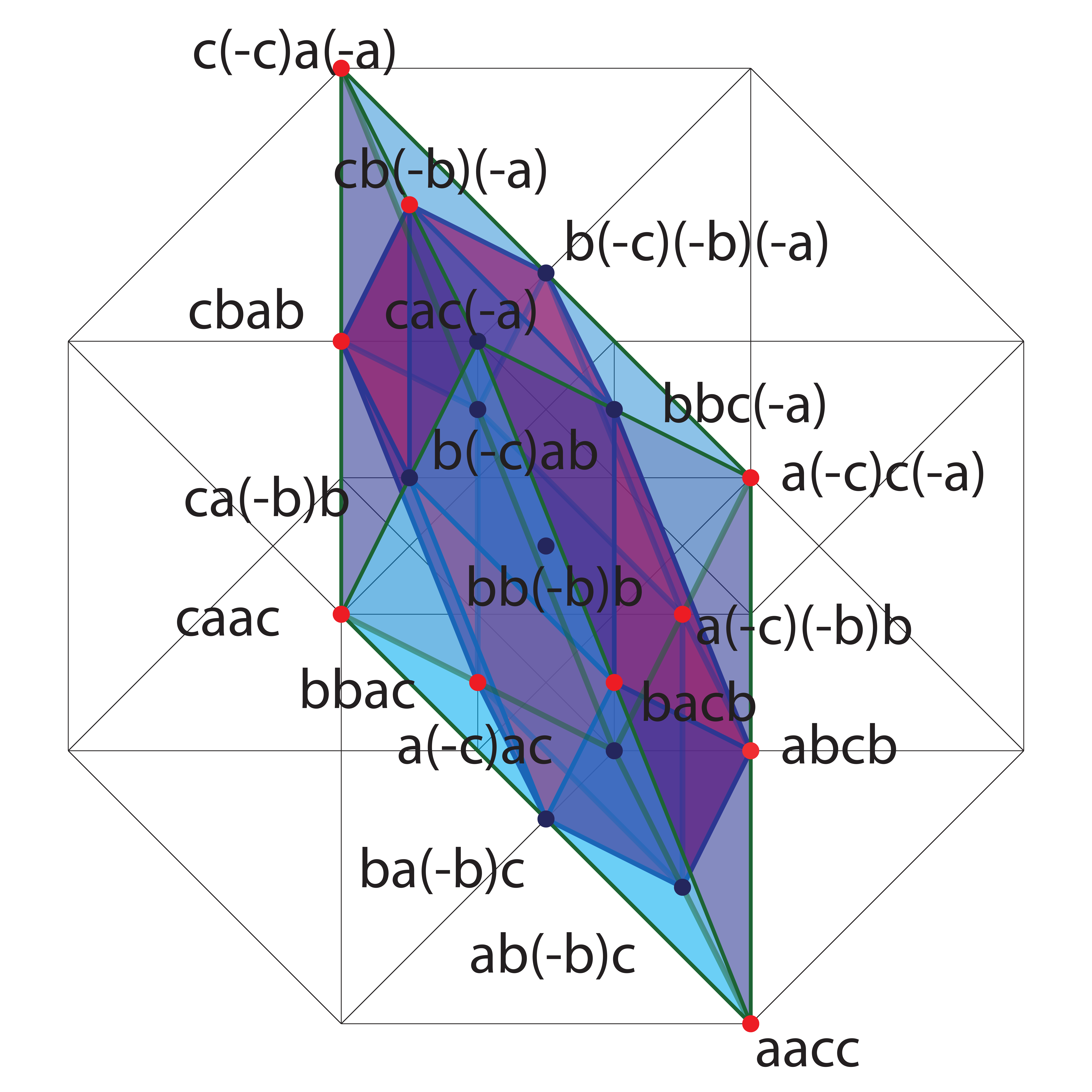}
\includegraphics[scale=.065]{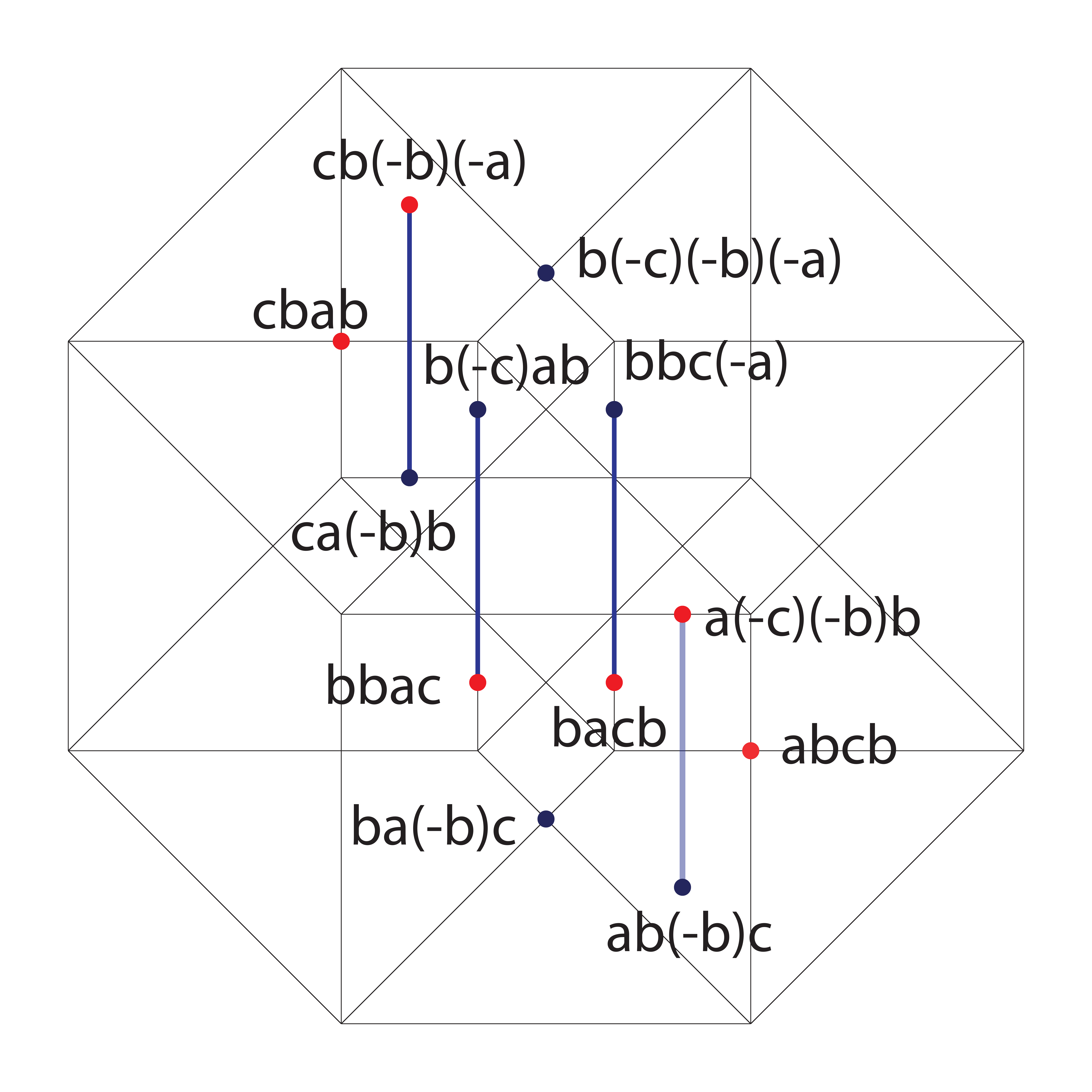}
\includegraphics[scale=.065]{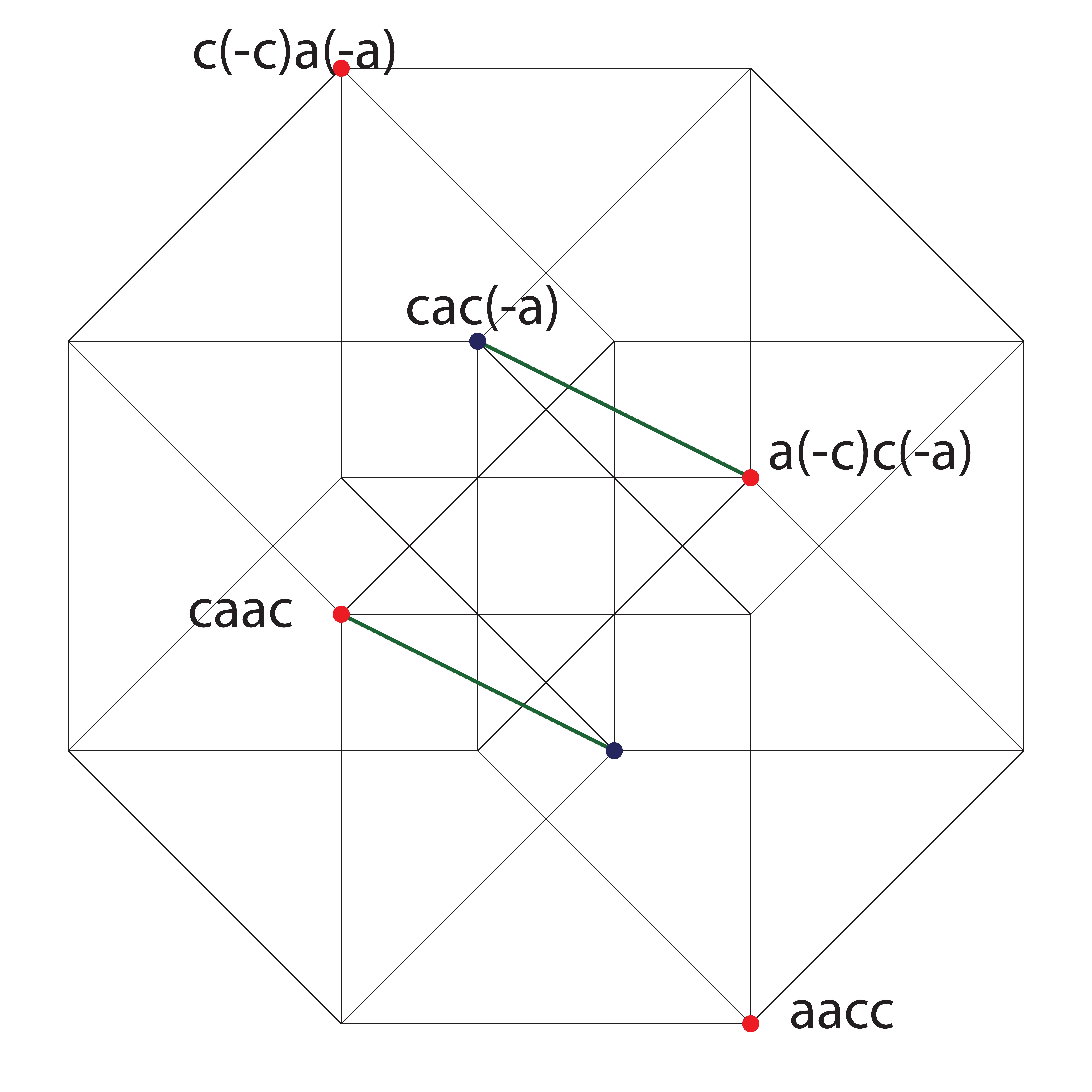}
\end{center}
\caption{The terms of the form $ab^2c$ and $a^2c^2$ in the expansion of RHS} 
\label{sudoko4}
\end{figure}

   \begin{figure}[htb]
   \begin{center}
\includegraphics[scale=.07]{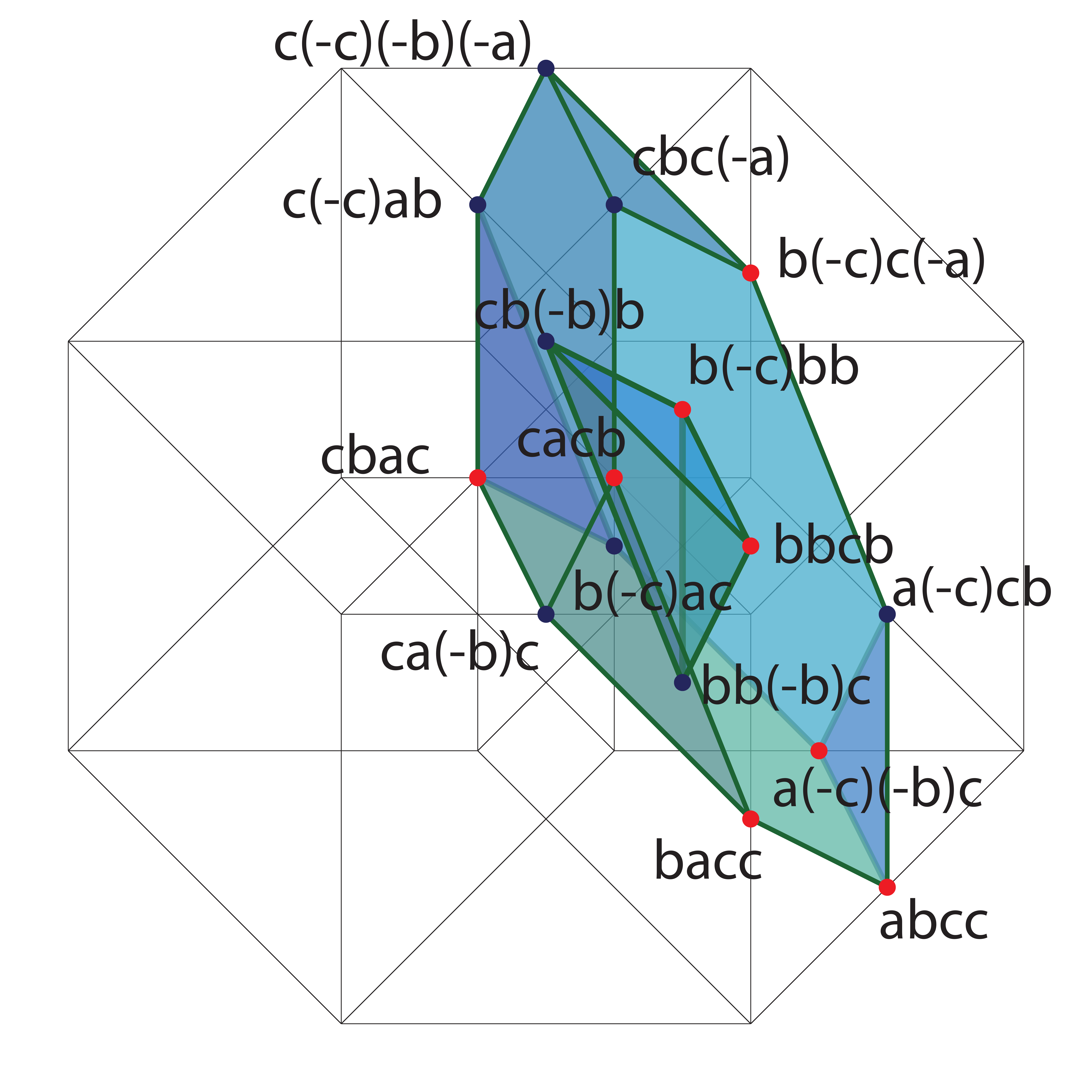}
\includegraphics[scale=.07]{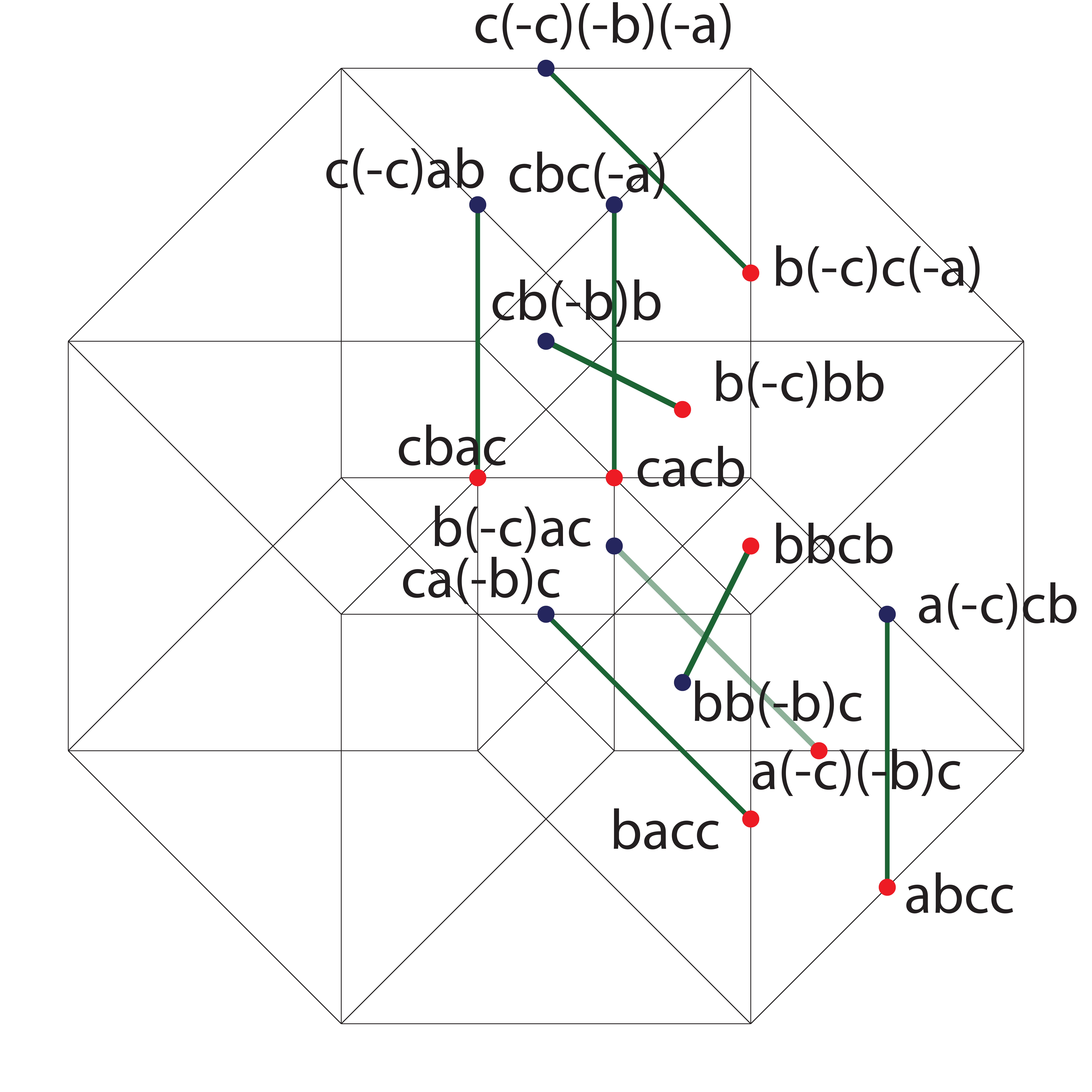}
\end{center}
\caption{The terms of the form $abc^2$ and $b^3c$ in the expansion of RHS} 
\label{sudoko5}
\end{figure}

 \begin{figure}[htb]
   \begin{center}
\includegraphics[scale=.065]{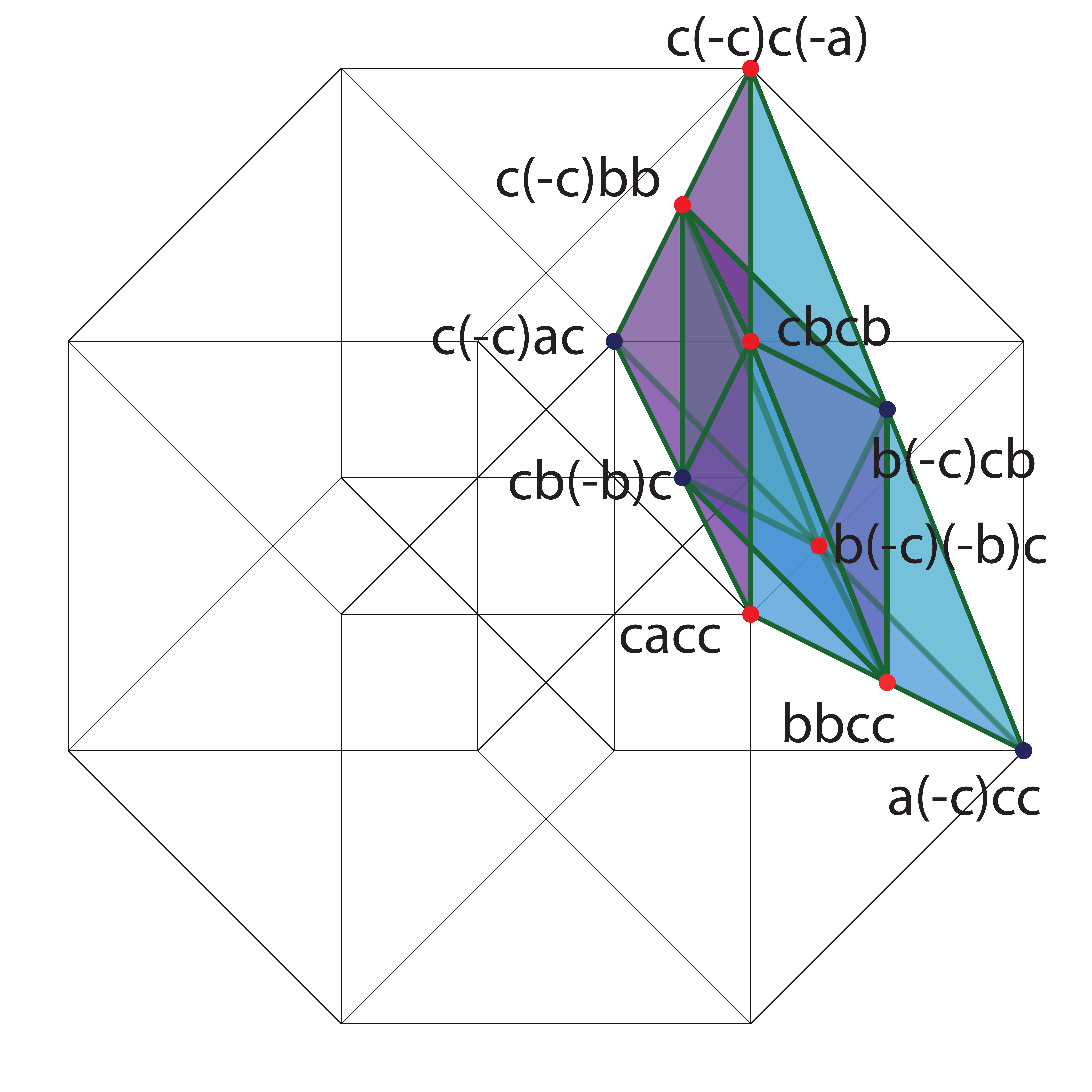}
\includegraphics[scale=.065]{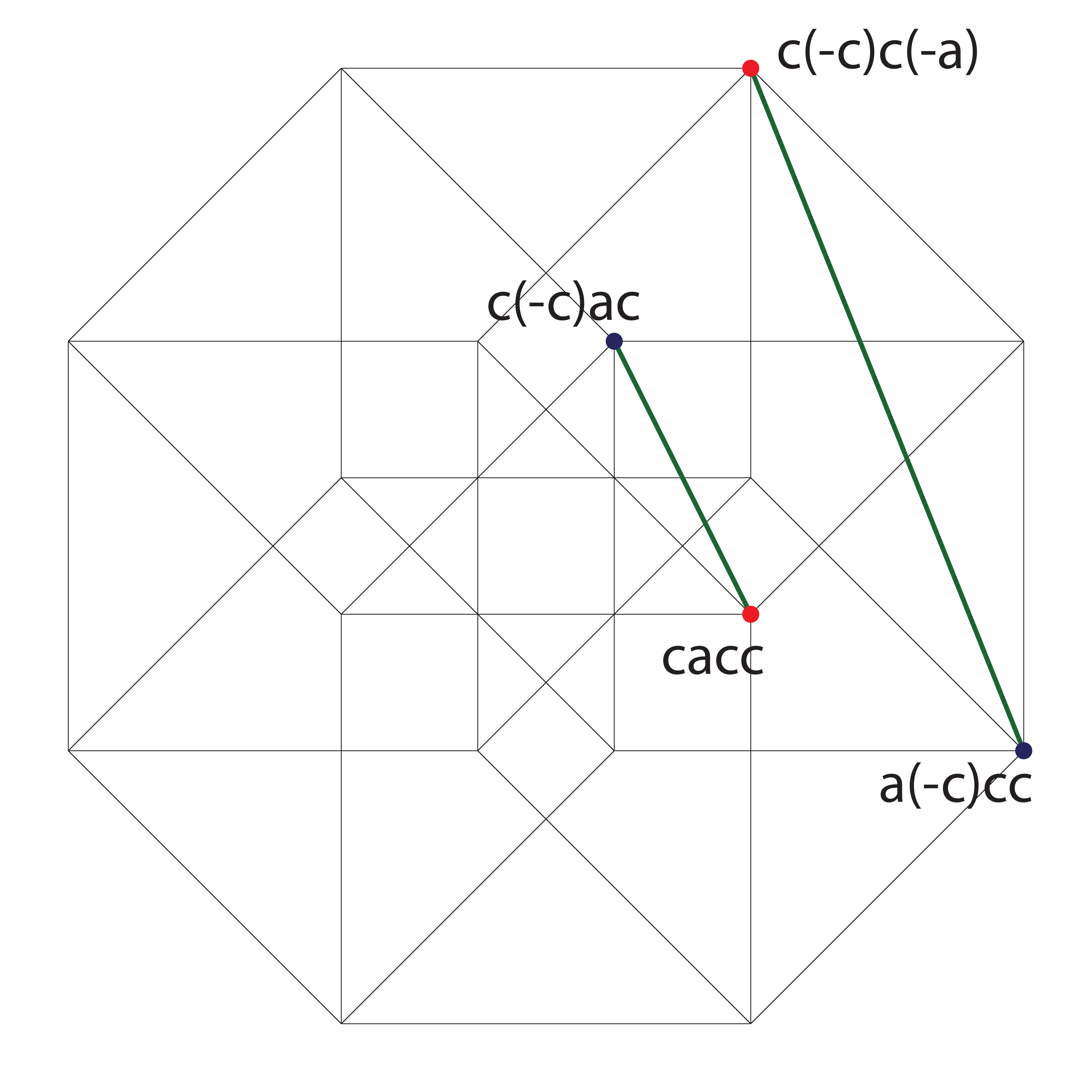}
\includegraphics[scale=.065]{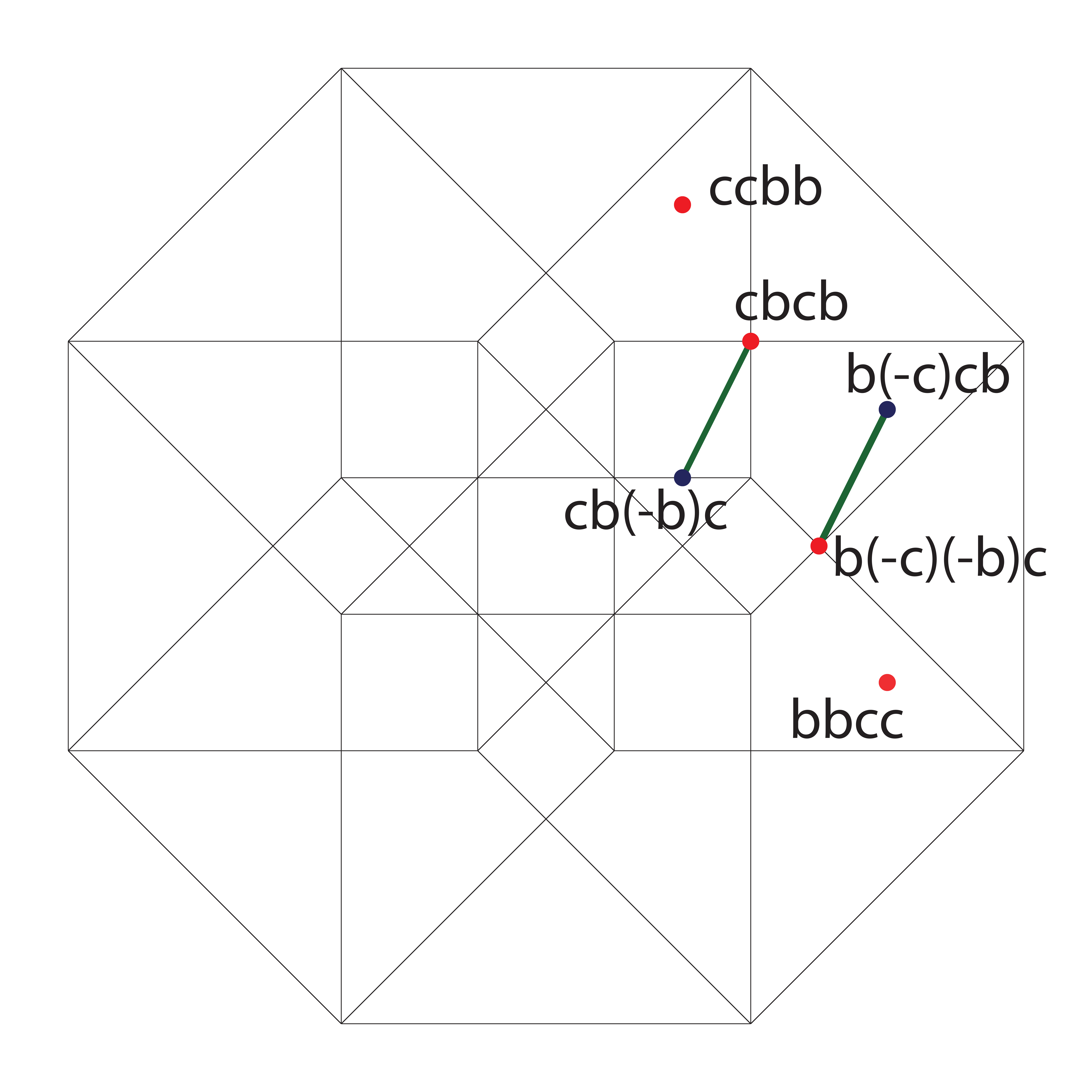}
\end{center}
\caption{The terms of the form $ac^3$ and $b^2c^2$ in the expansion of RHS} 
\label{sudoko6}
\end{figure}

\begin{figure}[htb]
   \begin{center}
\includegraphics[scale=.08]{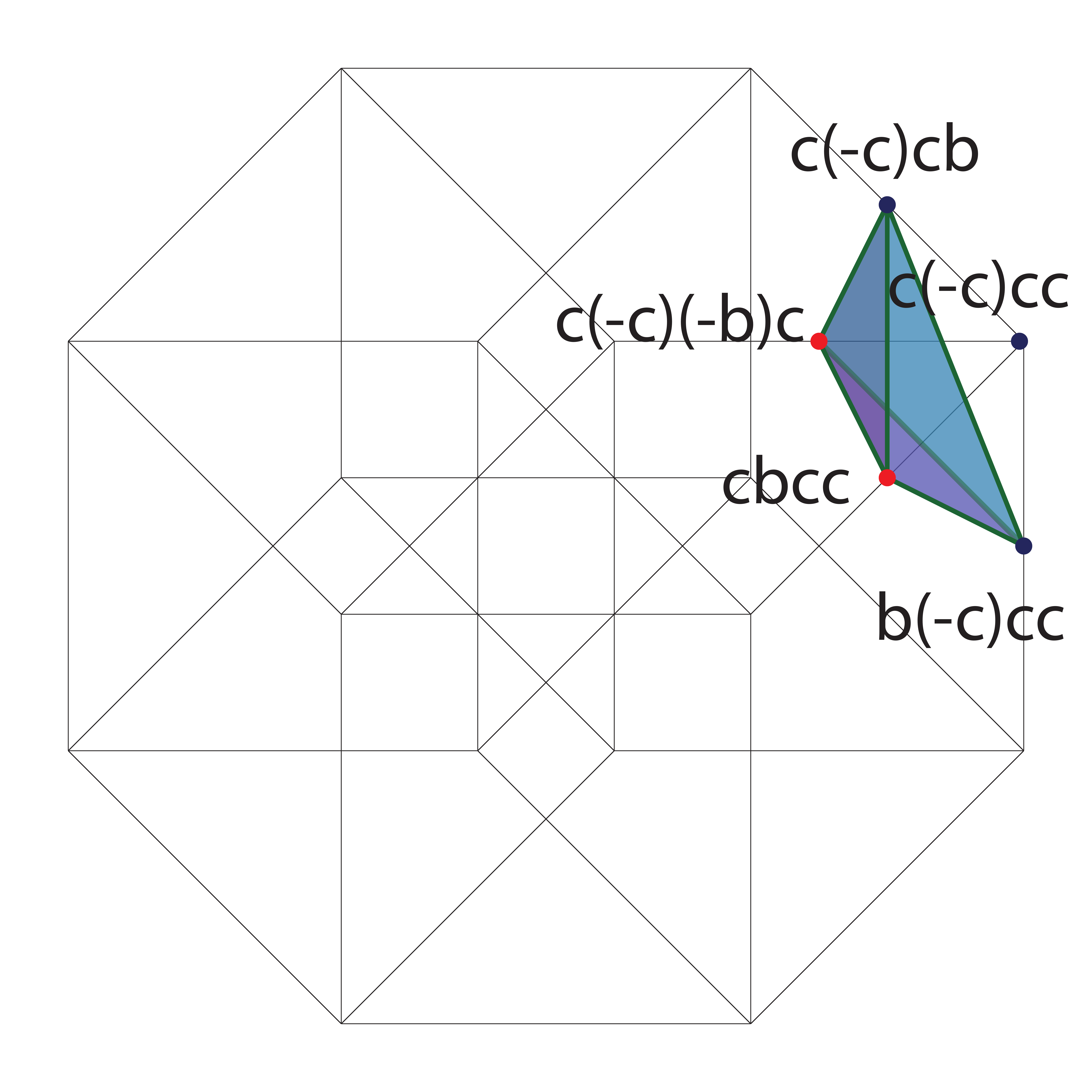}
\includegraphics[scale=.08]{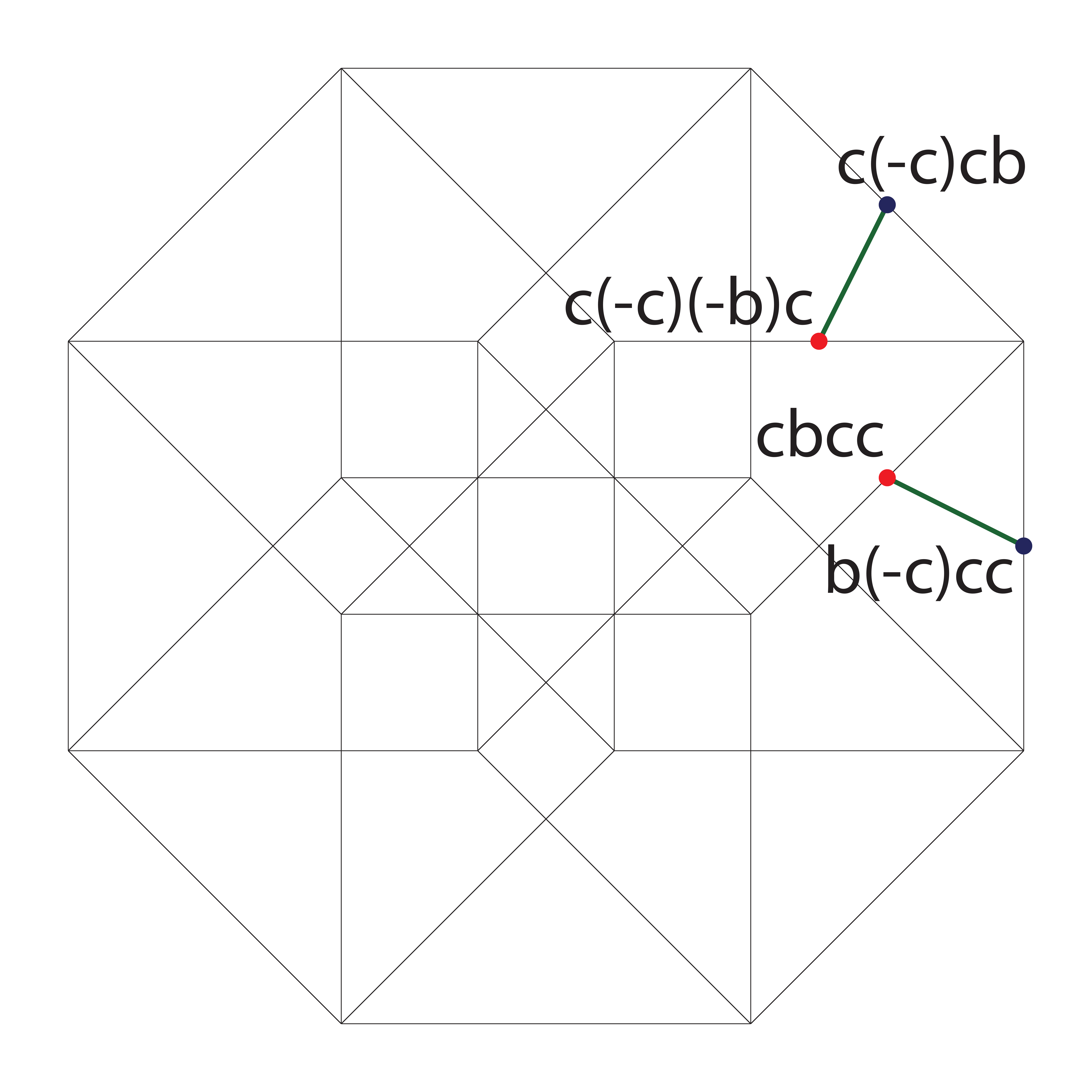}
\end{center}
\caption{The terms of the form $bc^3$ and $c^4$ in the expansion of RHS} 
\label{sudoko7}
\end{figure}


\begin{figure}[htb]
   \begin{center}
\includegraphics[scale=.1]{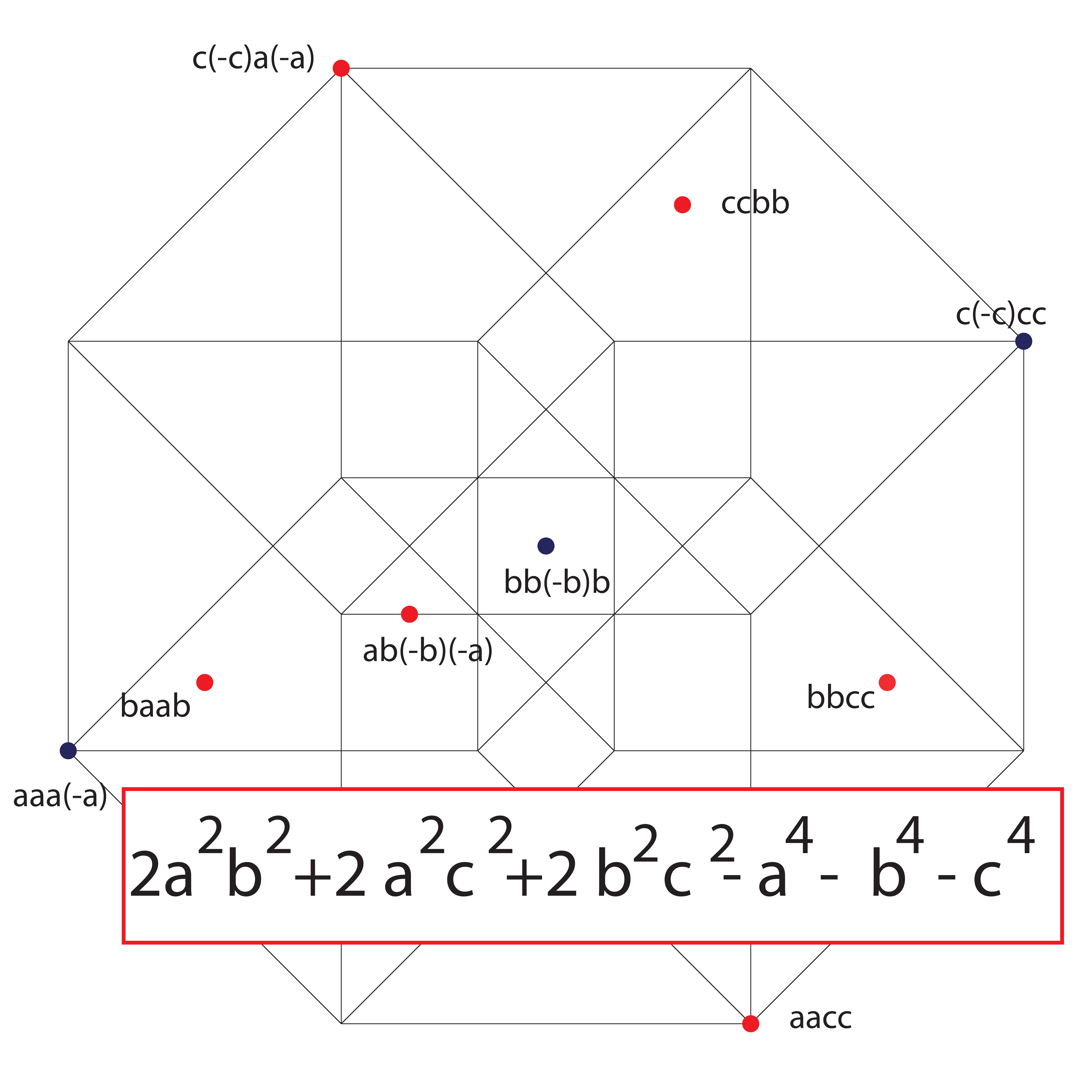}
\end{center}
\caption{The terms that remain in the expansion of RHS} 
\label{remaining}
\end{figure}

In the approach above, we kept track of each term in the expansion as a product --- the first factor in the product that represents a term was taken as a term from the first factor of the product. You may reread that sentence several times if you like, but it is better now simply to consider the product 
$$RH=(a+b+c)(a+b-c)(a-b+c)(-a+b+c).$$ 
In each term of the expansion, the first factor is positive, the second factor is negative if it is a $c$, and so forth. We will arrange the individual terms along the affine $3$-planes $x+y+z+w=K$ for $K=0,1,\ldots 8,$ and color the terms blue if they are negative and red if they are positive. Fig.~\ref{sudoko1} through Fig.~\ref{remaining} contain the illustrations. In this way, we need only consider up to $12$ terms at any one time, and only  among those terms will cancelations occur.

\section{The Pythagorean Theorem}
\label{pythag}

\begin{figure}[htb]
   \begin{center}
\includegraphics[scale=.1]{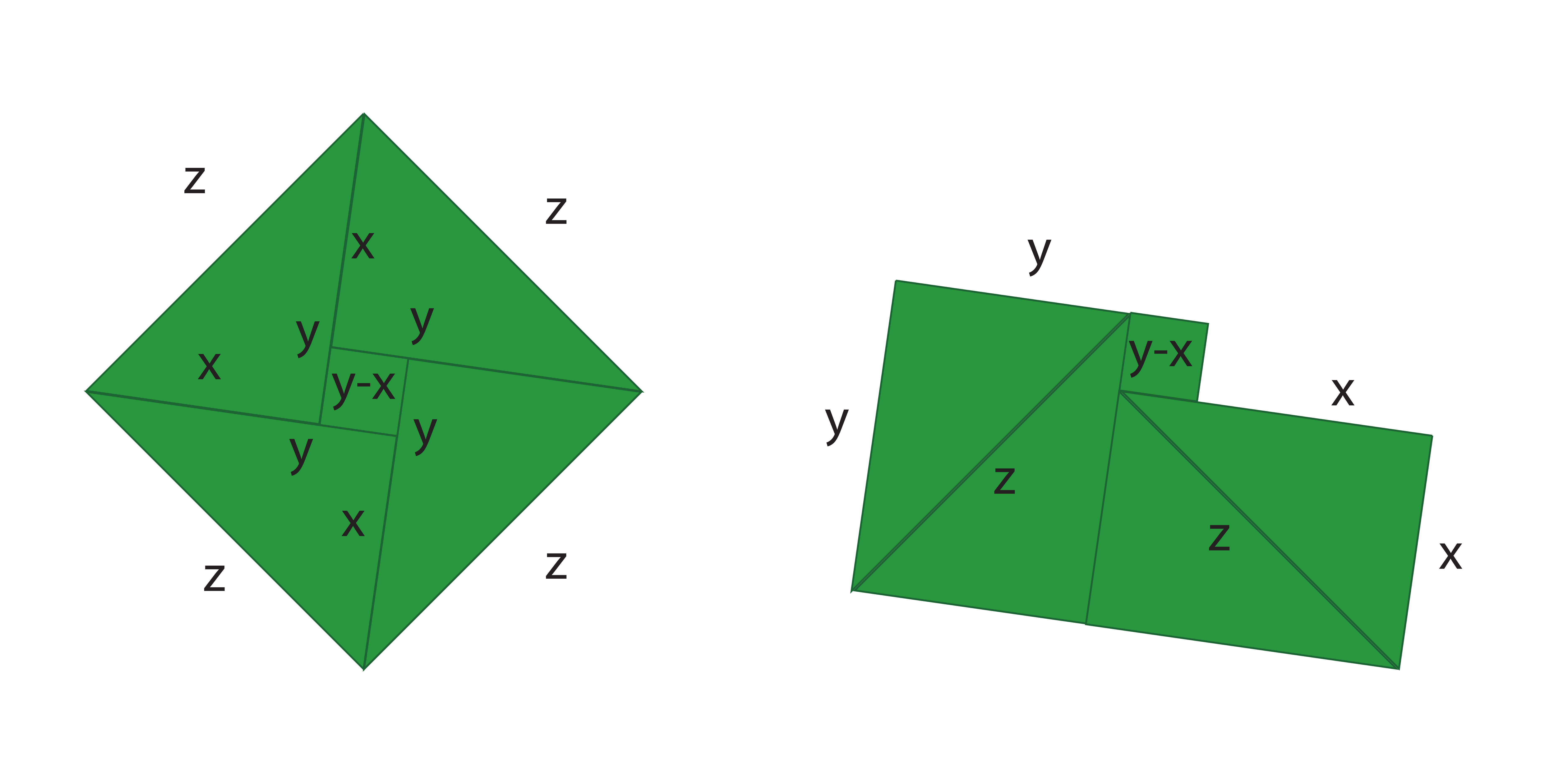}
\end{center}
\caption{The proof of the Pythagorean Theorem} 
\label{pythagproof}
\end{figure}

The illustration in Fig.~\ref{pythagproof} schematizes a proof of the Pythagorean Theorem. Four right triangles with legs of length $x$ and $y$ and hypotenuse $z$ are arranged to nearly fill a square whose edge length is $z$. The missing piece is a square of size $(y-x)\times(y-x)$, where, by convention, $x \le y \le z$. (In the case for which $x=y$, the central square is empty.) Let us call the copies of the triangles $NE$, $SE$, $SW$, and $NW$ depending on their positions with respect to the cardinal directions. The $NE$ and $SW$ triangles are reassembled to form a rectangle of size $x \times y$; similarly the $NW$ and $SE$ triangles form a congruent rectangle. In the reassembly, the triangles are only moved by translation. 
Then the central square is appended. The union on the central square and the $NW$ to $SE$ rectangle together with a portion of remaining rectangle form square of size $y \times y$. The remaining piece of the other rectangle is a square of size $x \times x$. Thus the area $z \times z$ can be cut into two squares of size $x\times x$ and $y \times y$.

\begin{figure}[htb]
\includegraphics[scale=.08]{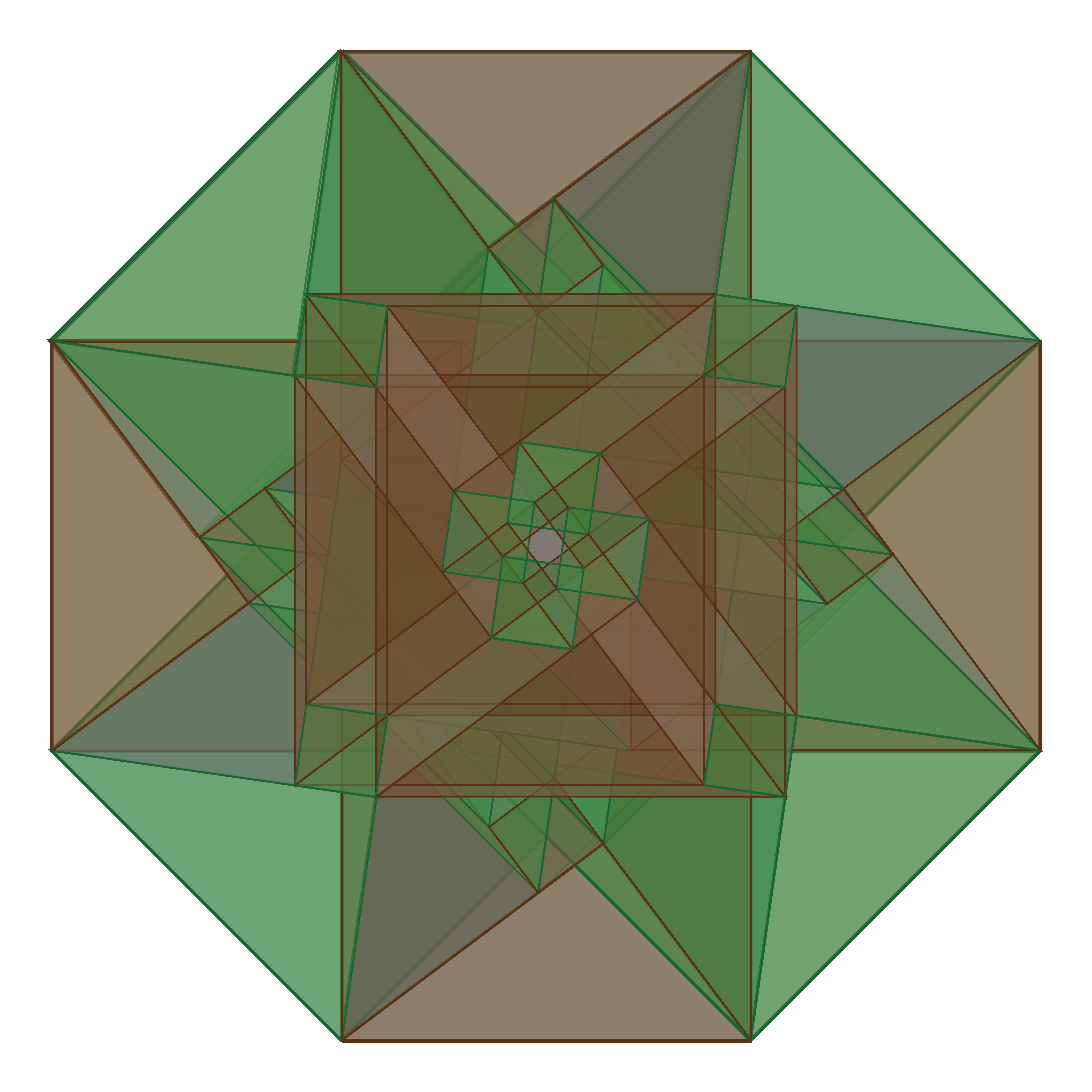}
\includegraphics[scale=.08]{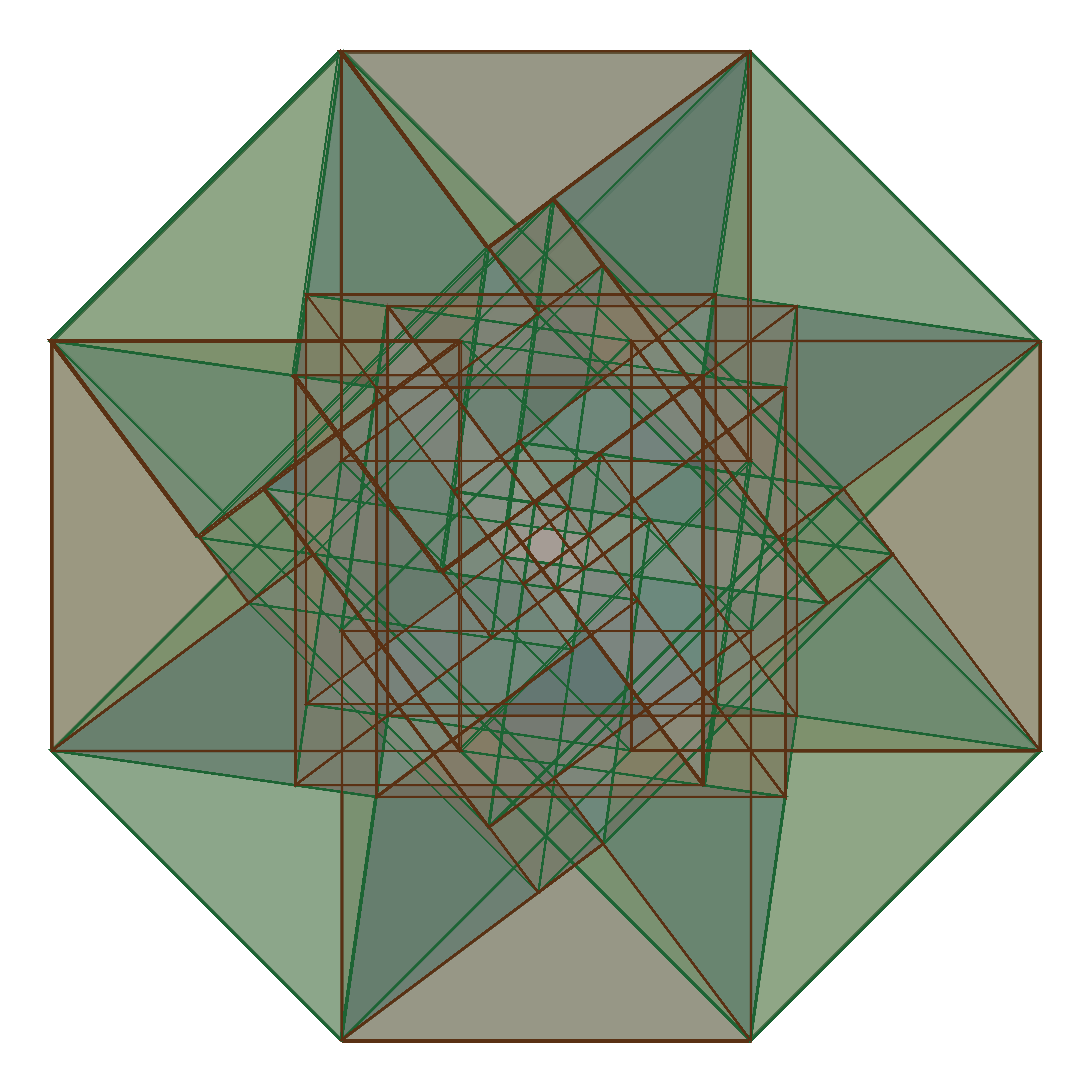}
\caption{Two views of the hyper-rectangle  $R_{zzww}$} 
\label{stereo}
\end{figure}

Now consider a pair of right triangles. The first has legs of length $x$ and $y$ with $x\le y$. The second has legs of length $u$ and $v$ with $u\le v$. The length of the hypotenuses are $z$ and $w$ respectively. We'll now consider a hyper-rectangle $R_{zzww}=[0,z]\times [0.z] \times [0,w]\times [0,w],$ and decompose both the $(z\times z)$ and $(w\times w)$ squares into five pieces (four triangles and one small square) as in the proof of the Pythagorean theorem. The hyper-rectangle $R_{zzww}$ is cut into $25$ pieces. 
Two views of this object are illustrated in Fig.~\ref{stereo}. We encourage the reader to attempt a cross-eyed stereo-opsis to imagine some of the depth in this illustration.

\begin{figure}[htb]
 \begin{center}
\includegraphics[scale=.05]{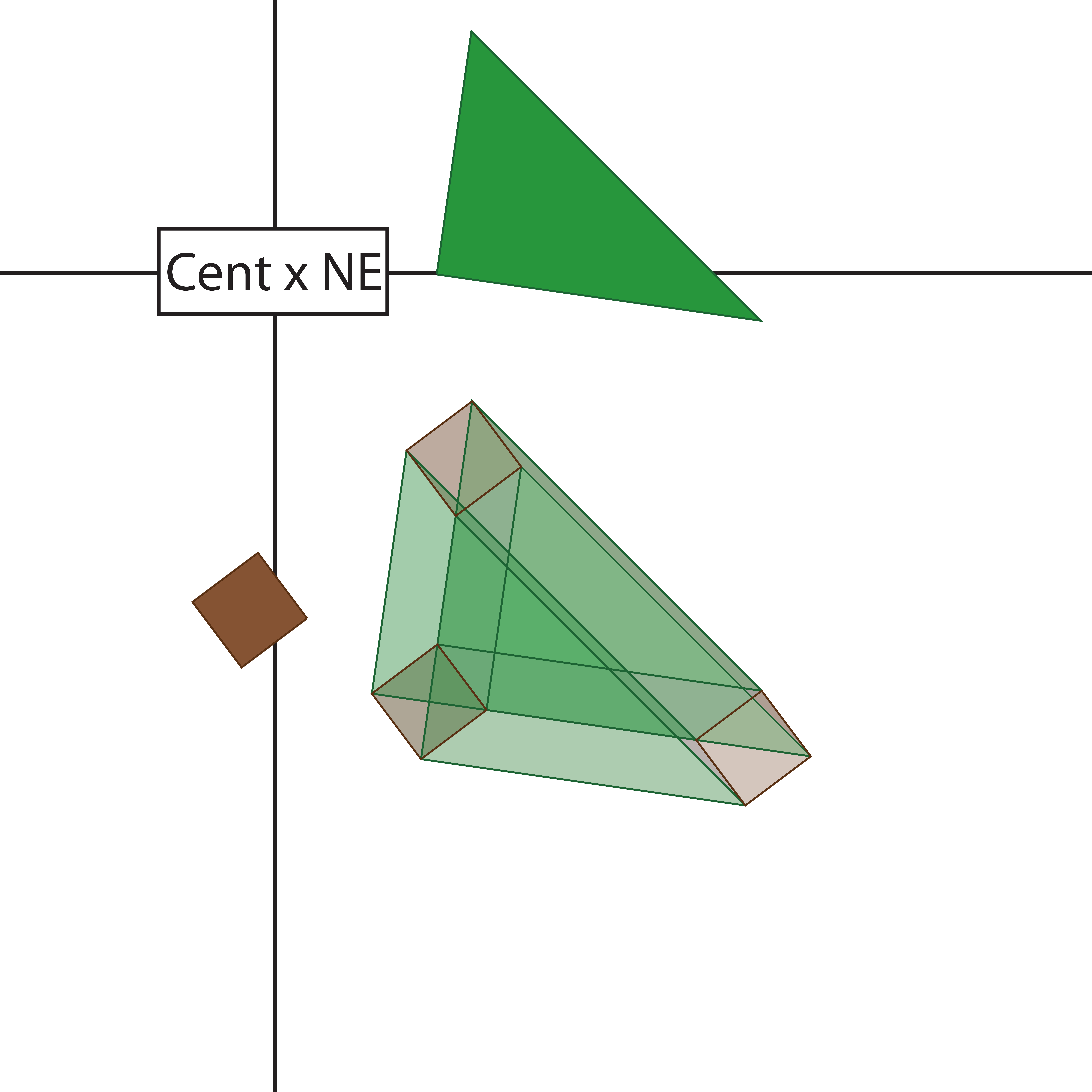}
\includegraphics[scale=.05]{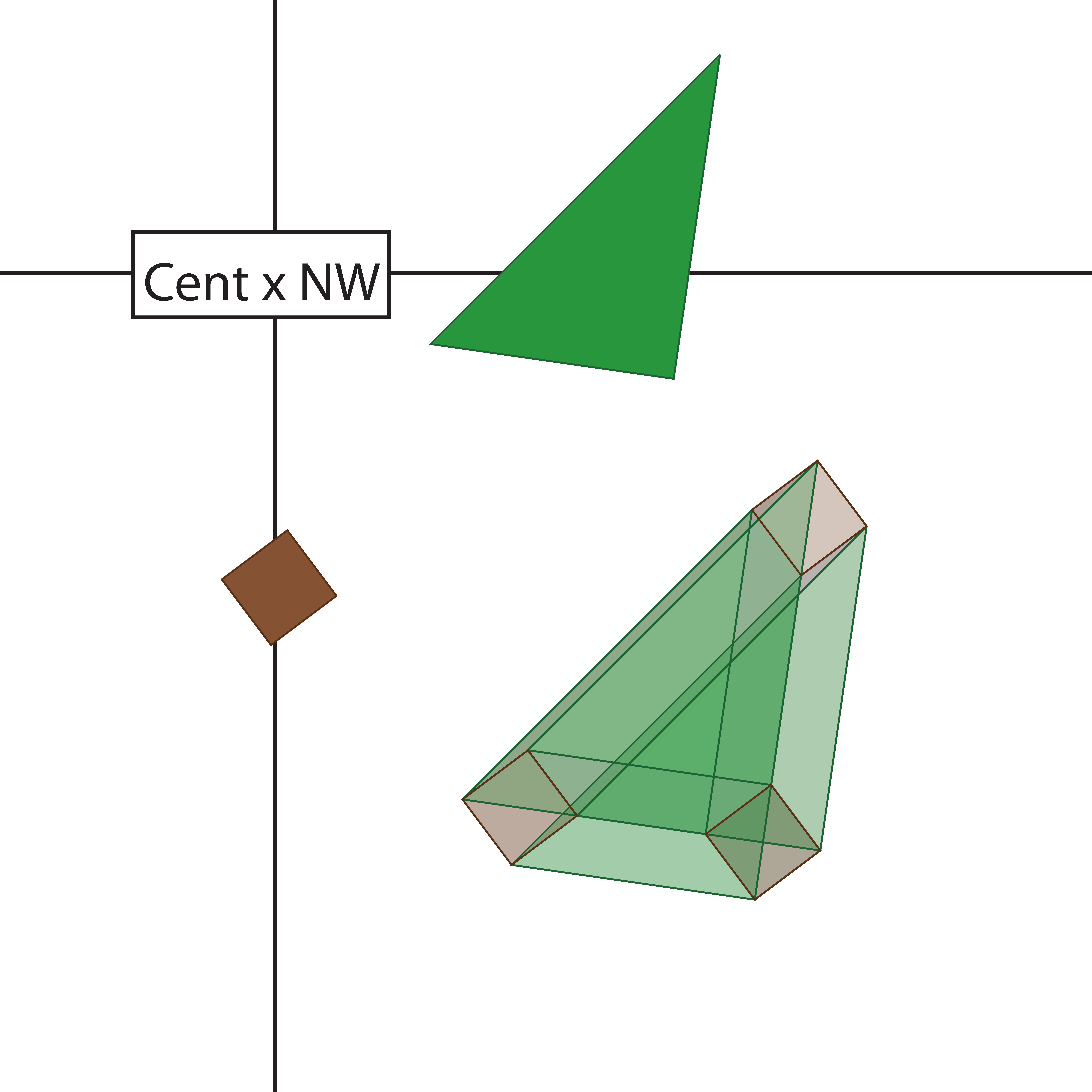}
\includegraphics[scale=.05]{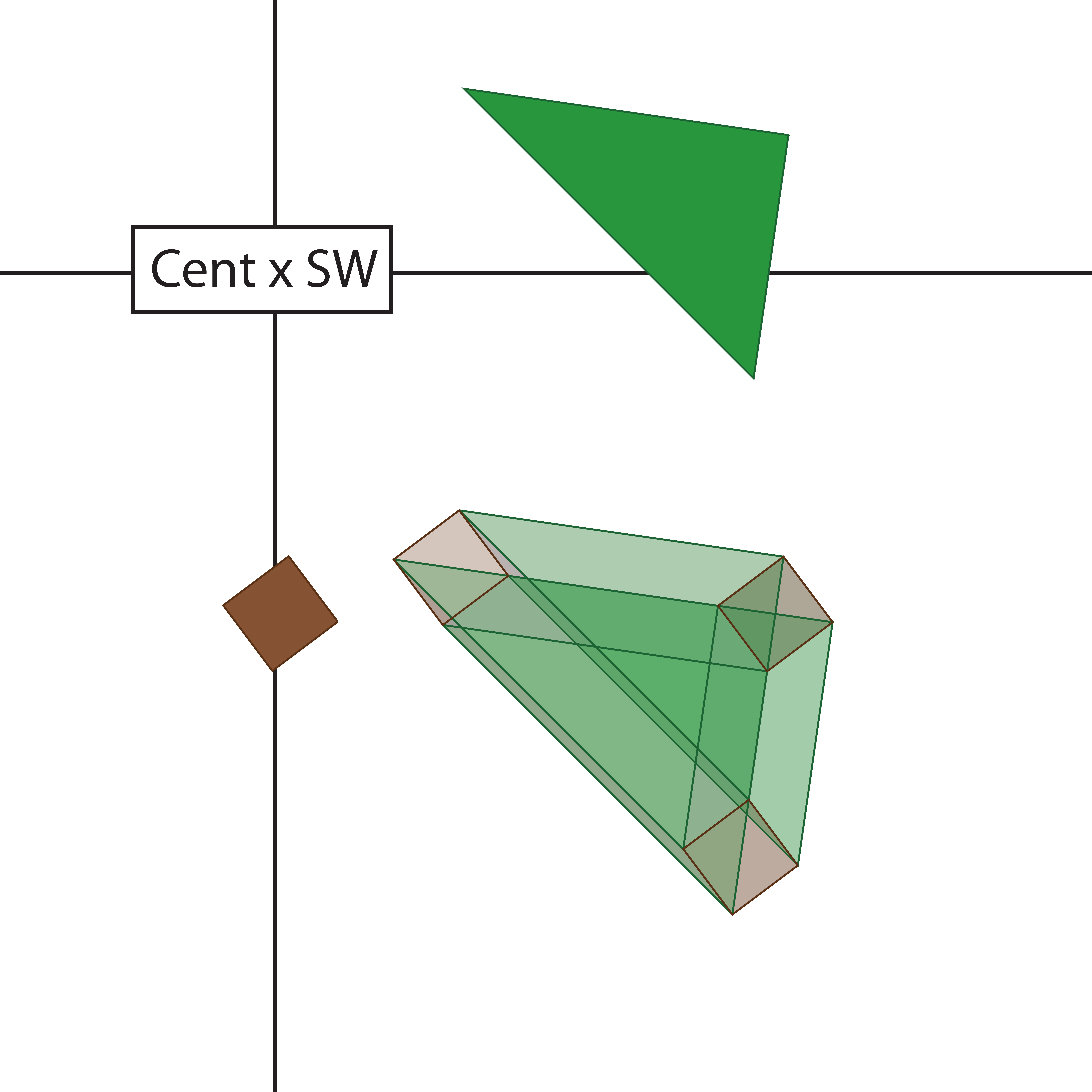}
\includegraphics[scale=.05]{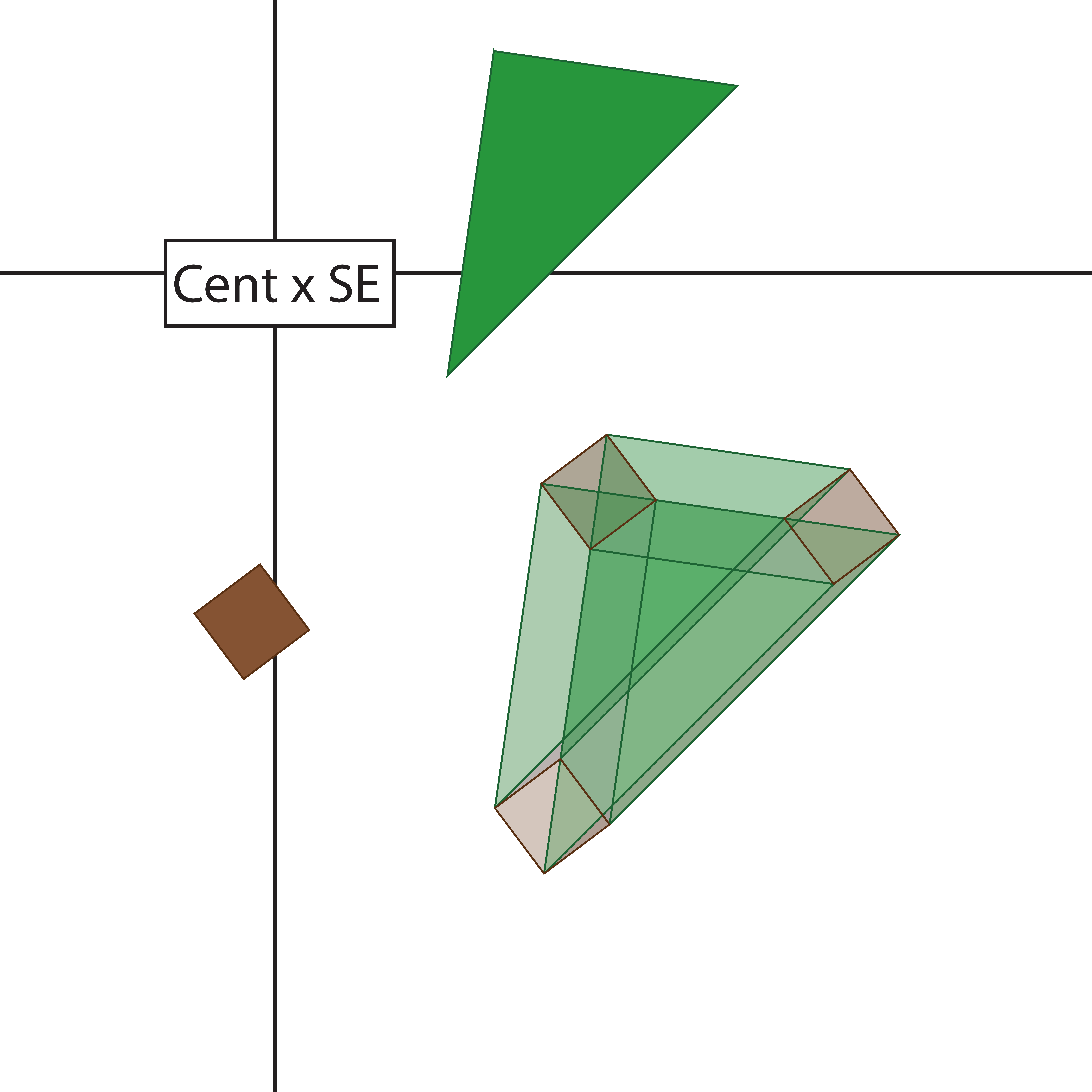}
\includegraphics[scale=.05]{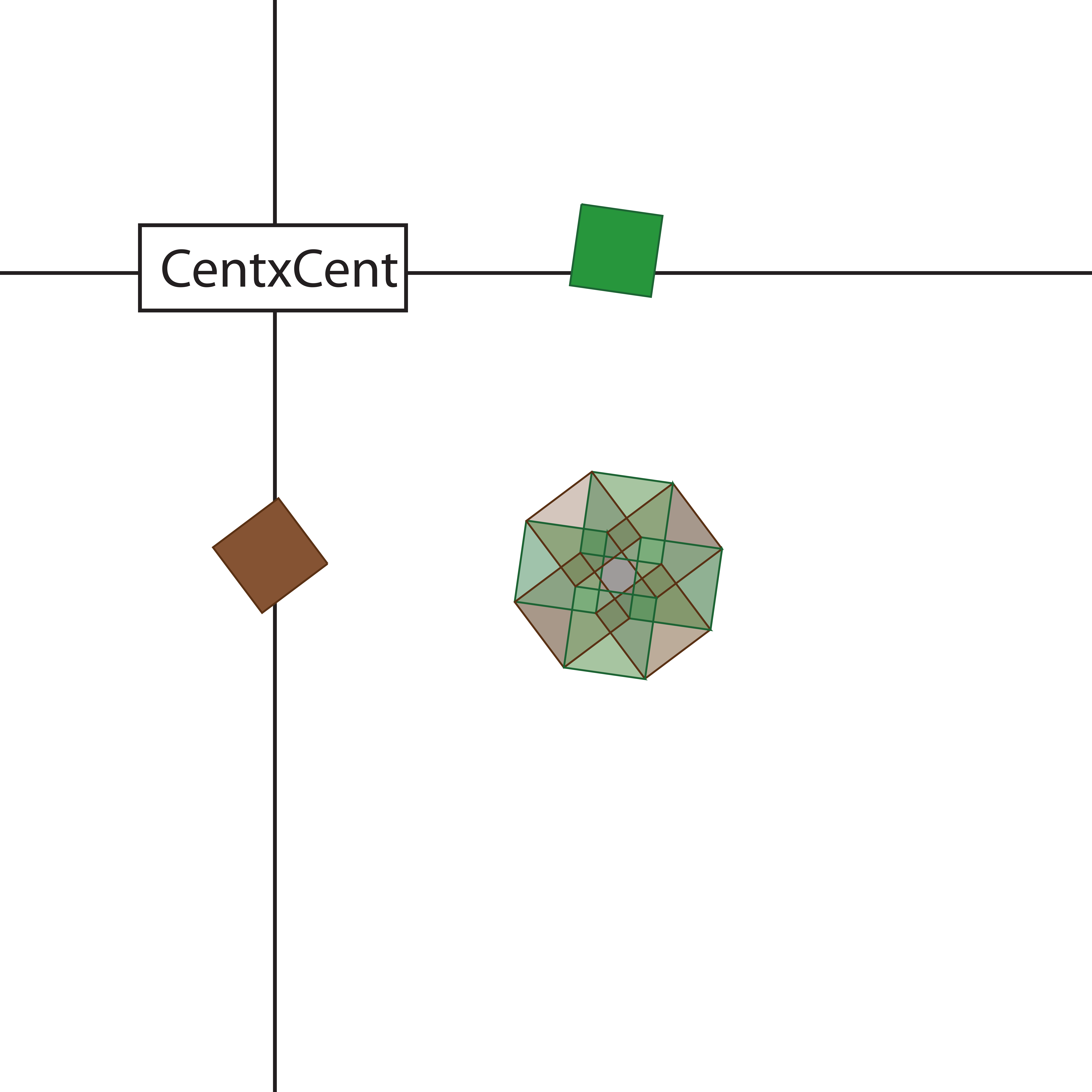}
\end{center}
\caption{The central pieces of the decomposition of the hyper-rectangle} 
\label{Central}
\end{figure}

\begin{figure}[htb]
 \begin{center}
\includegraphics[scale=.05]{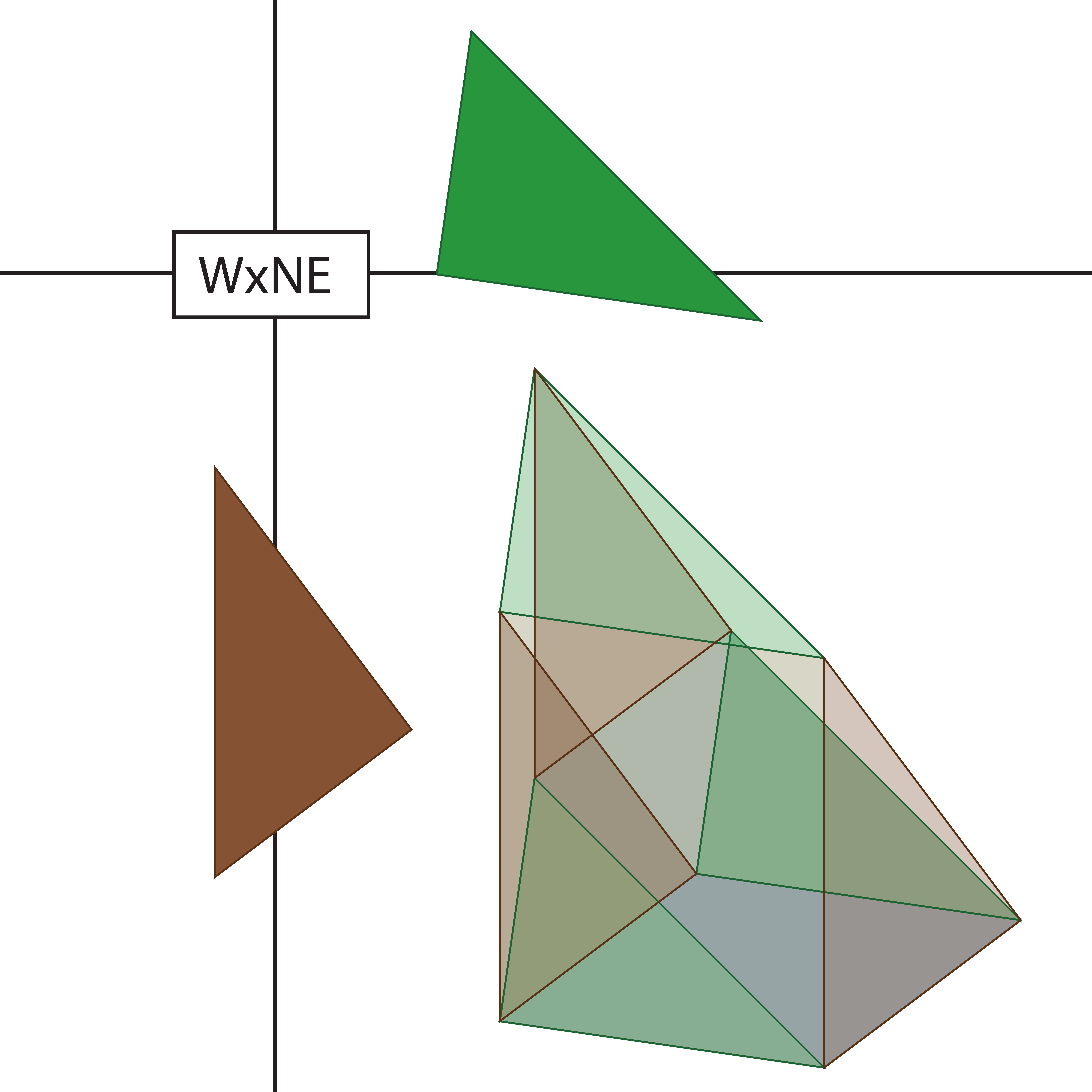}
\includegraphics[scale=.05]{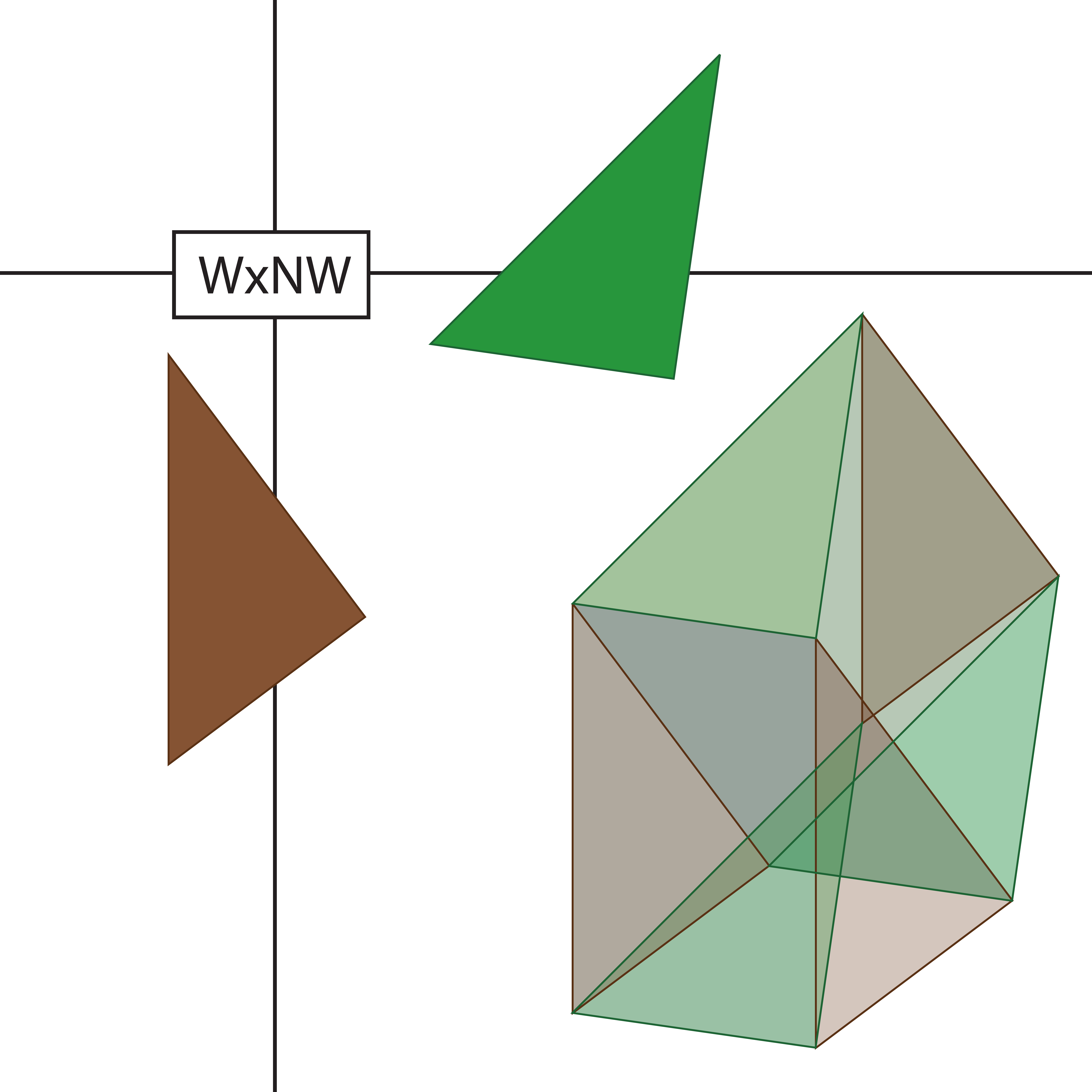}
\includegraphics[scale=.05]{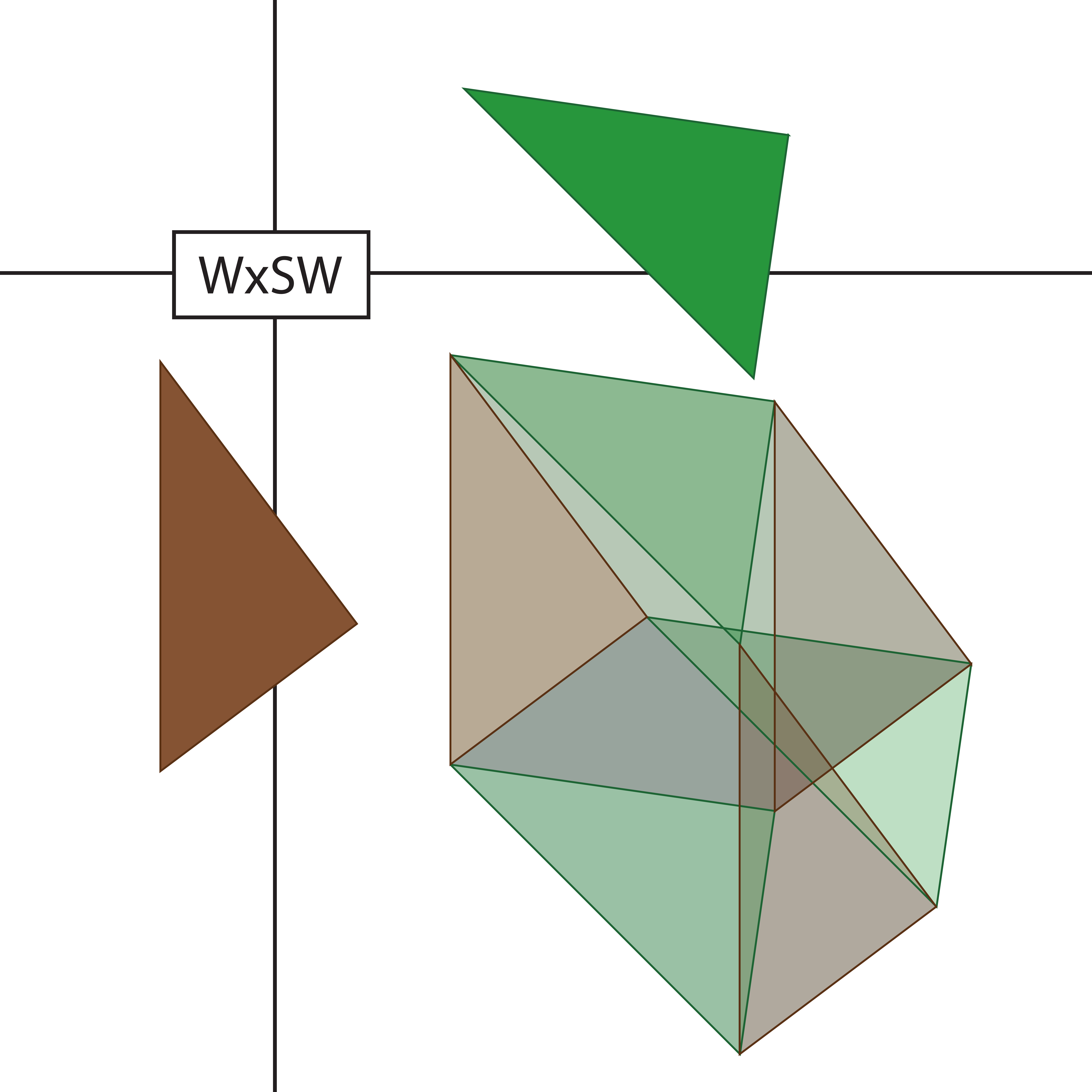}
\includegraphics[scale=.05]{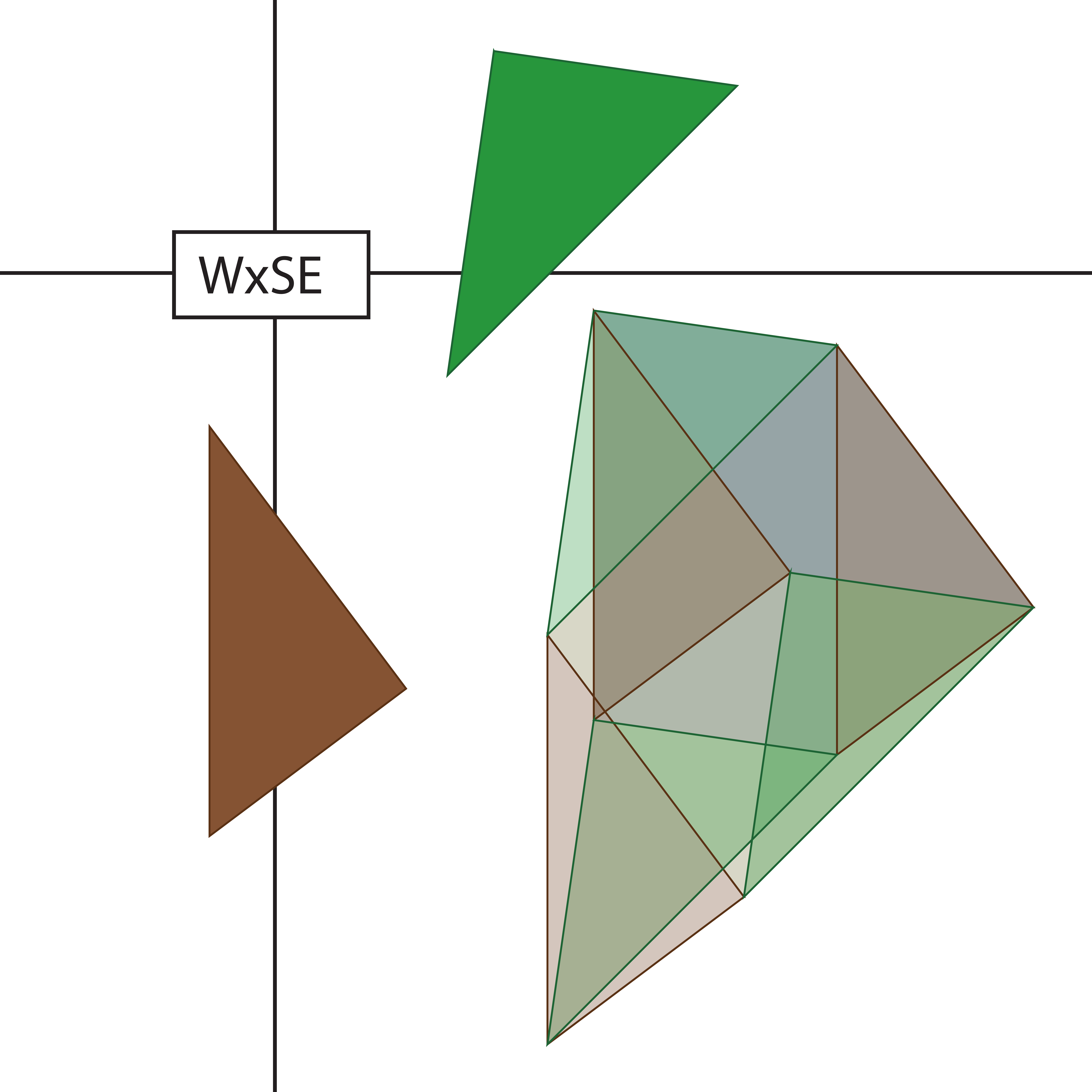}
\includegraphics[scale=.05]{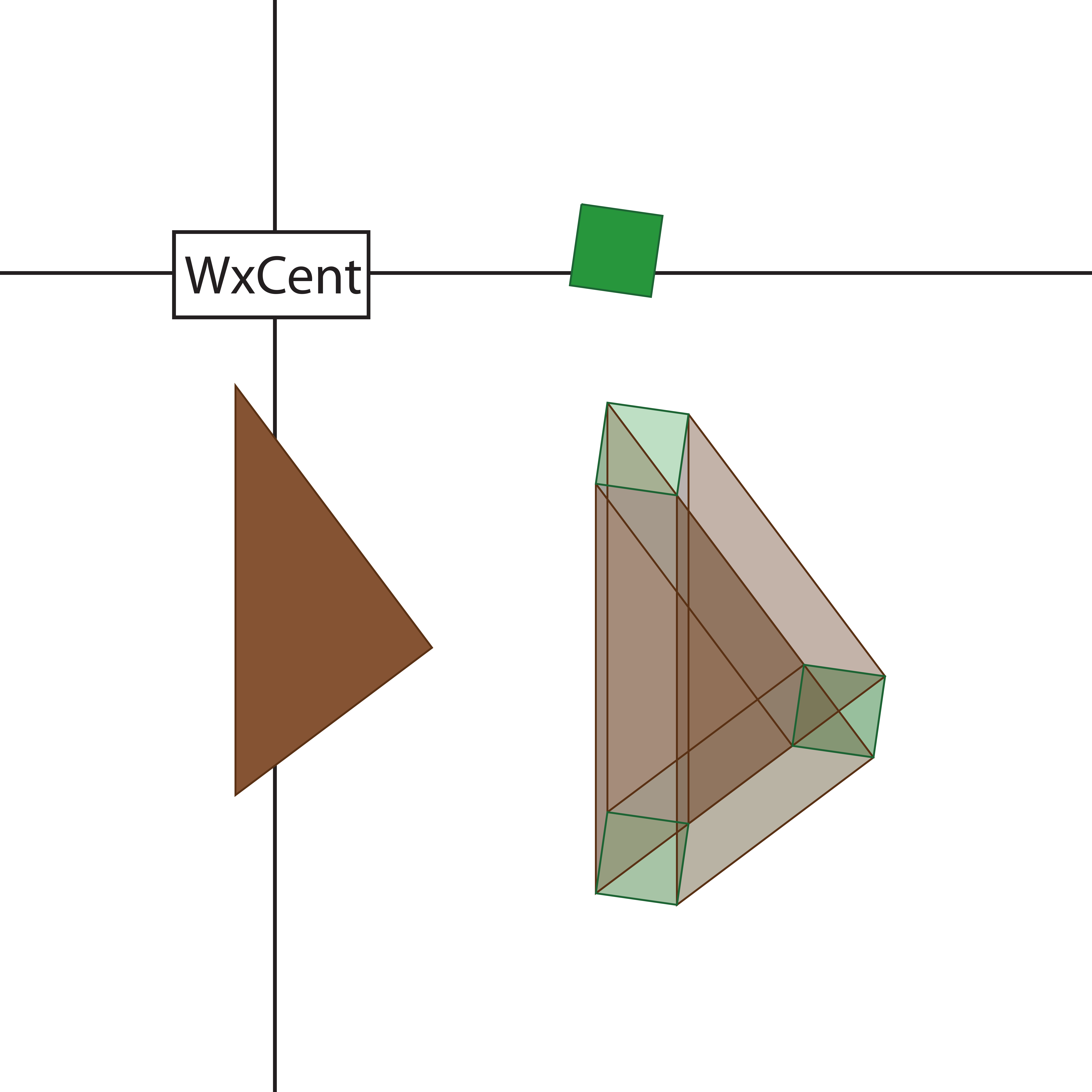}
\end{center}
\caption{The western pieces of the decomposition of the hyper-rectangle} 
\label{West}
\end{figure}

\begin{figure}[htb]
 \begin{center}
\includegraphics[scale=.05]{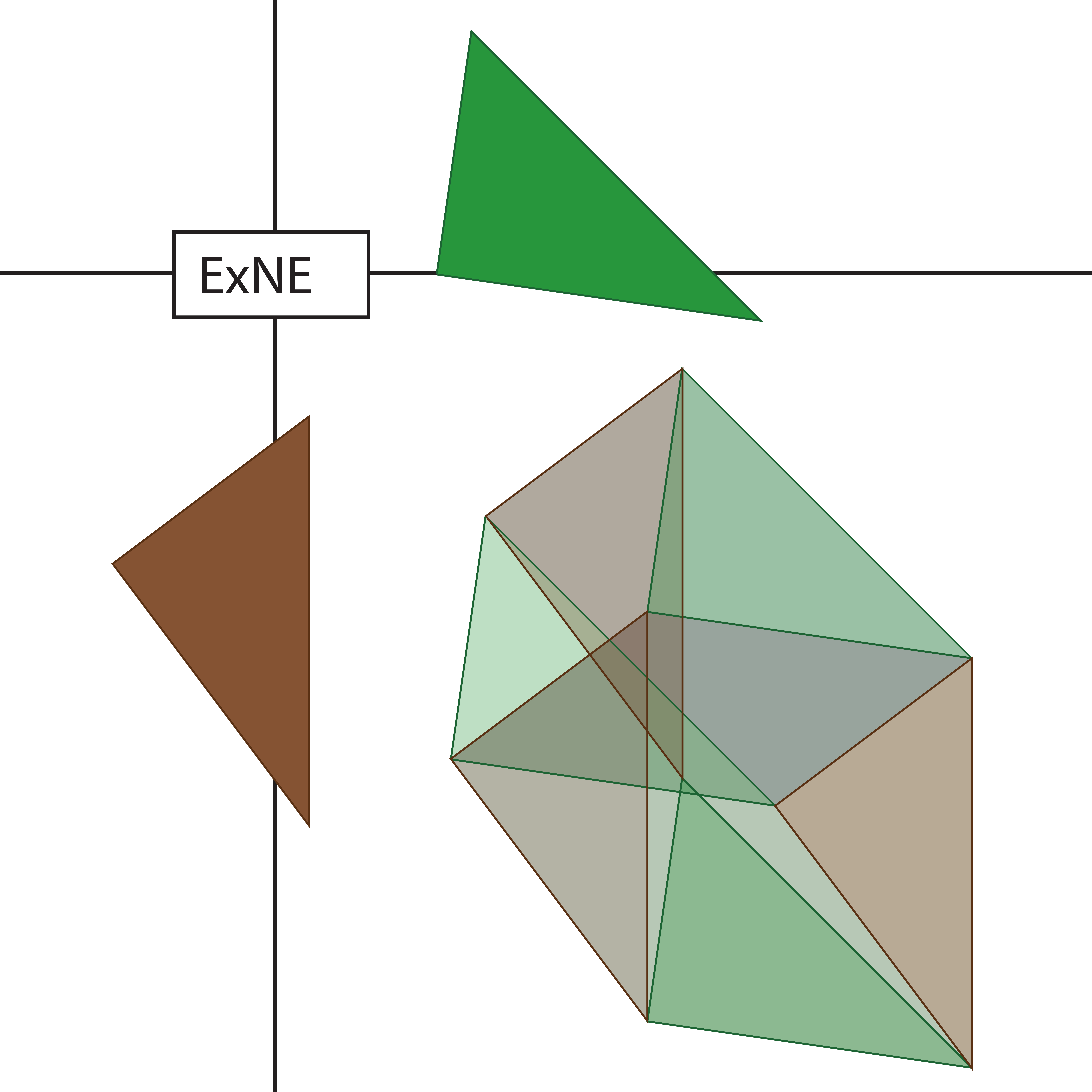}
\includegraphics[scale=.05]{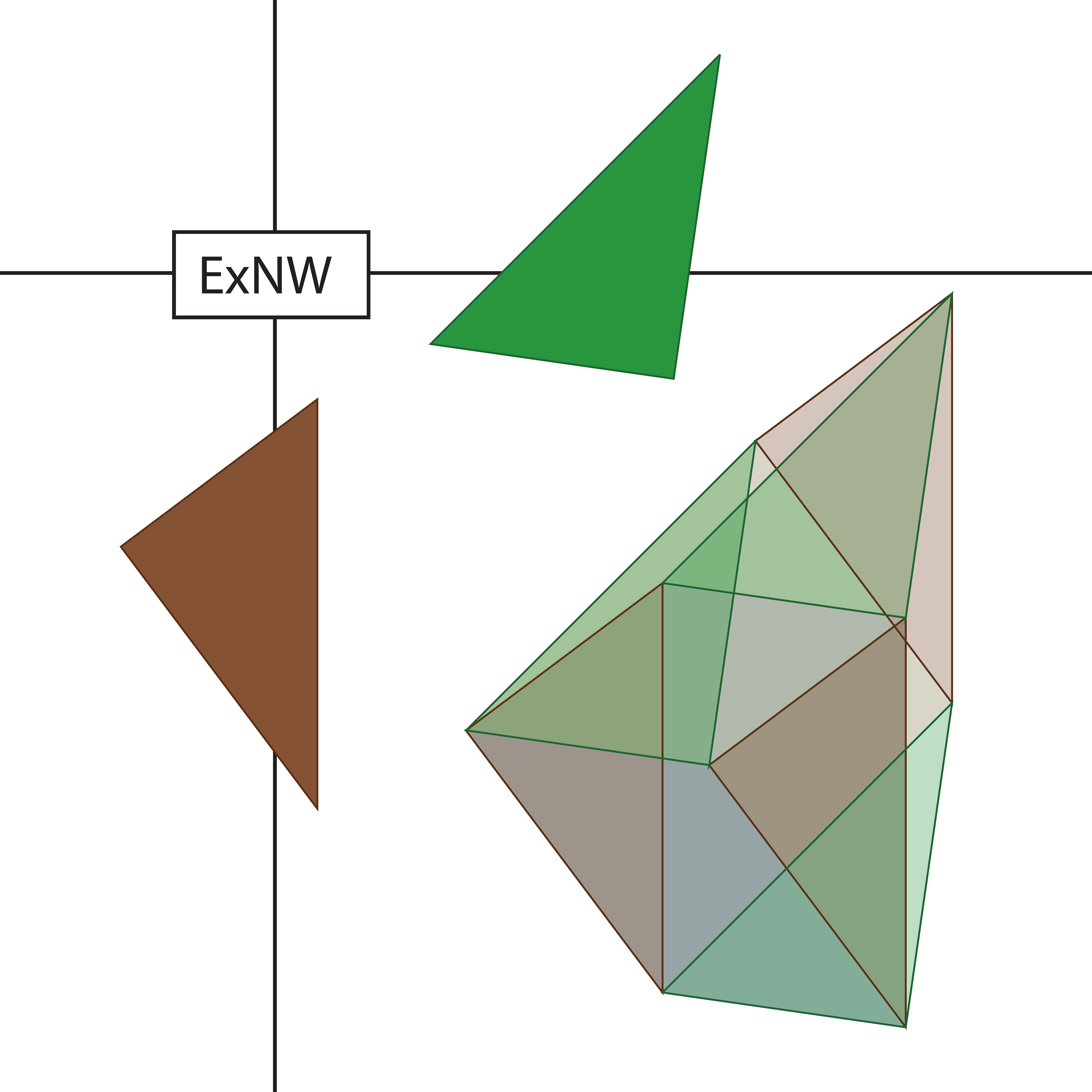}
\includegraphics[scale=.05]{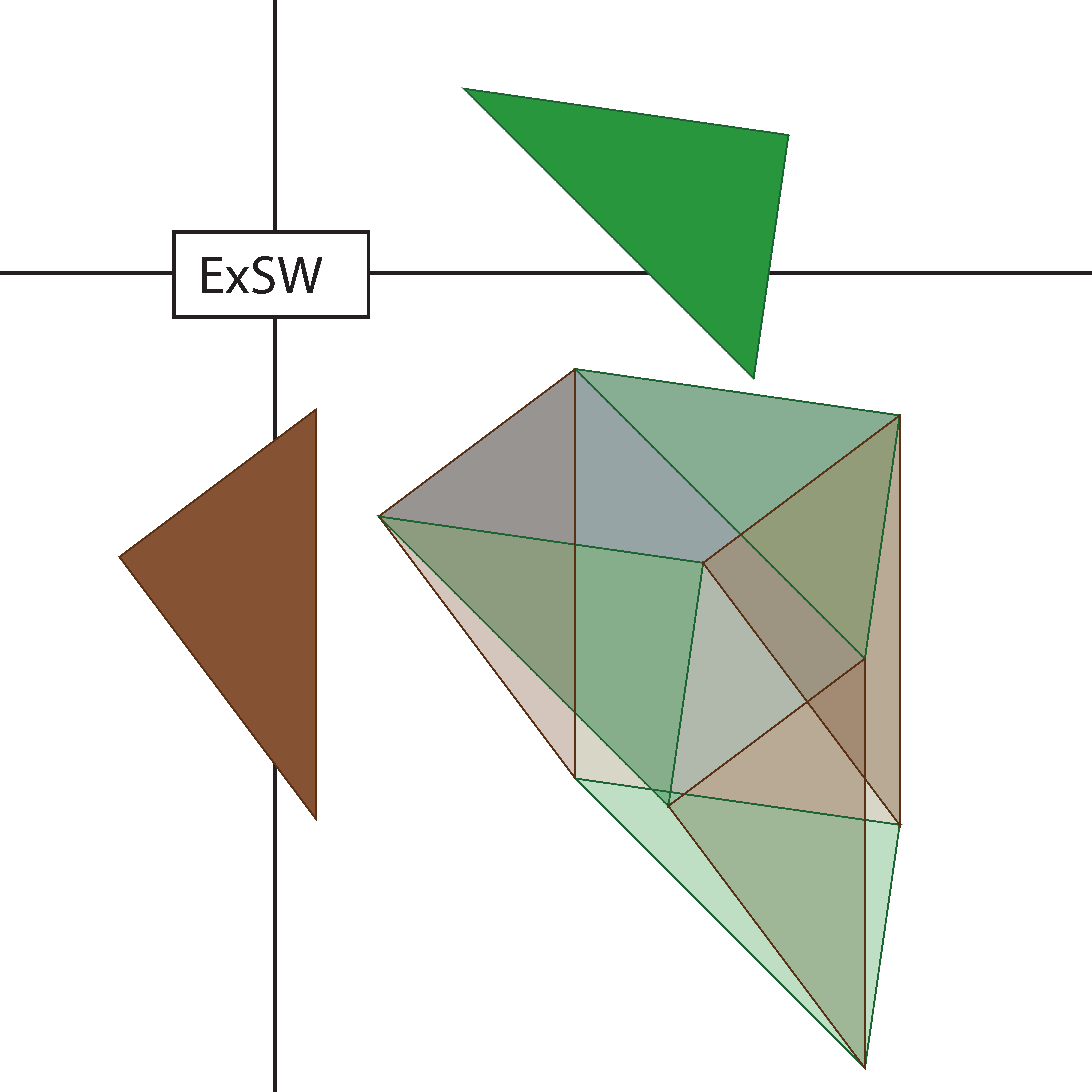}
\includegraphics[scale=.05]{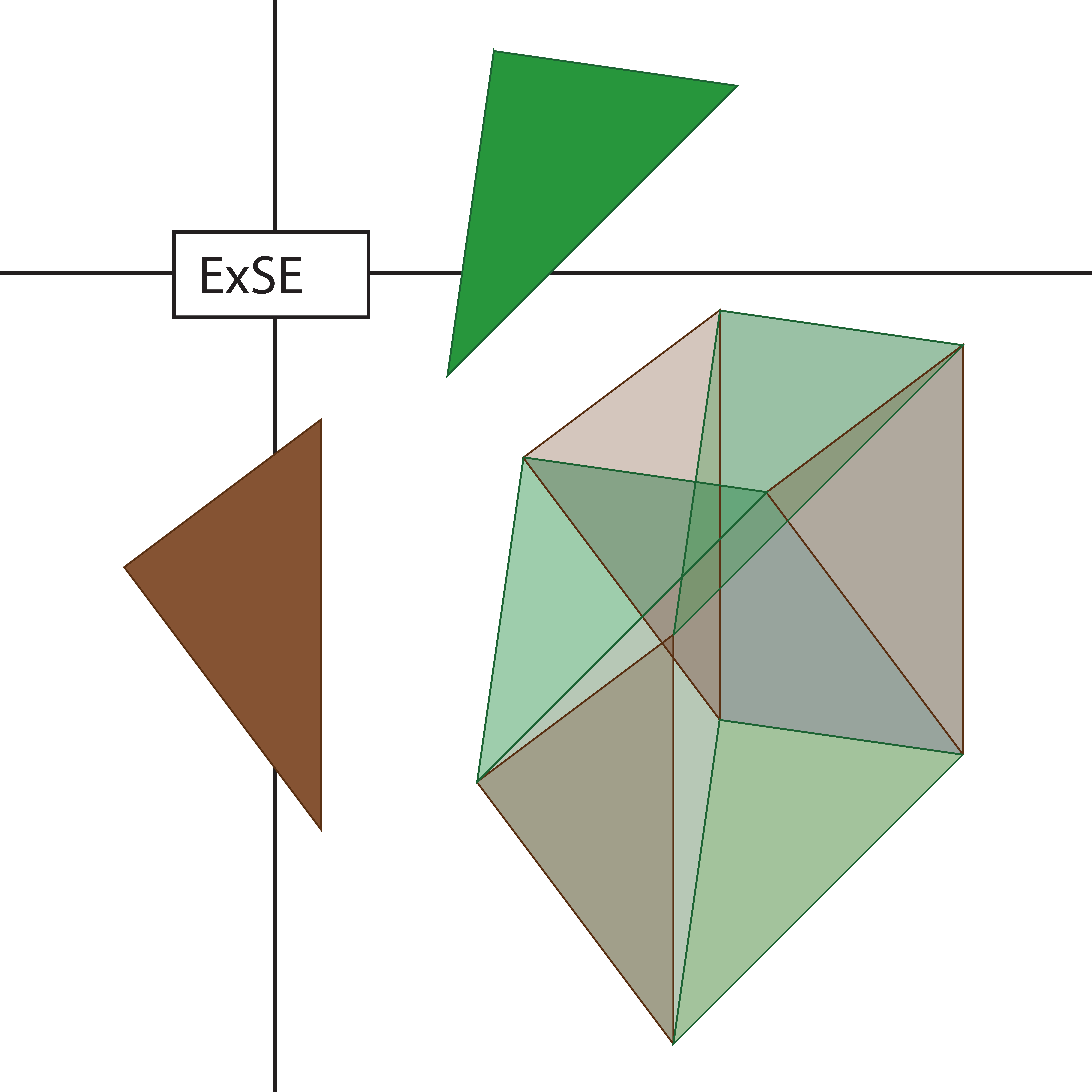}
\includegraphics[scale=.05]{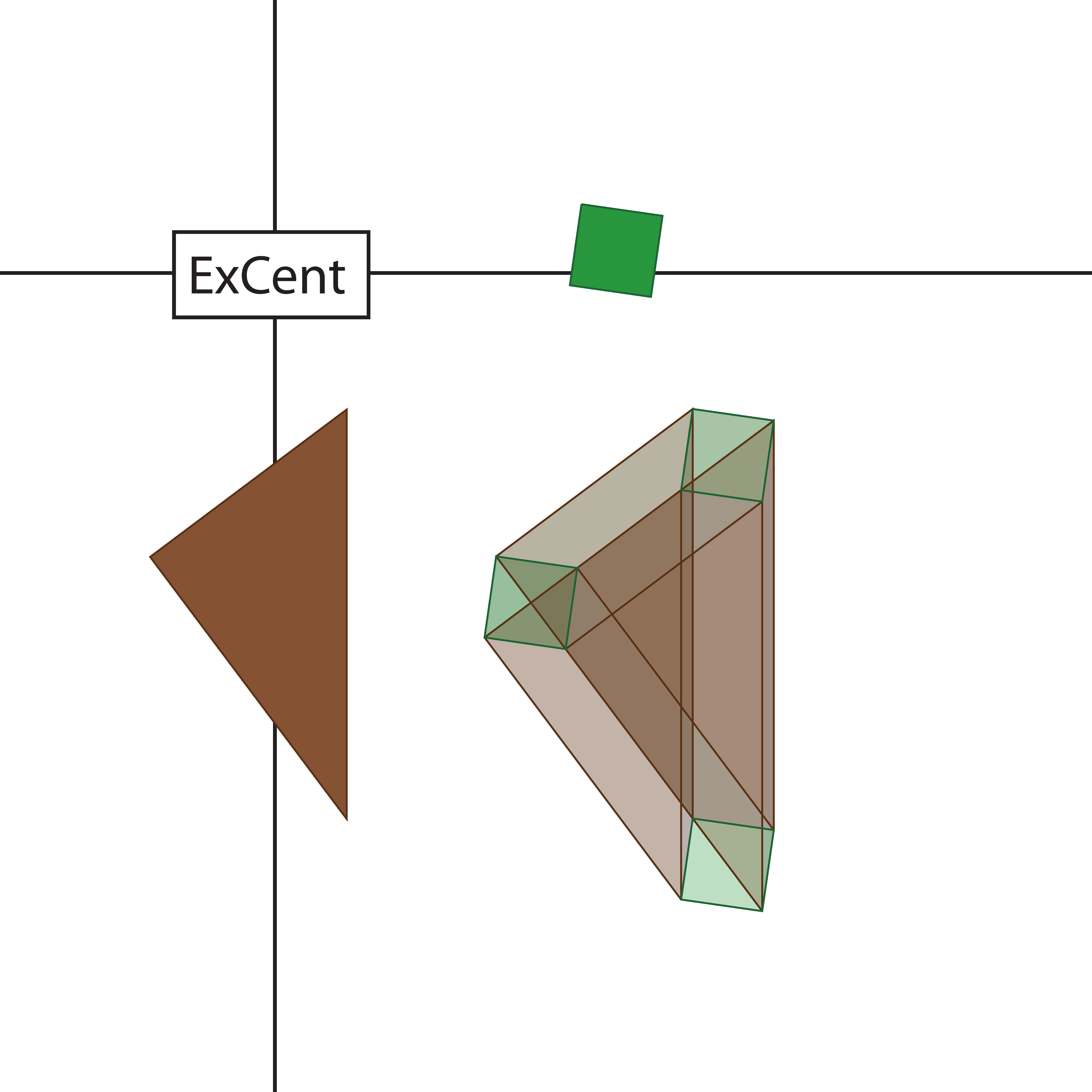}
\end{center}
\caption{The Eastern pieces of the  decomposition of the hyper-rectangle} 
\label{East}
\end{figure}

\begin{figure}[htb]
 \begin{center}
\includegraphics[scale=.05]{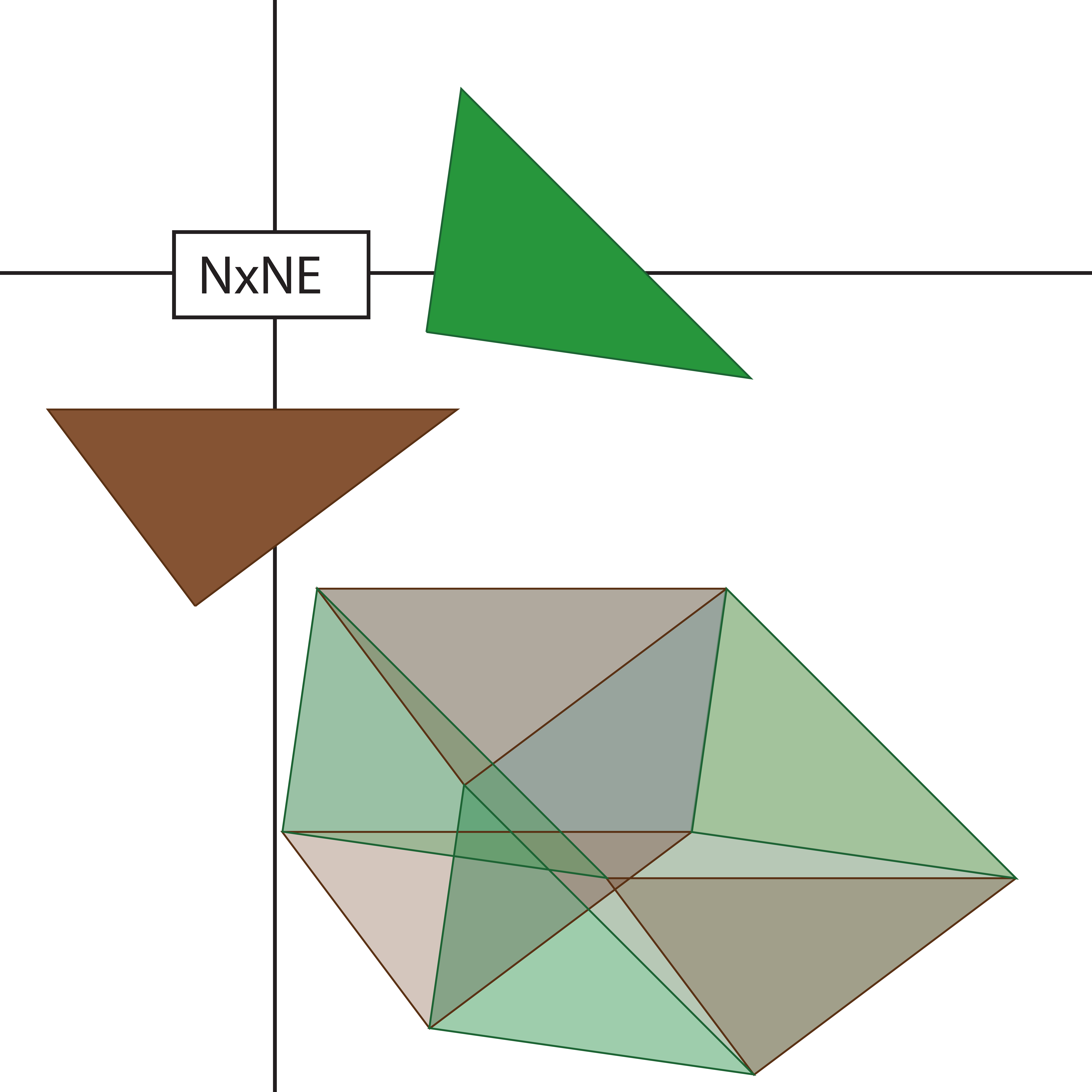}
\includegraphics[scale=.05]{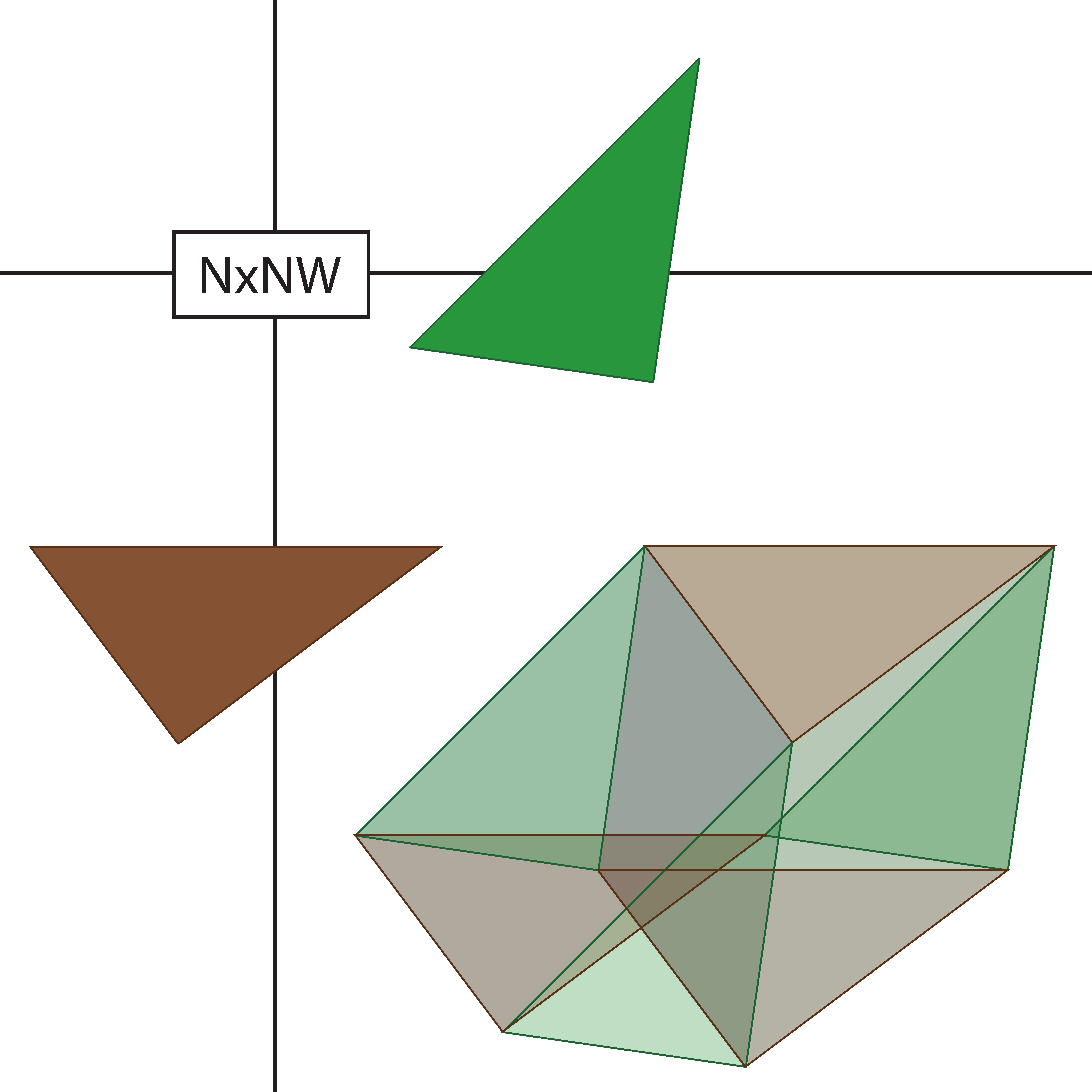}
\includegraphics[scale=.05]{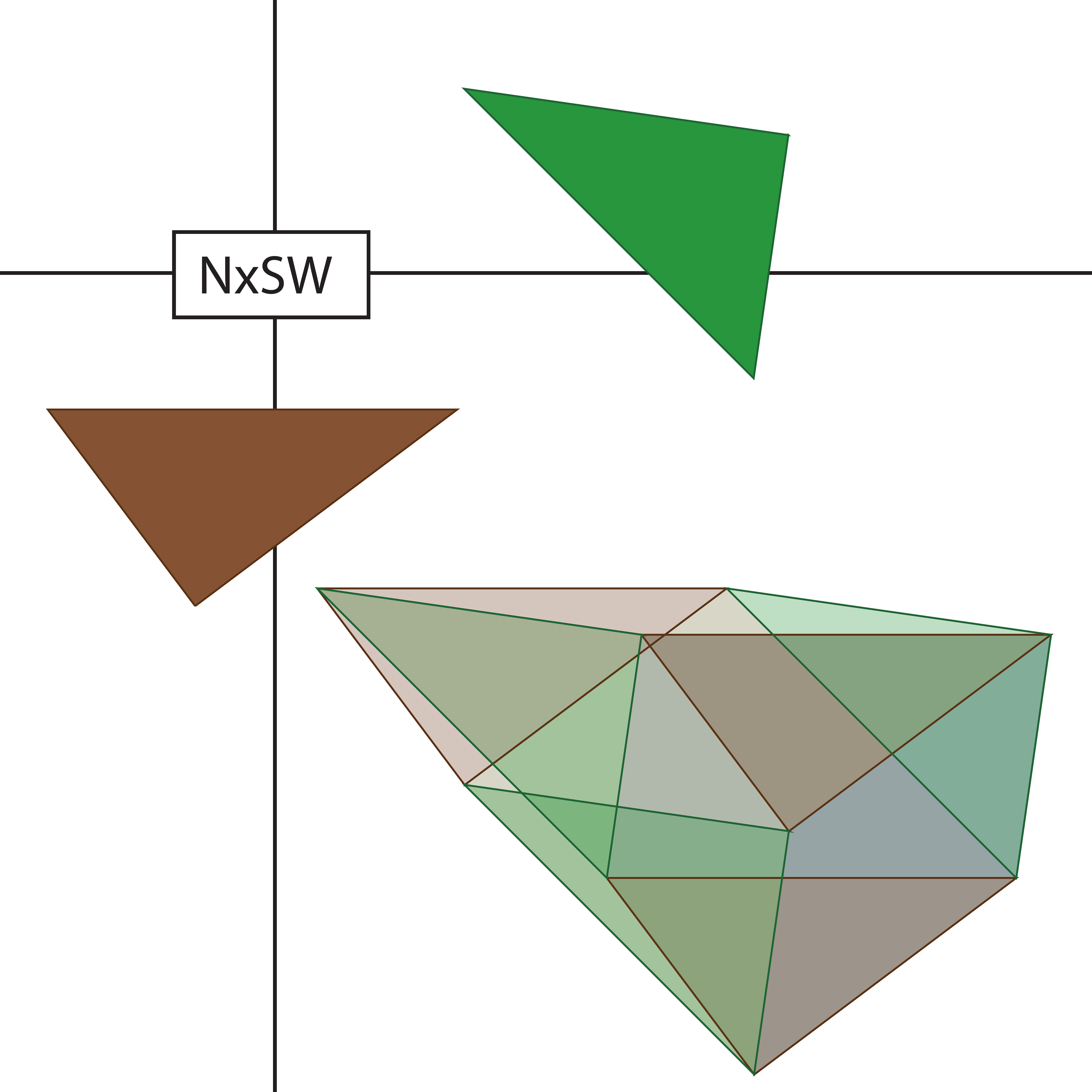}
\includegraphics[scale=.05]{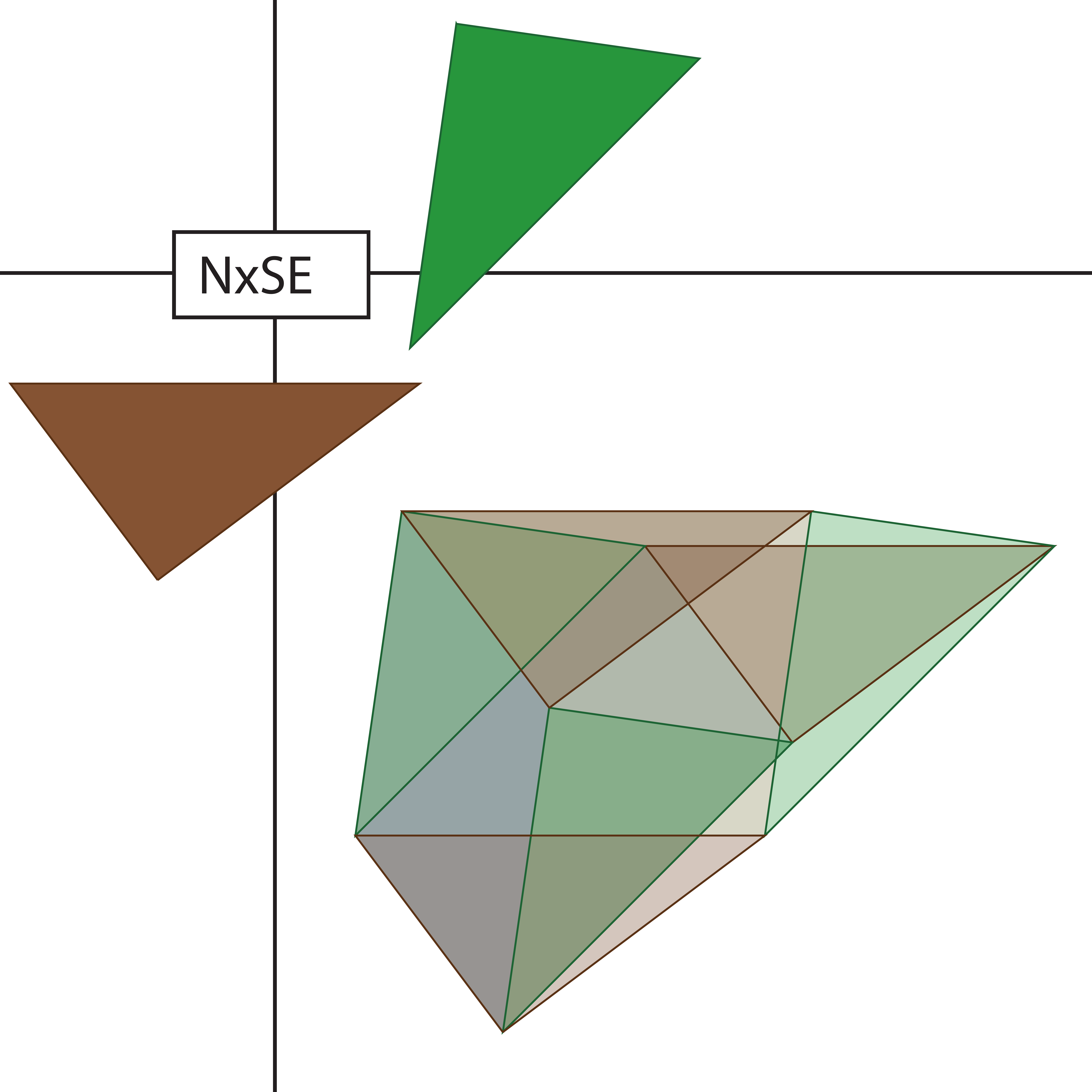}
\includegraphics[scale=.05]{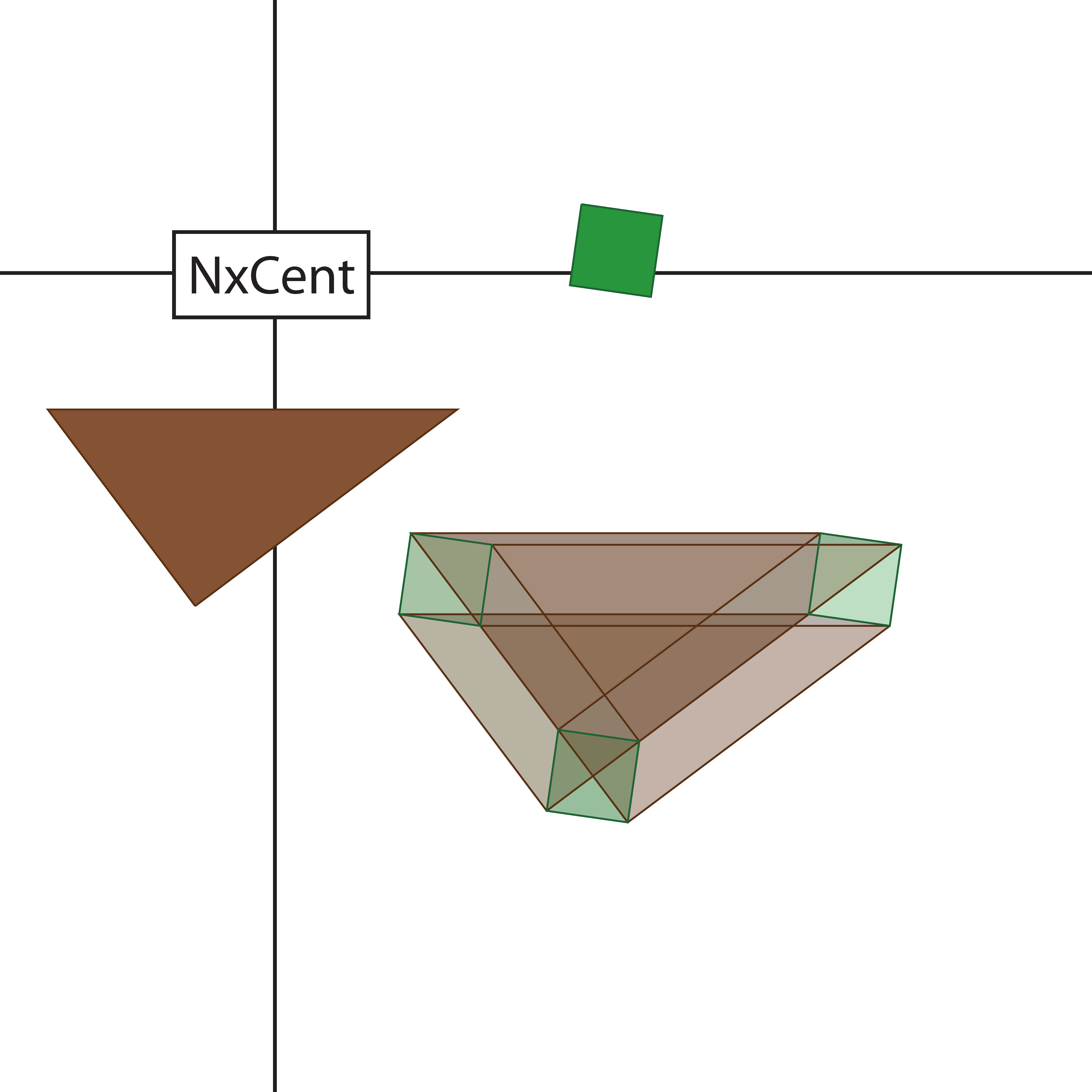}
\end{center}
\caption{The Northern pieces of the  decomposition of the hyper-rectangle} 
\label{North}
\end{figure}

\begin{figure}[htb]
 \begin{center}
\includegraphics[scale=.05]{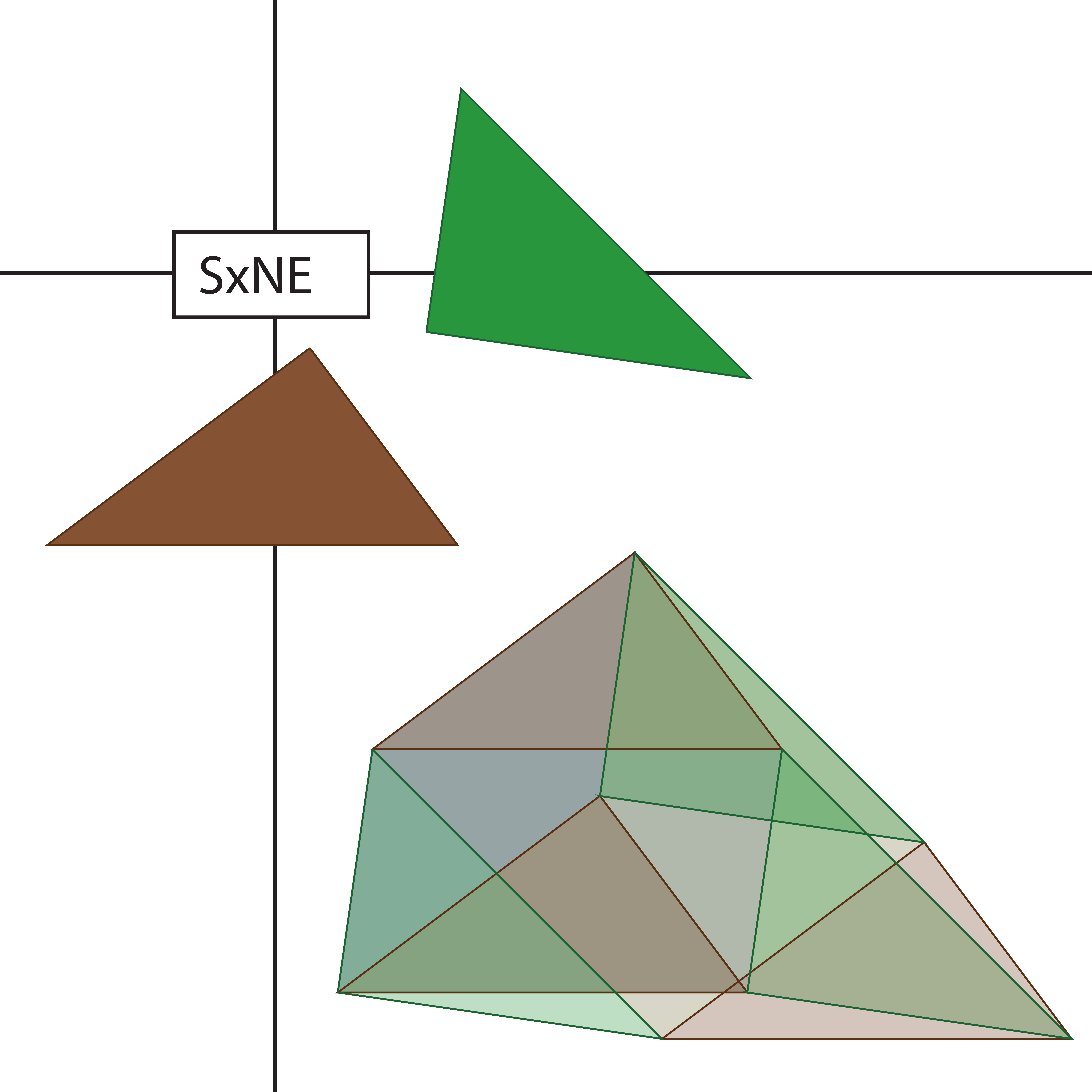}
\includegraphics[scale=.05]{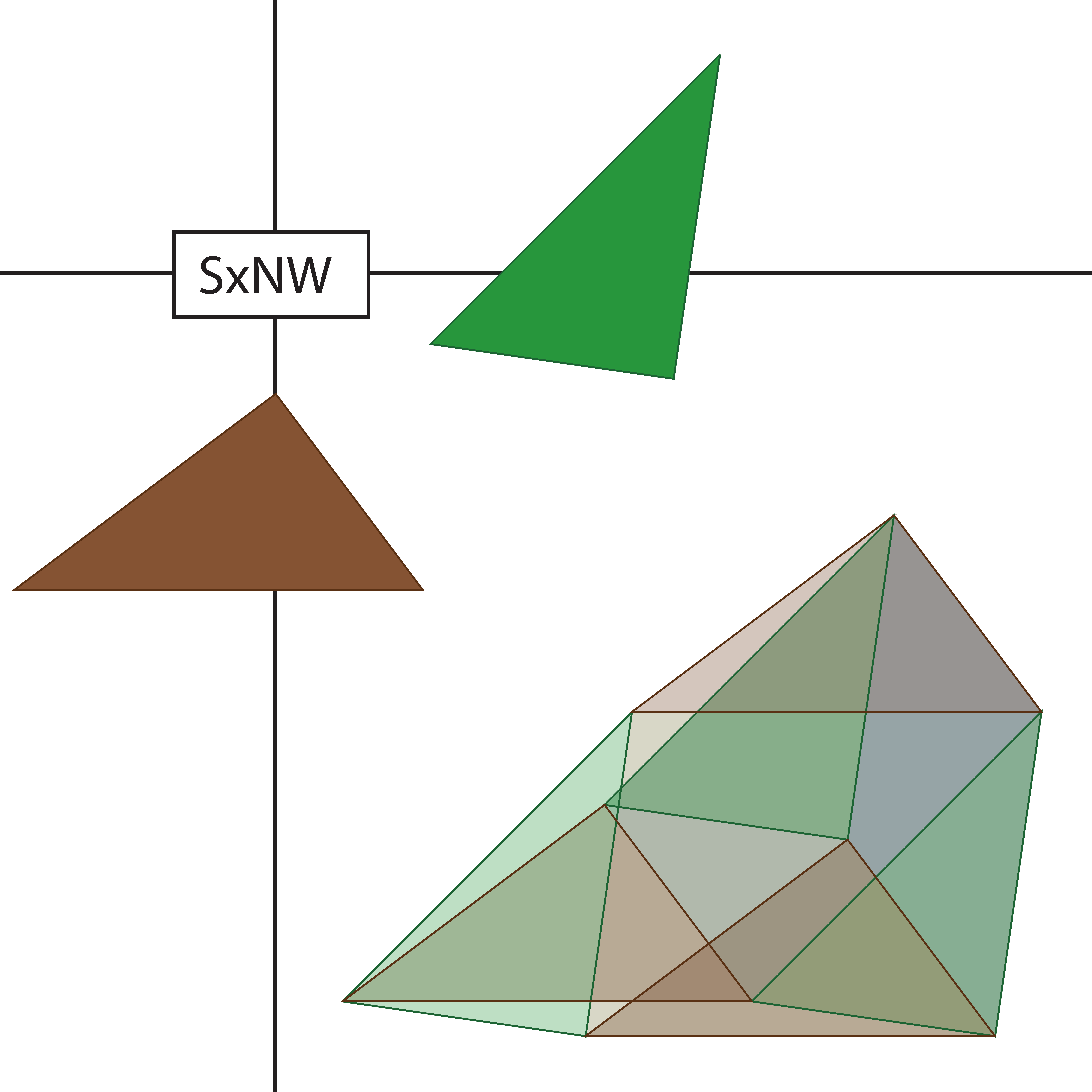}
\includegraphics[scale=.05]{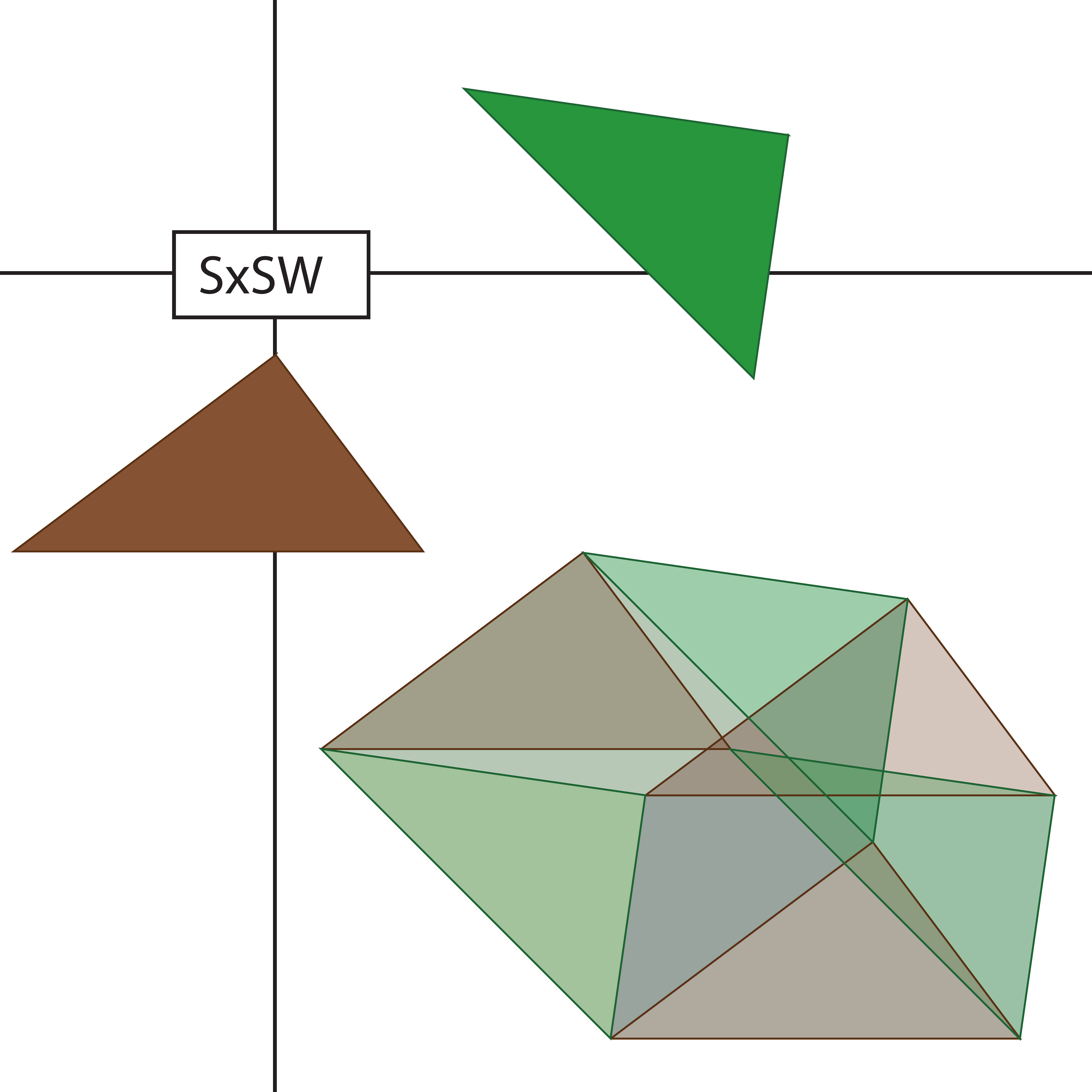}
\includegraphics[scale=.05]{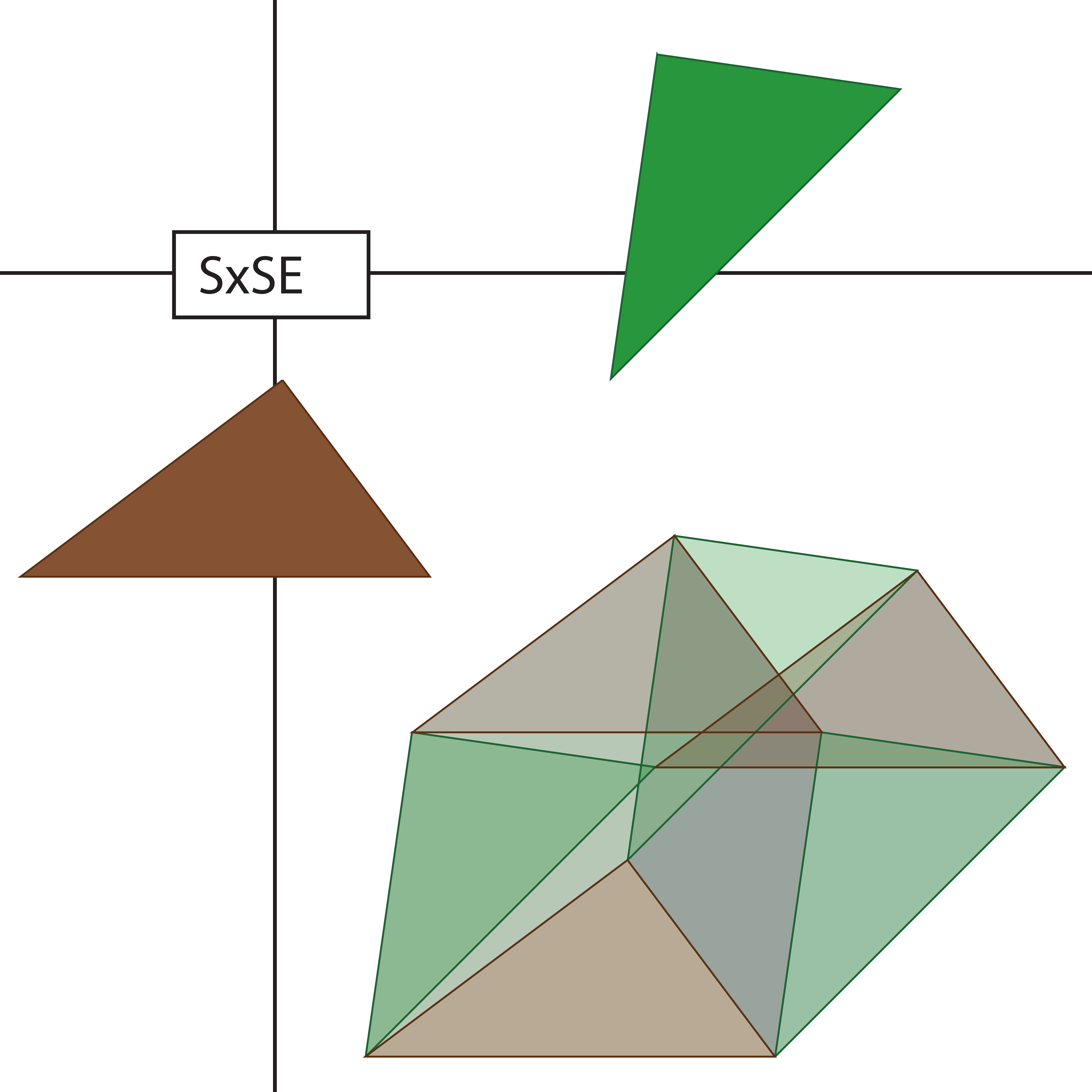}
\includegraphics[scale=.05]{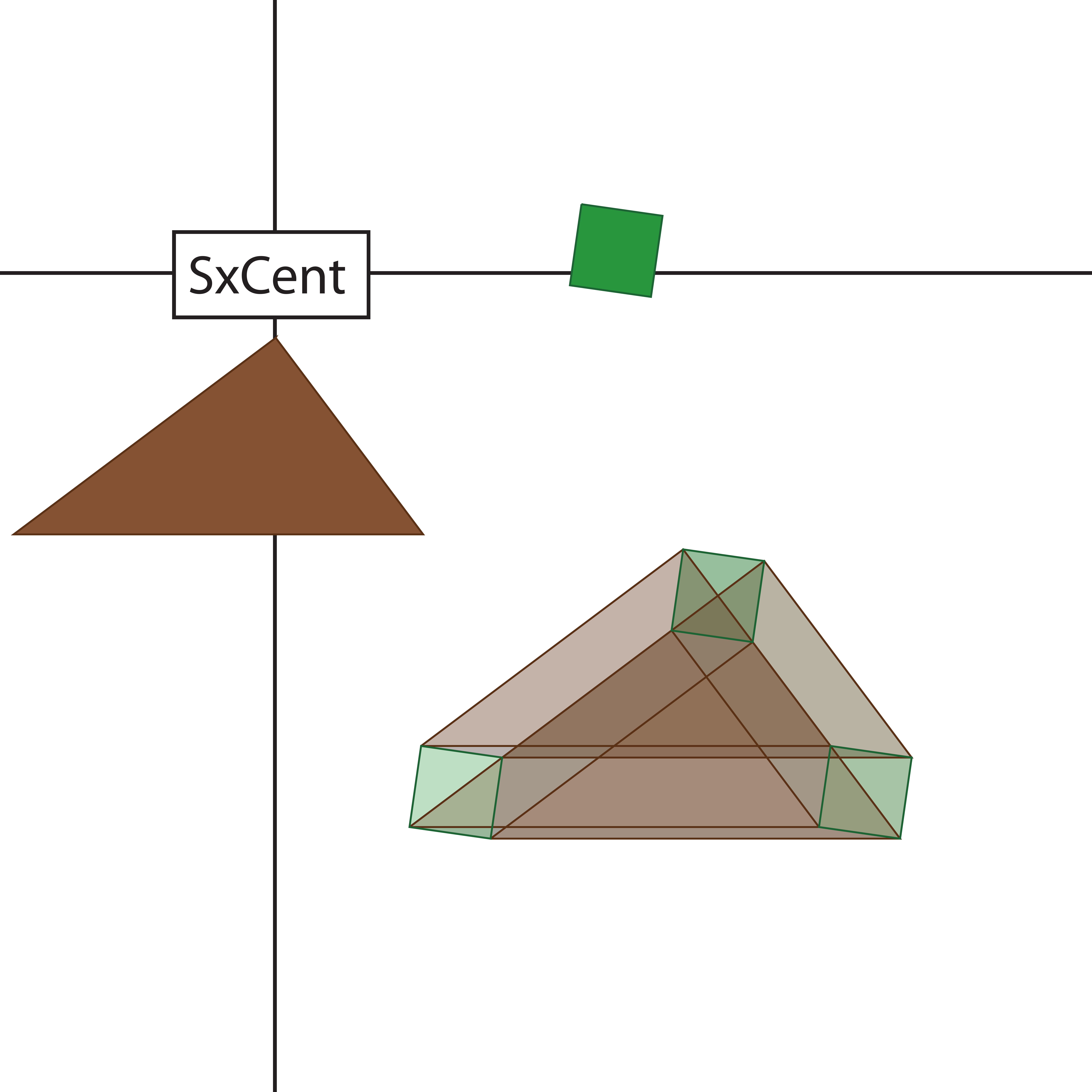}
\end{center}
\caption{The Southern pieces of the  decomposition of the hyper-rectangle} 
\label{South}
\end{figure}

In our illustrations, we project the $zz$ directions into the horizontal and vertical directions of the plane defined by the paper (computer screen), and project the $ww$ directions into directions at a $45^\circ$ angle. In this way, we will speak of the $N$,$S$, $E$ and $W$ triangles, and the $NW$, $NE$, $SW$ and $SE$ triangles. Also, we'll consider the cardinal directions first, and the off-cardinal directions second. Then, in Fig.~\ref{Central} through \ref{South}, we illustrate the $25$ pieces that will result as planar projections of solids of the form (triangle $\times$  triangle), (triangle $\times$ square), (square $\times$ triangle), and (square $\times$ square). 

Next we reassemble these pieces following the reassembly of the Pythagorean Theorem.
Some of these reassemblies are illustrated in Fig.~\ref{West} through Fig.~\ref{South}.  In Fig.~\ref{Central}, the north/south rectangle, the east/west rectangle, and the central brown square are crossed with the central green square. There is a similar shape made from the brown central square and the off-cardinal green rectangles. Of course, the small central hypercube is only used once.

\begin{figure}[htb]
 \begin{center}
\includegraphics[scale=.03]{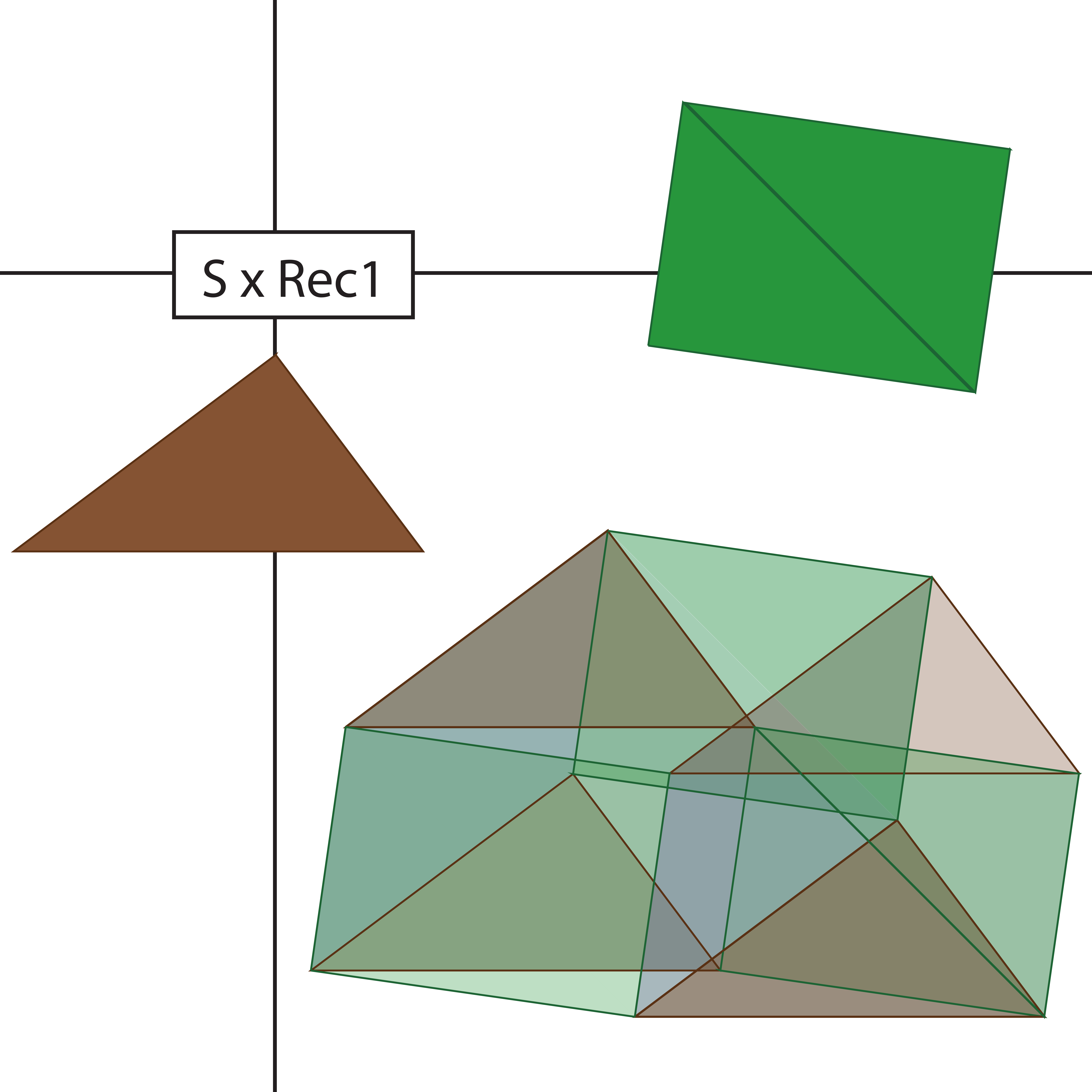}
\includegraphics[scale=.03]{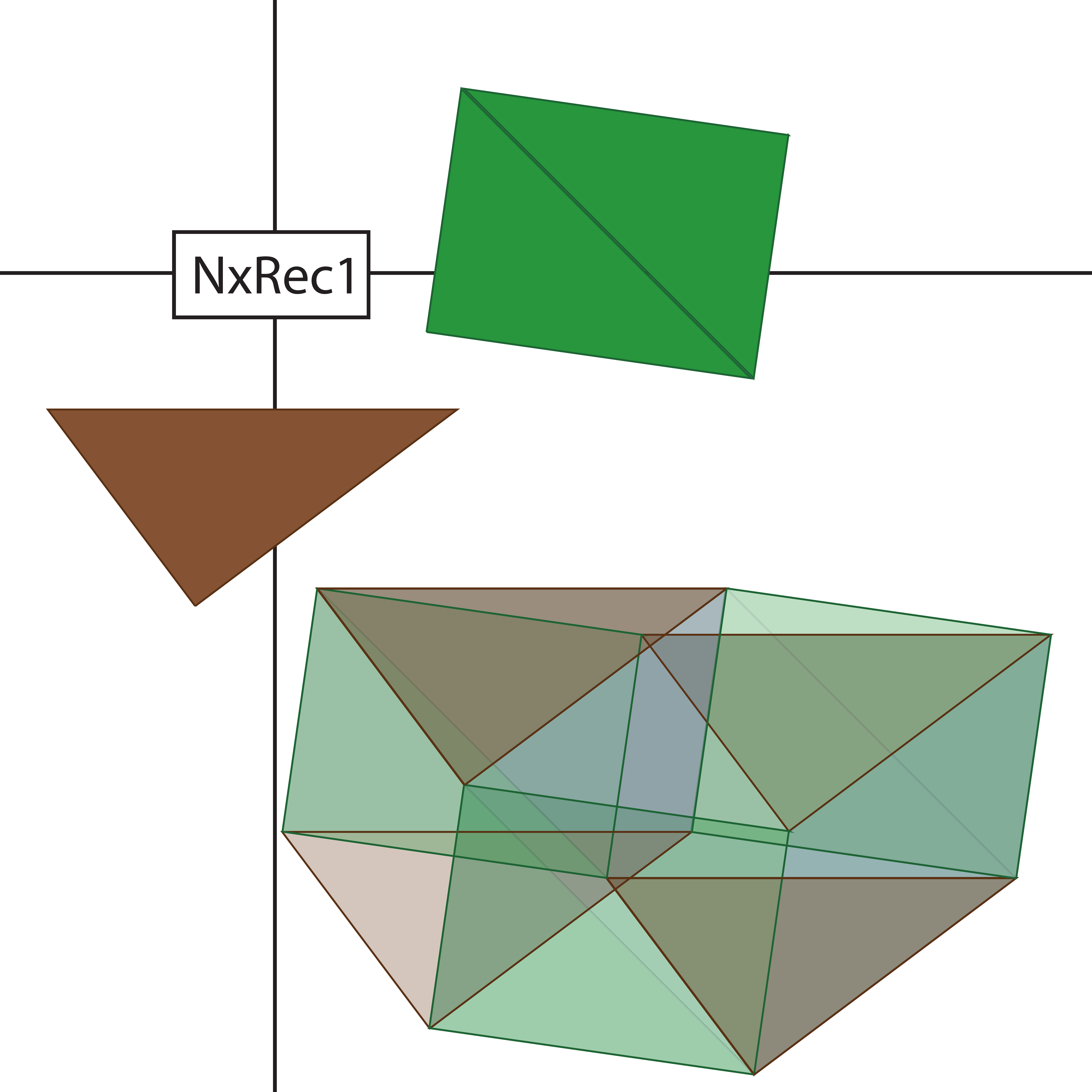}
\includegraphics[scale=.03]{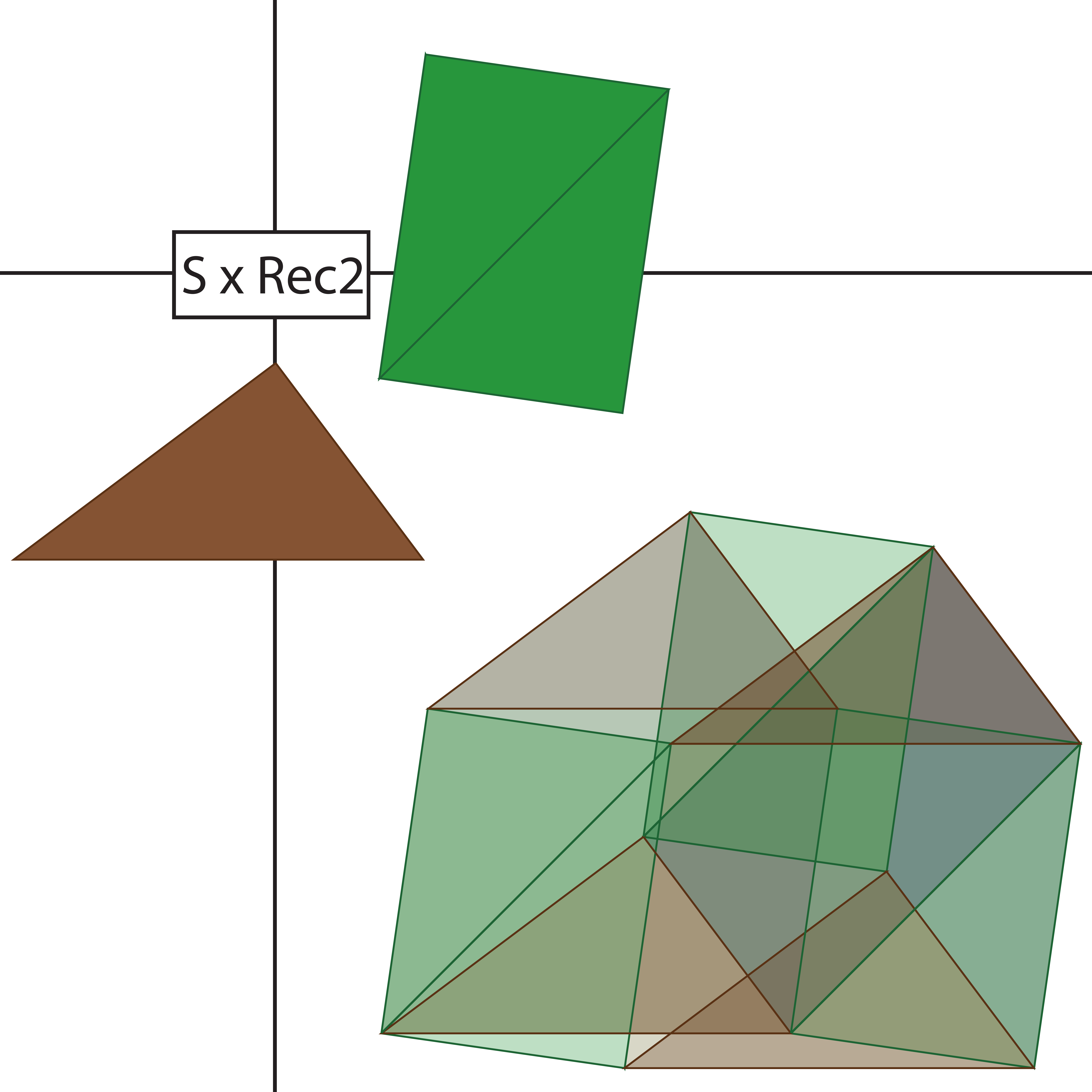}
\includegraphics[scale=.03]{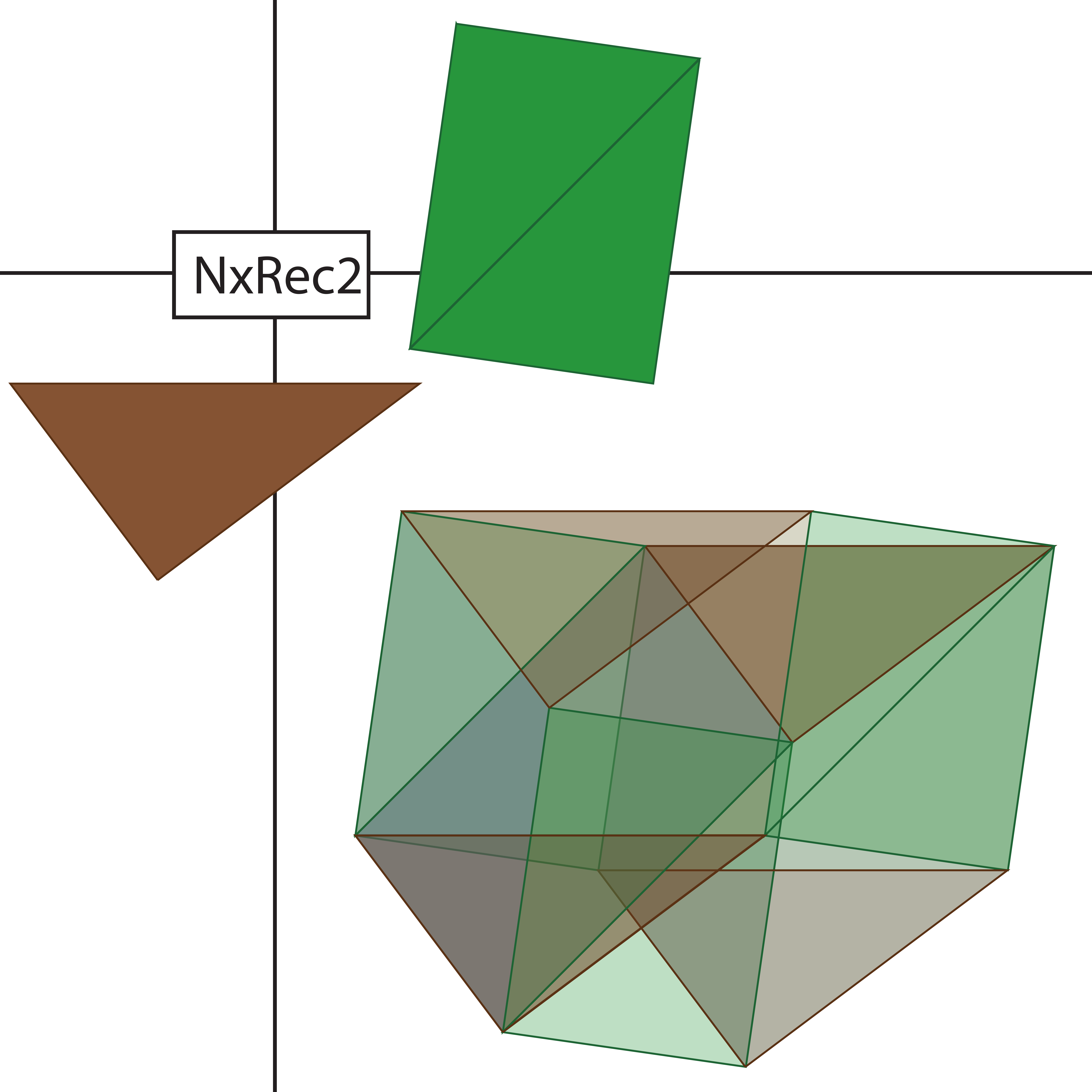}
\end{center}
\caption{Reassembling the eastern and western triangles times rectangles} 
\label{4Rec}
\end{figure}

\begin{figure}[htb]
 \begin{center}
\includegraphics[scale=.1]{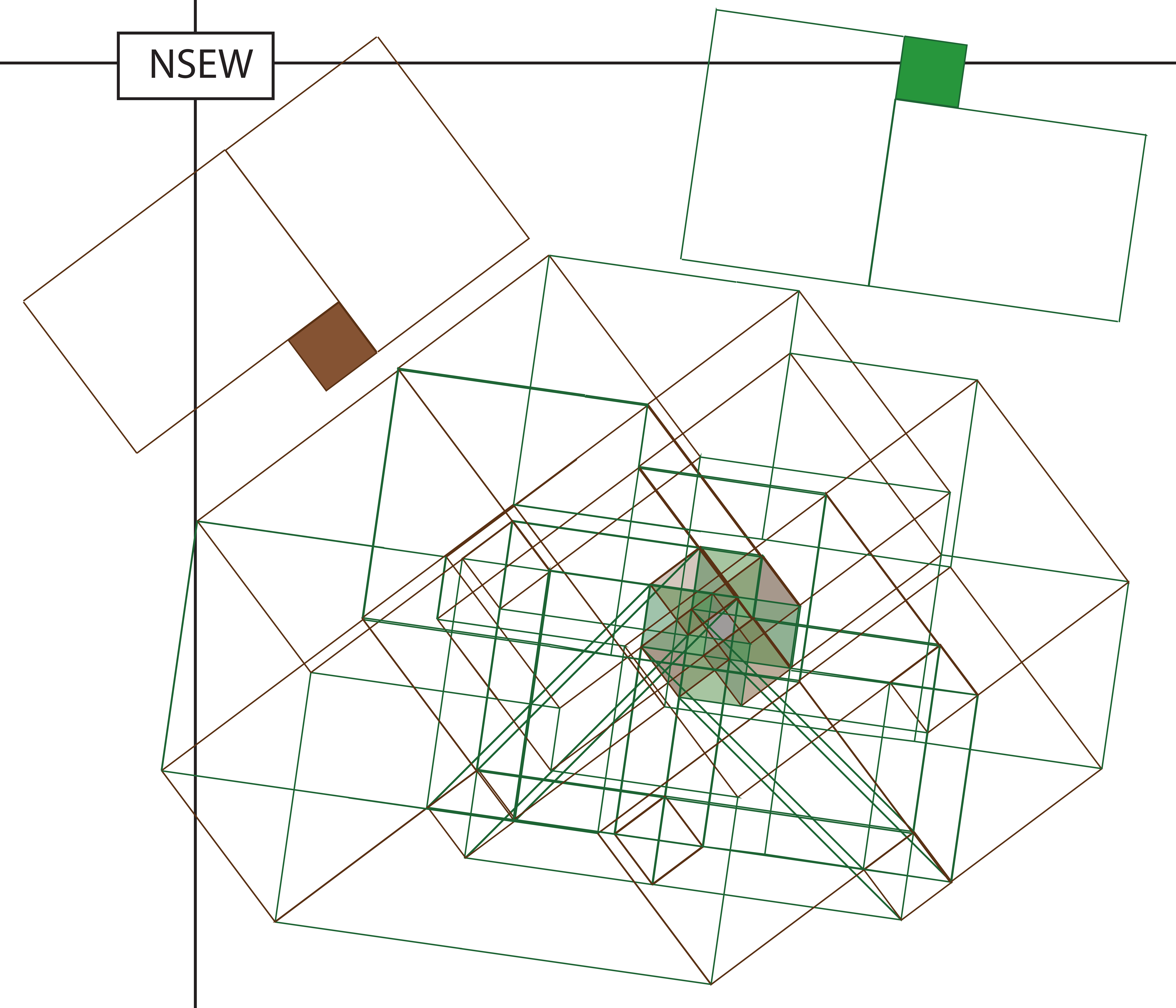}
\end{center}
\caption{The small central square times itself as a small piece} 
\label{NSEW2}
\end{figure}

\begin{figure}[htb]
 \begin{center}
\includegraphics[scale=.1]{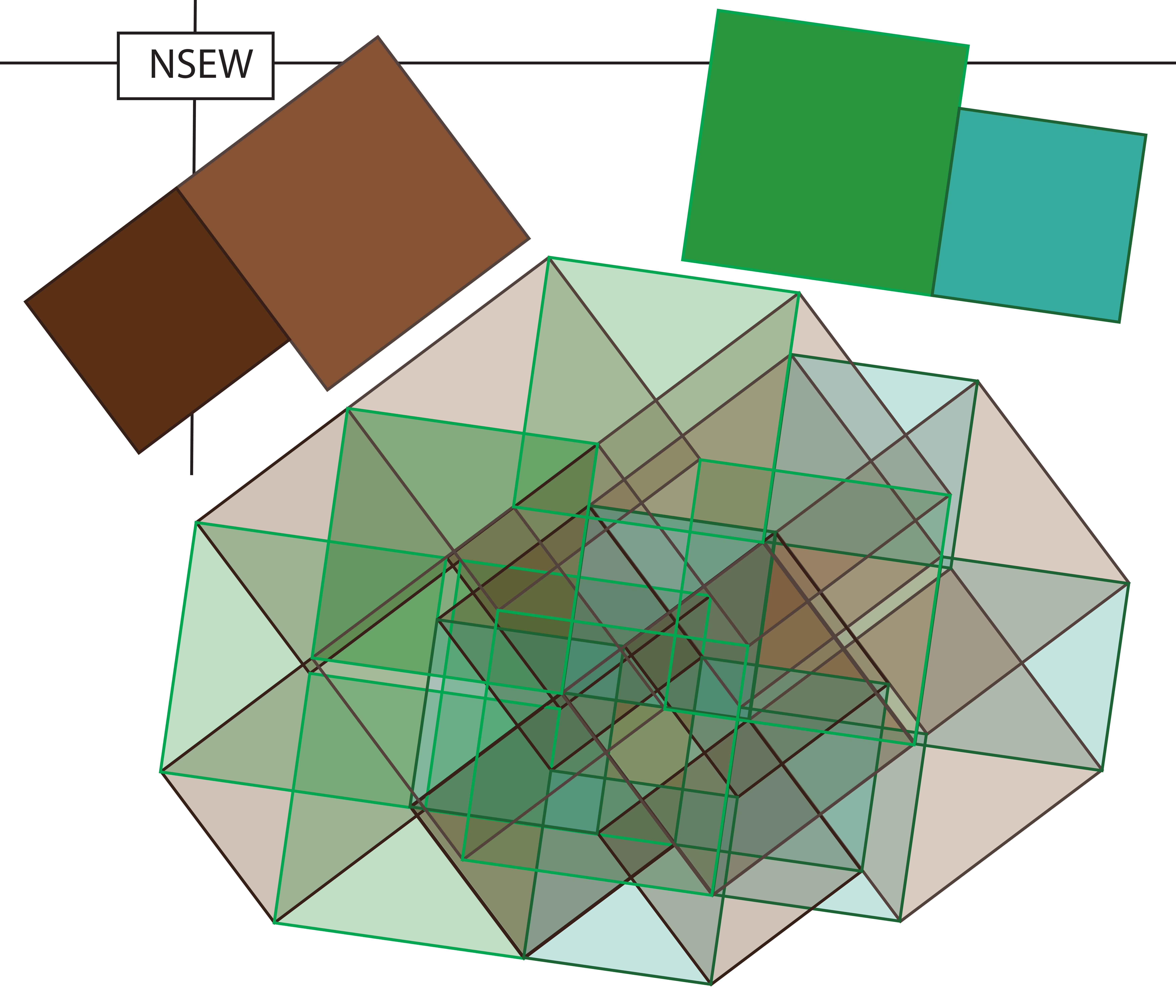}
\end{center}
\caption{Pieces of the sum of squares by sum of squares} 
\label{big}
\end{figure}

The reassemblies continue as the triangles-by-rectangles are reassembled into four rectangles-by-rectangles (Fig.~\ref{4Rec}). Then in Fig.~\ref{NSEW2}, the small central hypercube is placed in the $1$-skeleton of the full reassembly. 
The Fig.~\ref{4hyp} demonstrates the outlines of the four hypercubes fit into this conglomeration. Finally, Fig.~\ref{big} shows all four of these square-by-squares inside this conglomeration. 

In this way, we have illustrated that the algebraic identity 
$$z^2w^2=(x^2+y^2)(u^2+v^2)= x^2u^2 + y^2u^2 +x^2 v^2 + y^2 v^2$$
is realized as a decomposition of the original hyper-rectangle $R_{zzww}$ into the union of four hyper-rectangles $R_{xxuu}$, $R_{yyuu}$, $R_{xxvv}$, and $R_{yyvv}$.

\begin{figure}[htb]
 \begin{center}
 \includegraphics[scale=.03]{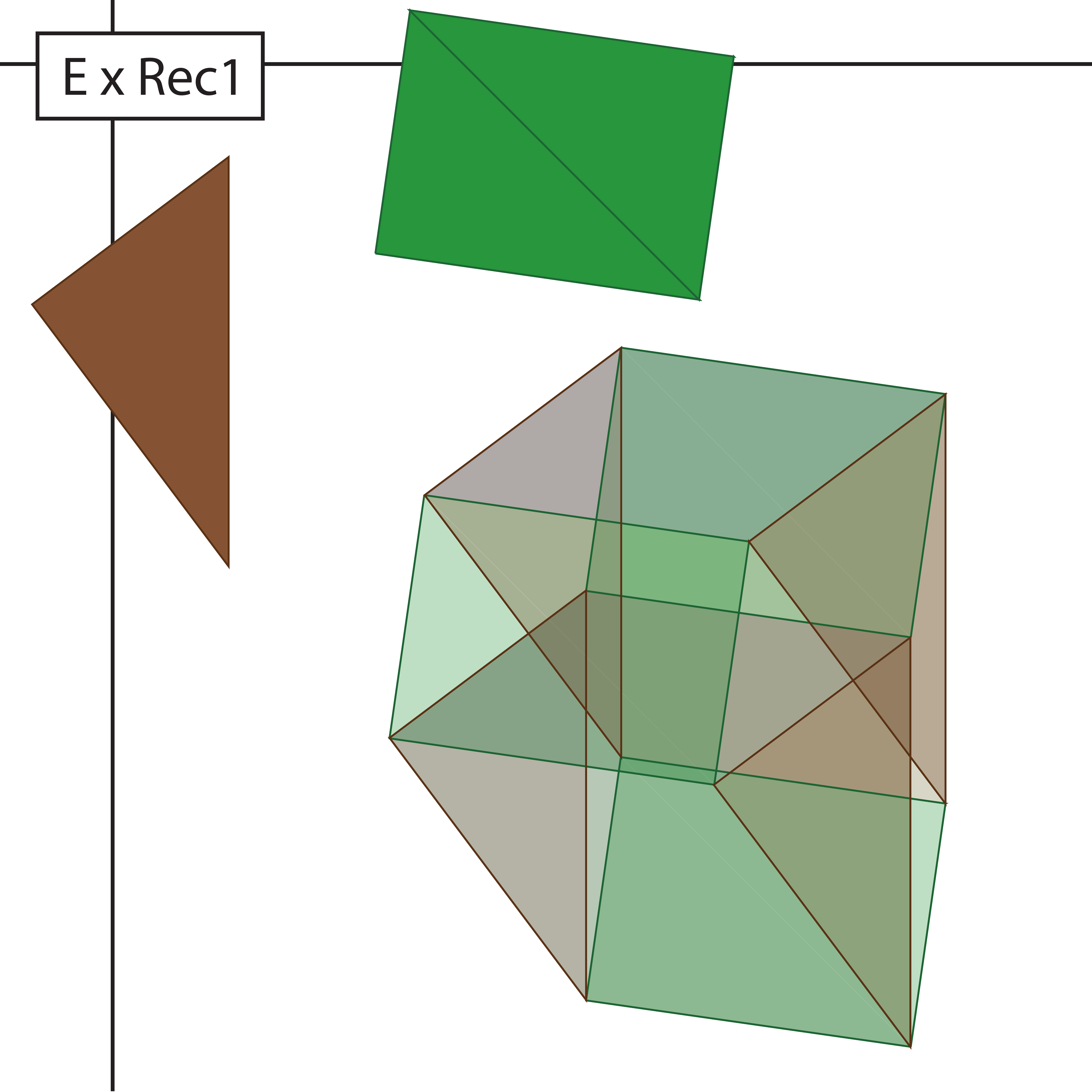}
\includegraphics[scale=.03]{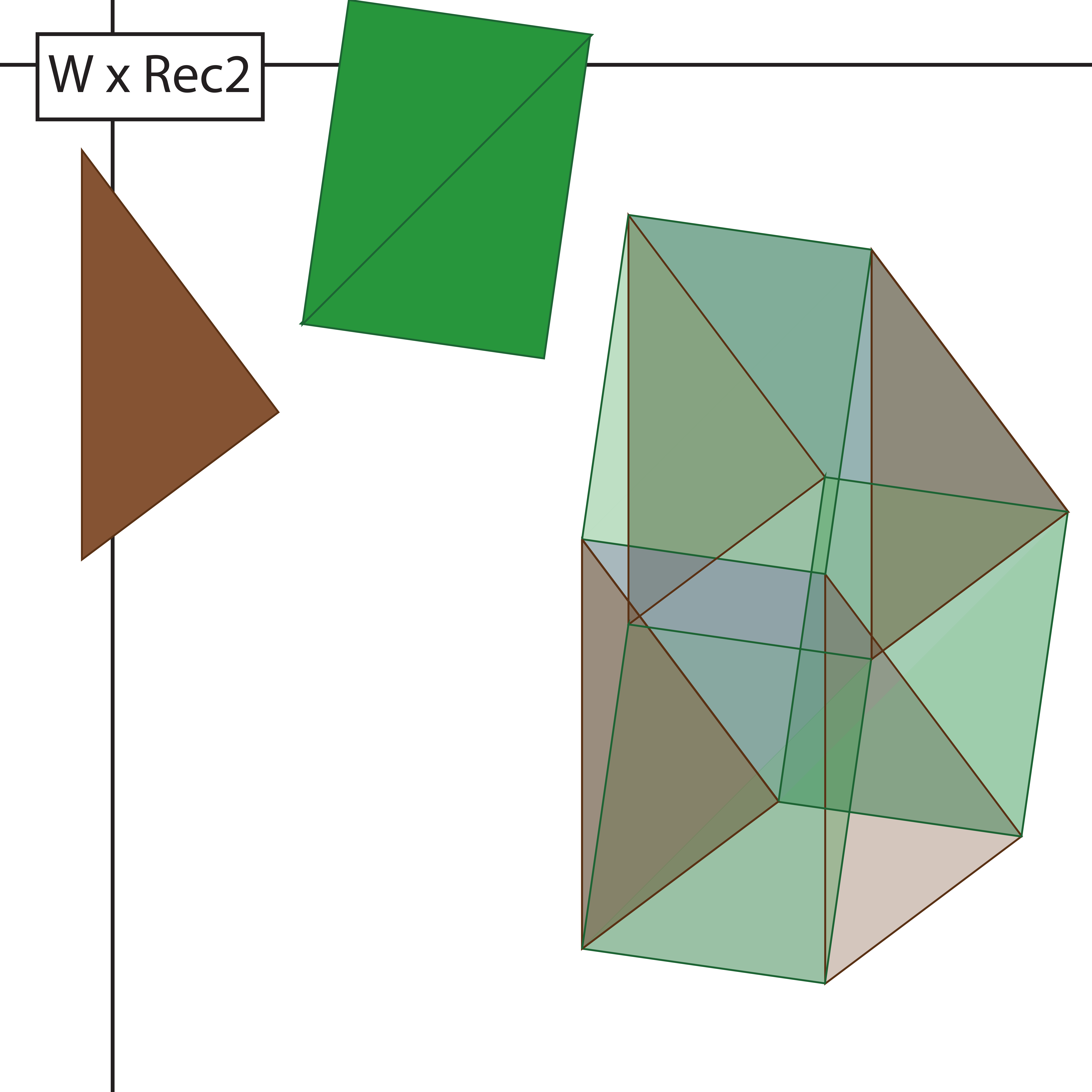}
\includegraphics[scale=.03]{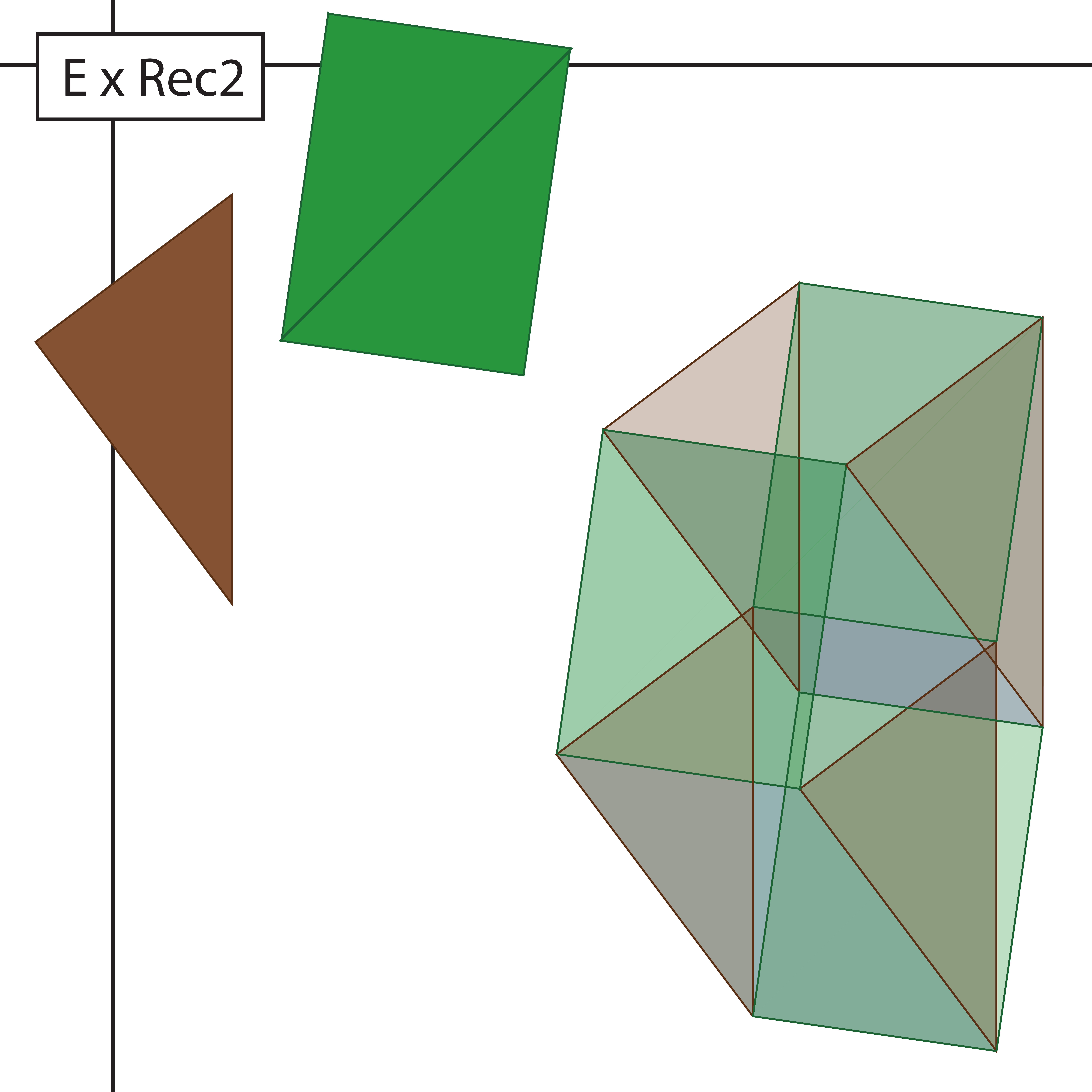}
\includegraphics[scale=.03]{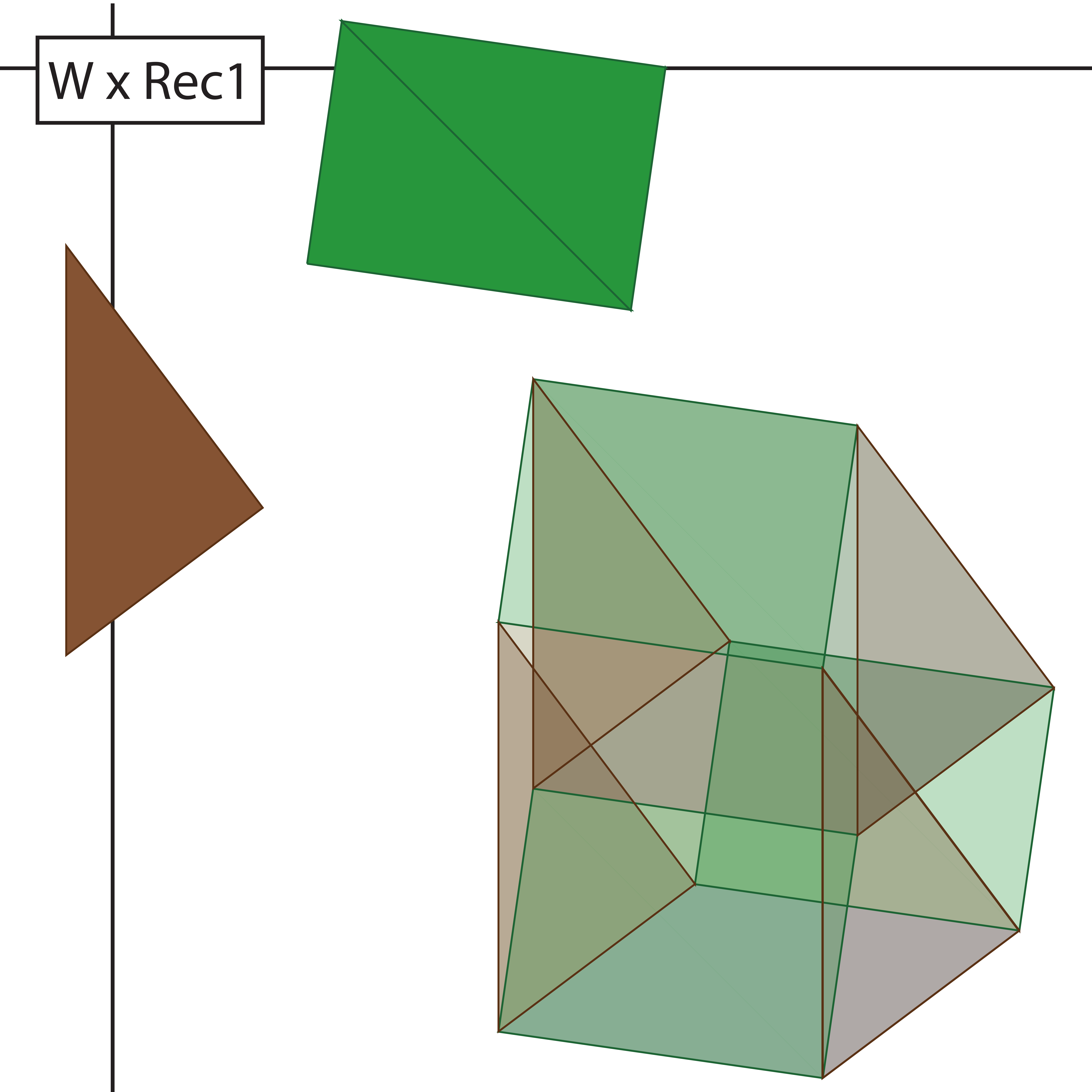}
\end{center}
\caption{Reassembling the eastern and western triangles times rectangles} 
\label{MoreWest}
\end{figure}


\begin{figure}[htb]
 \begin{center}
\includegraphics[scale=.04]{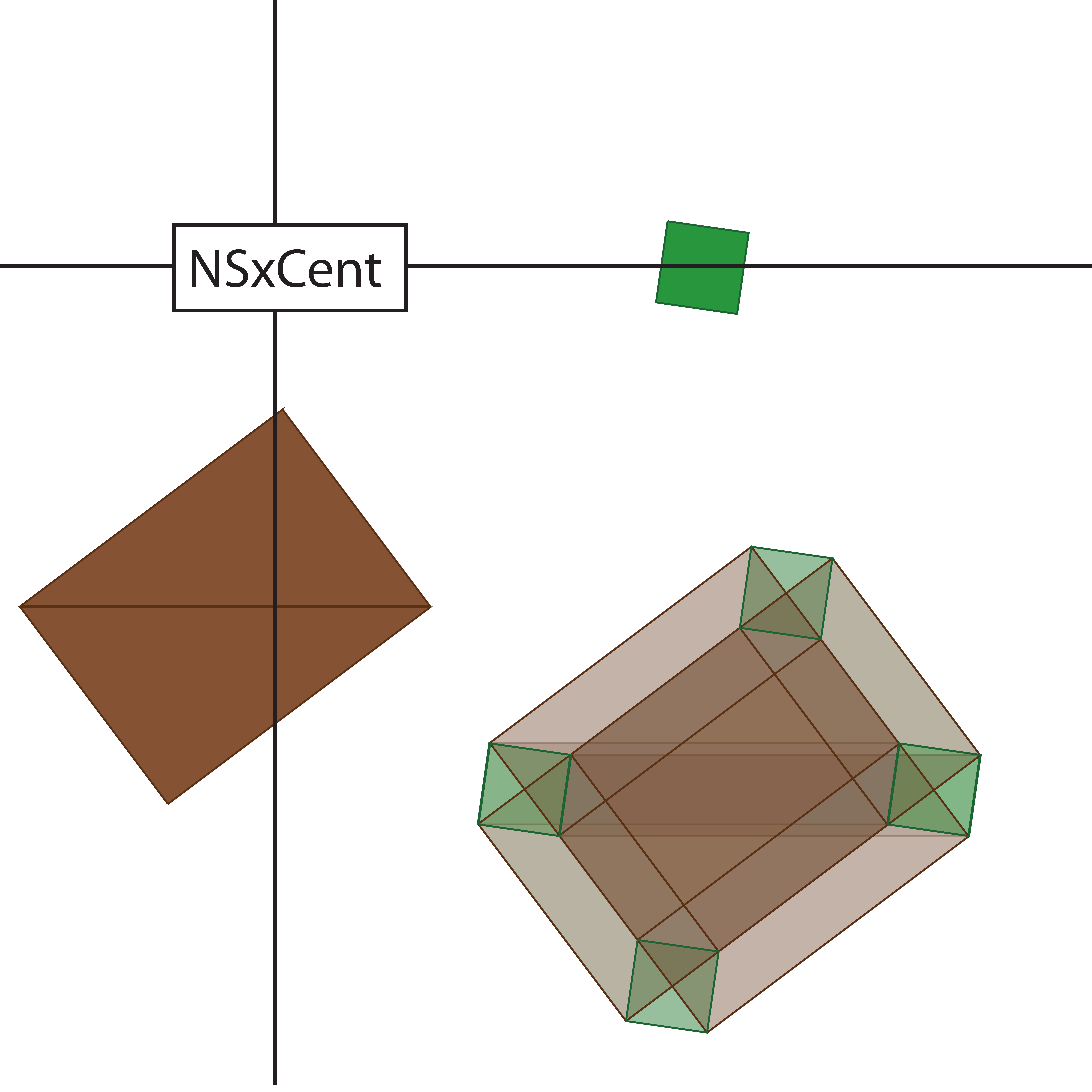}
\includegraphics[scale=.04]{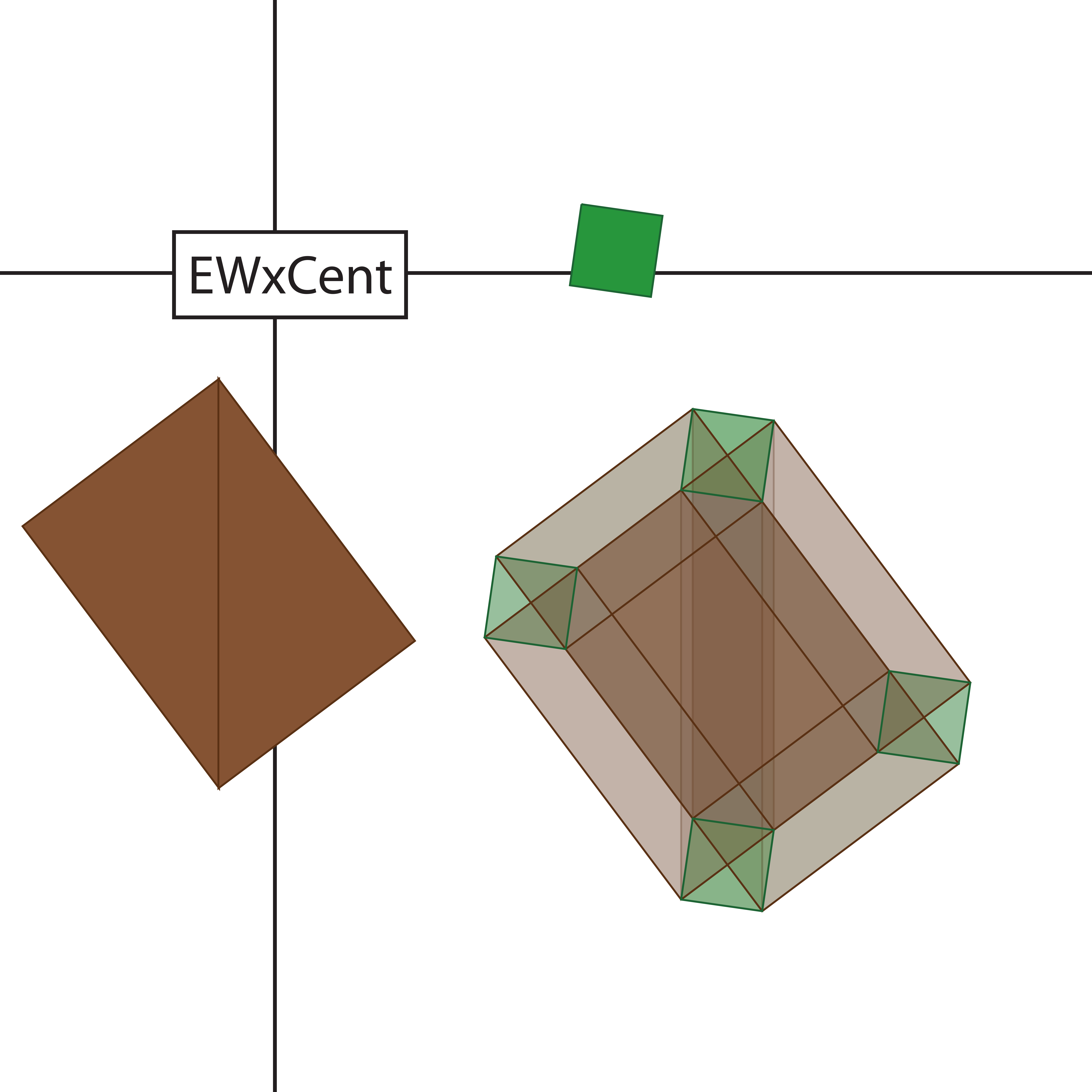}
\includegraphics[scale=0.04]{CenterxCenter}
\end{center}
\caption{The north/south, east/west, and the central square times the center} 
\label{NSCenter}
\end{figure}


\begin{figure}[htb]
 \begin{center}
\includegraphics[scale=.03]{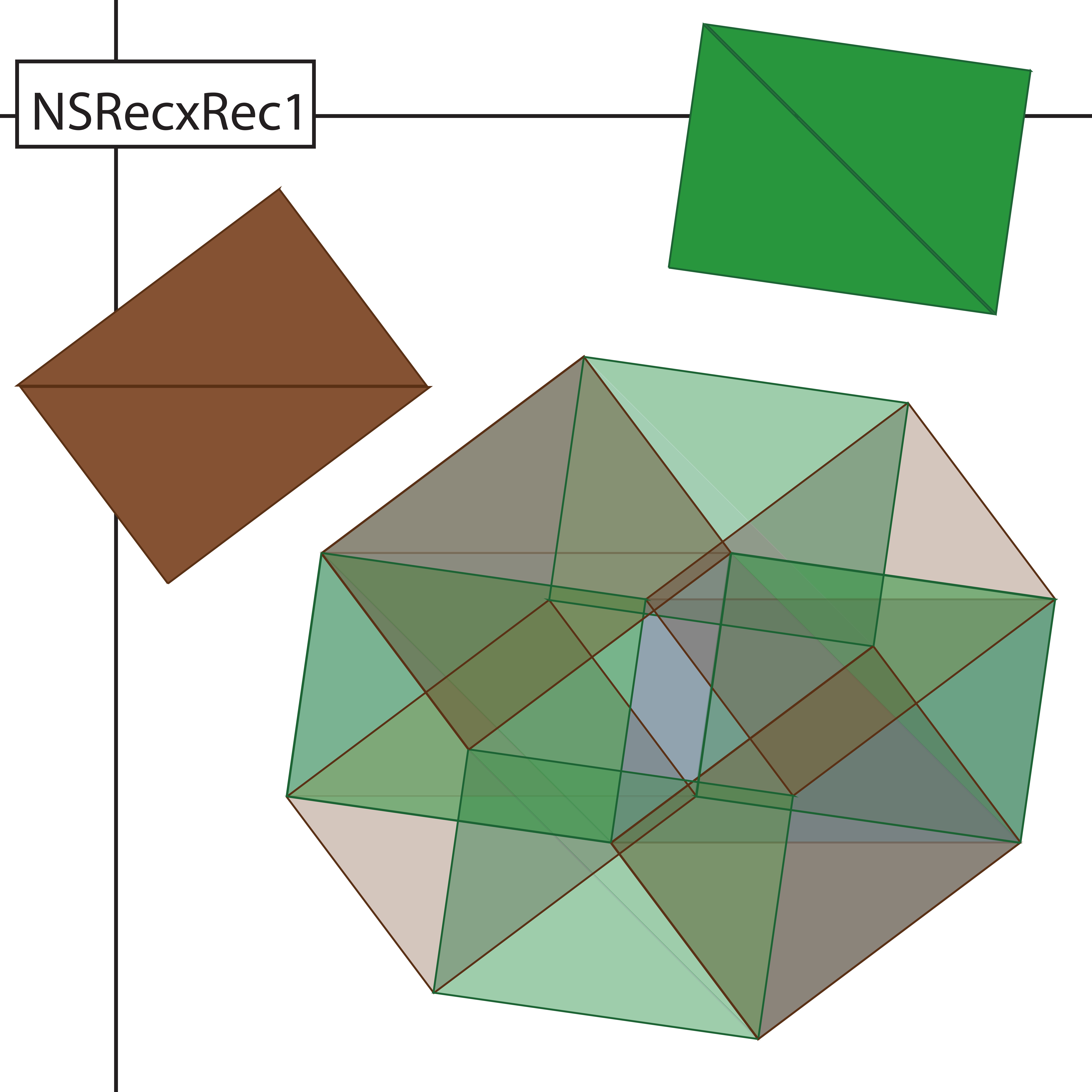}
\includegraphics[scale=.03]{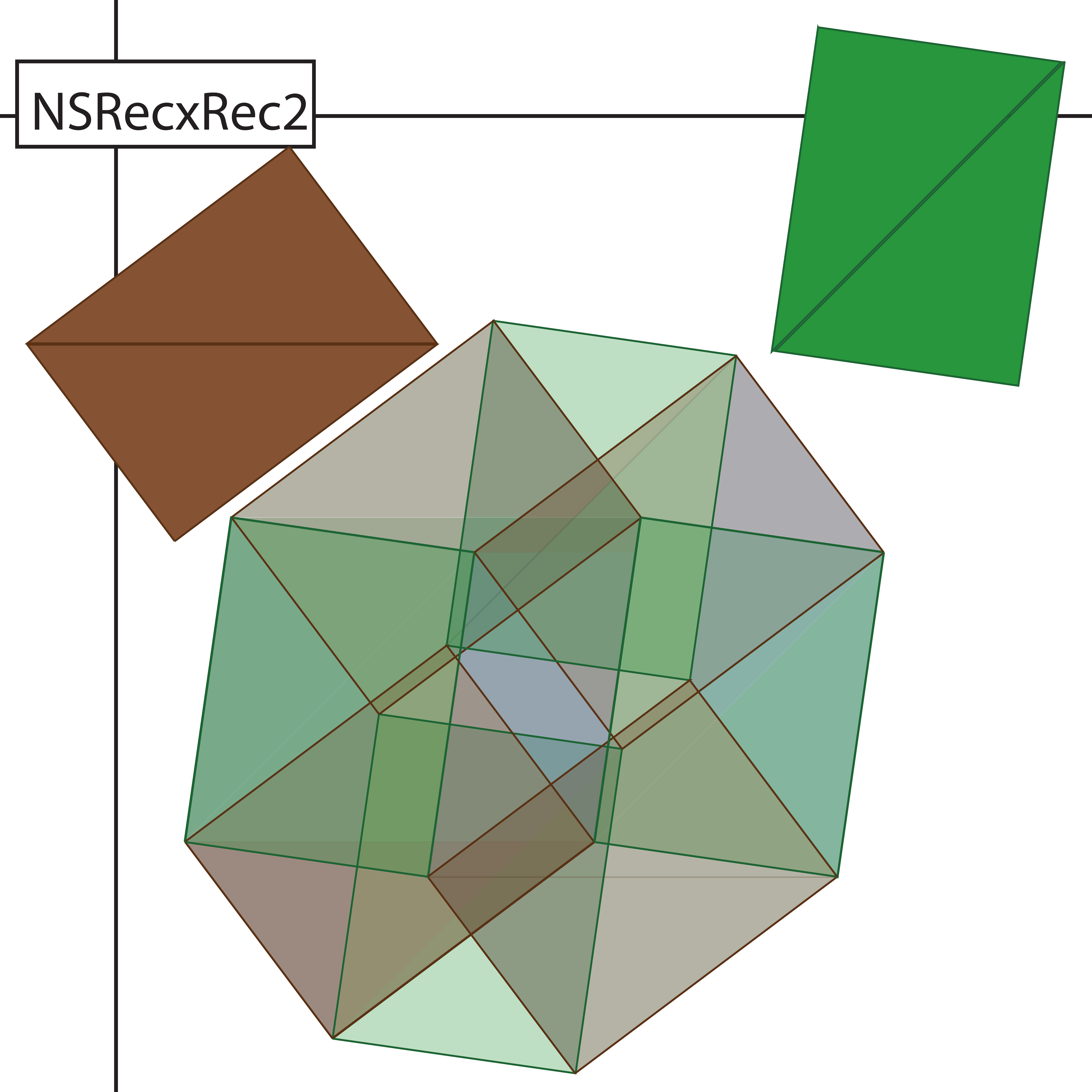}
\includegraphics[scale=.03]{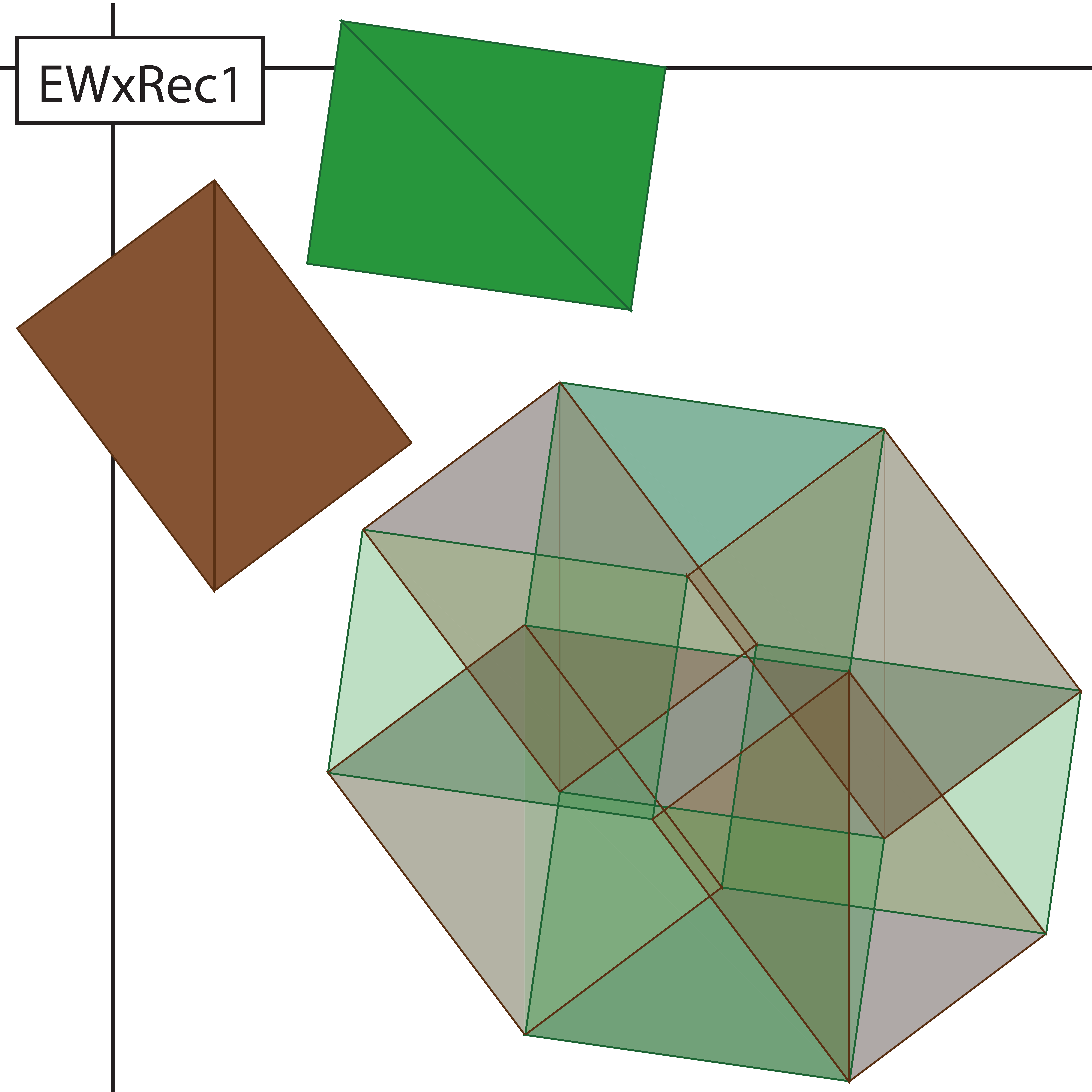}
\includegraphics[scale=.03]{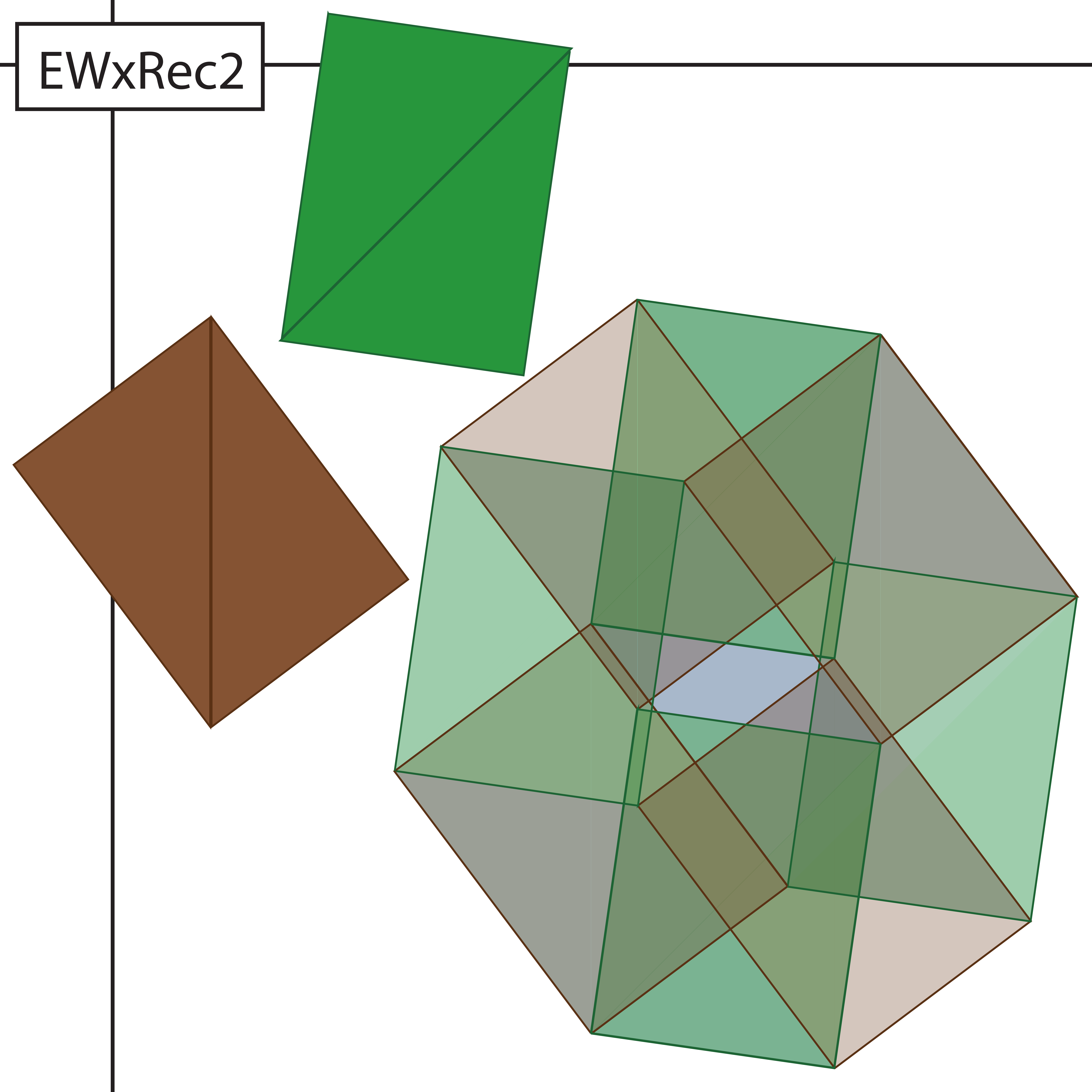}
\end{center}
\caption{Products of rectangles} 
\label{4hyp}
\end{figure}

\begin{figure}[htb]
 \begin{center}
\includegraphics[scale=.045]{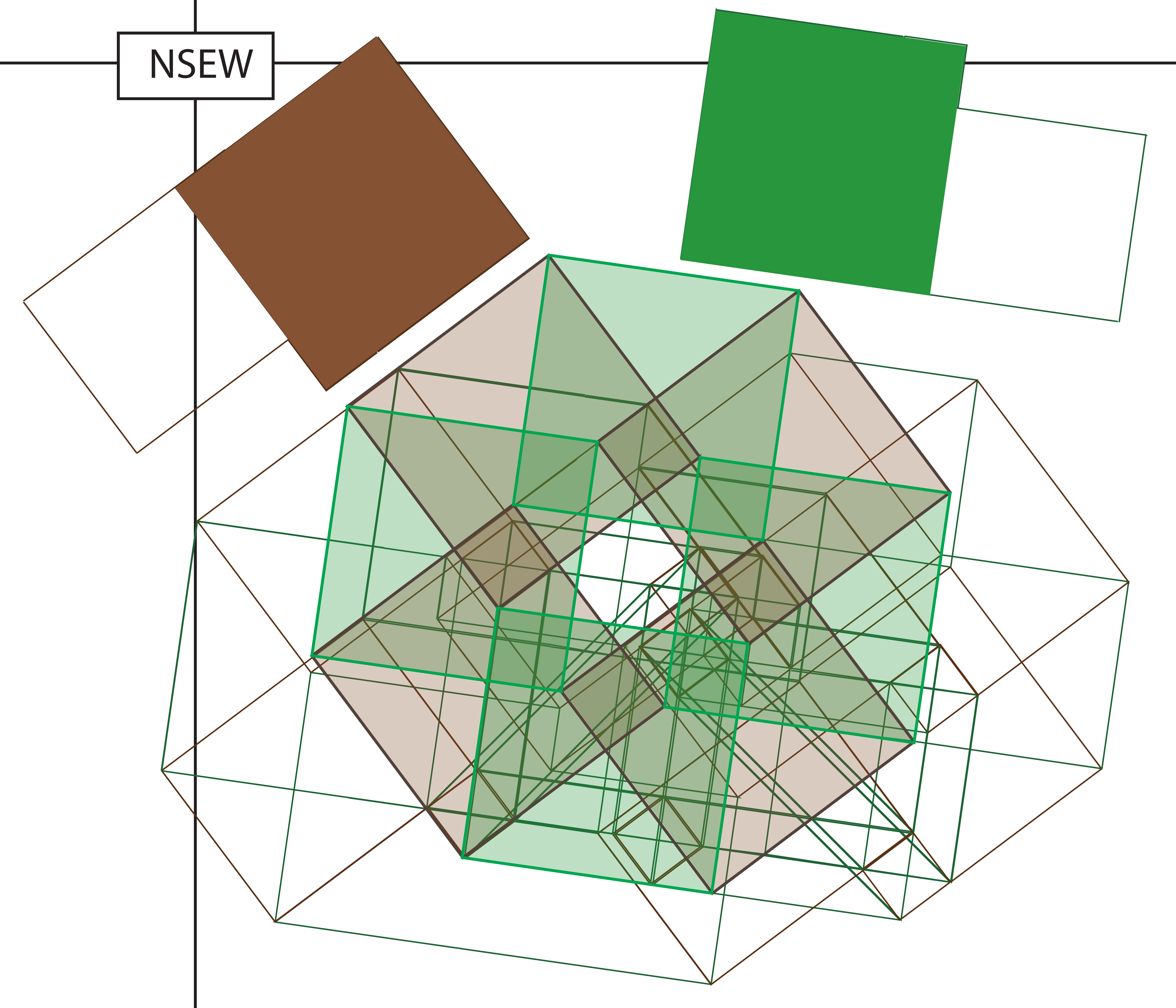} \rule{.5cm}{0cm}
\includegraphics[scale=.045]{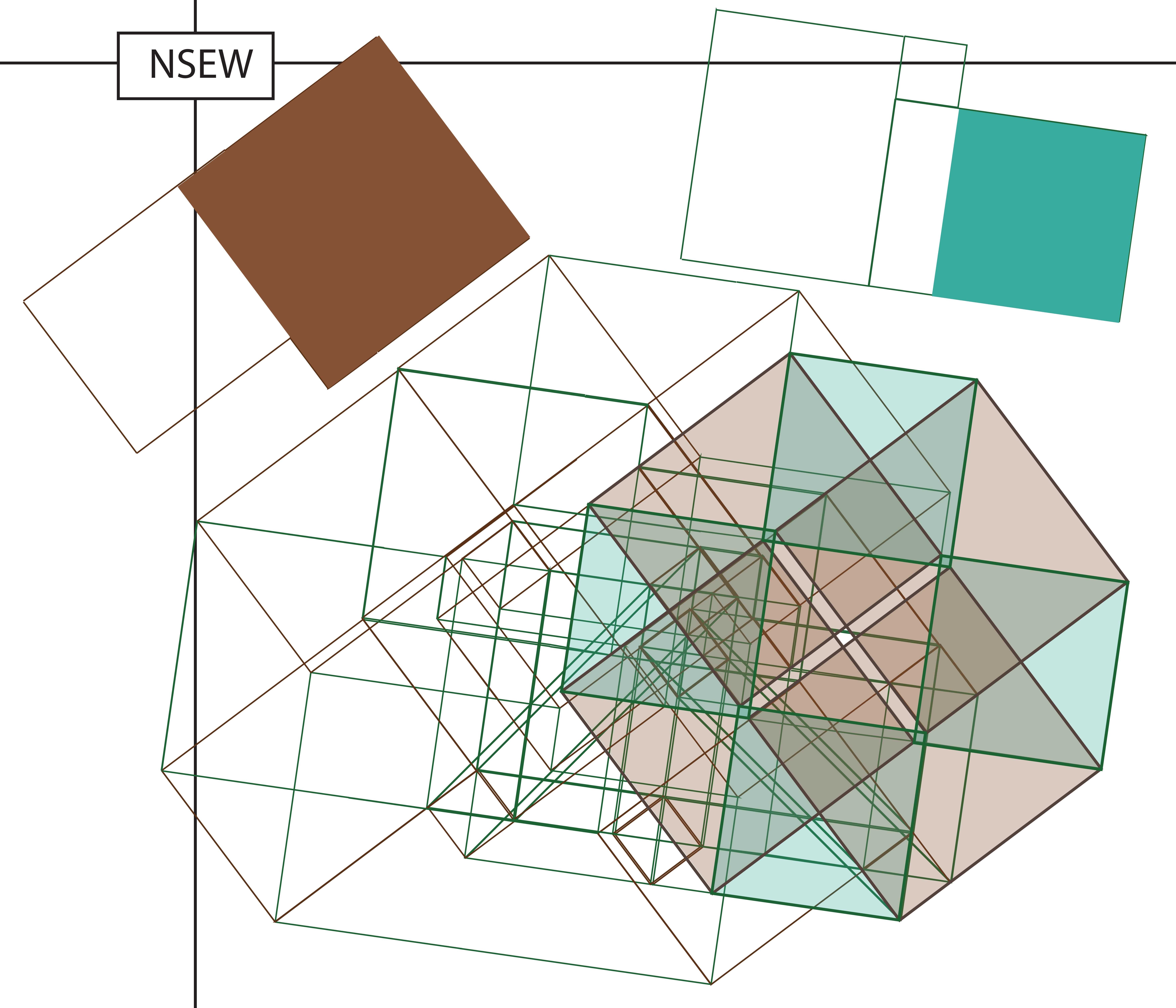}
\includegraphics[scale=.045]{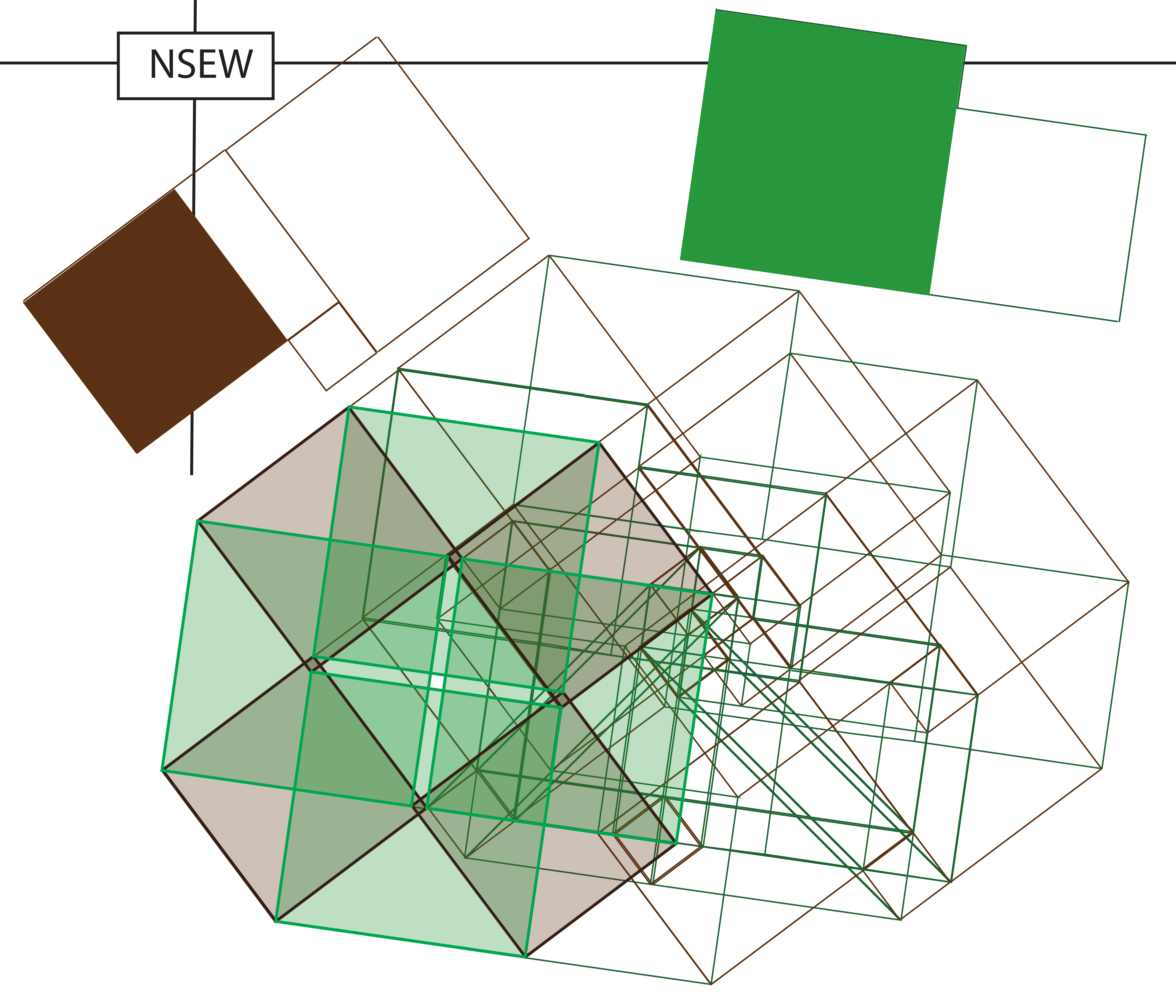} \rule{.5cm}{0cm}
\includegraphics[scale=.045]{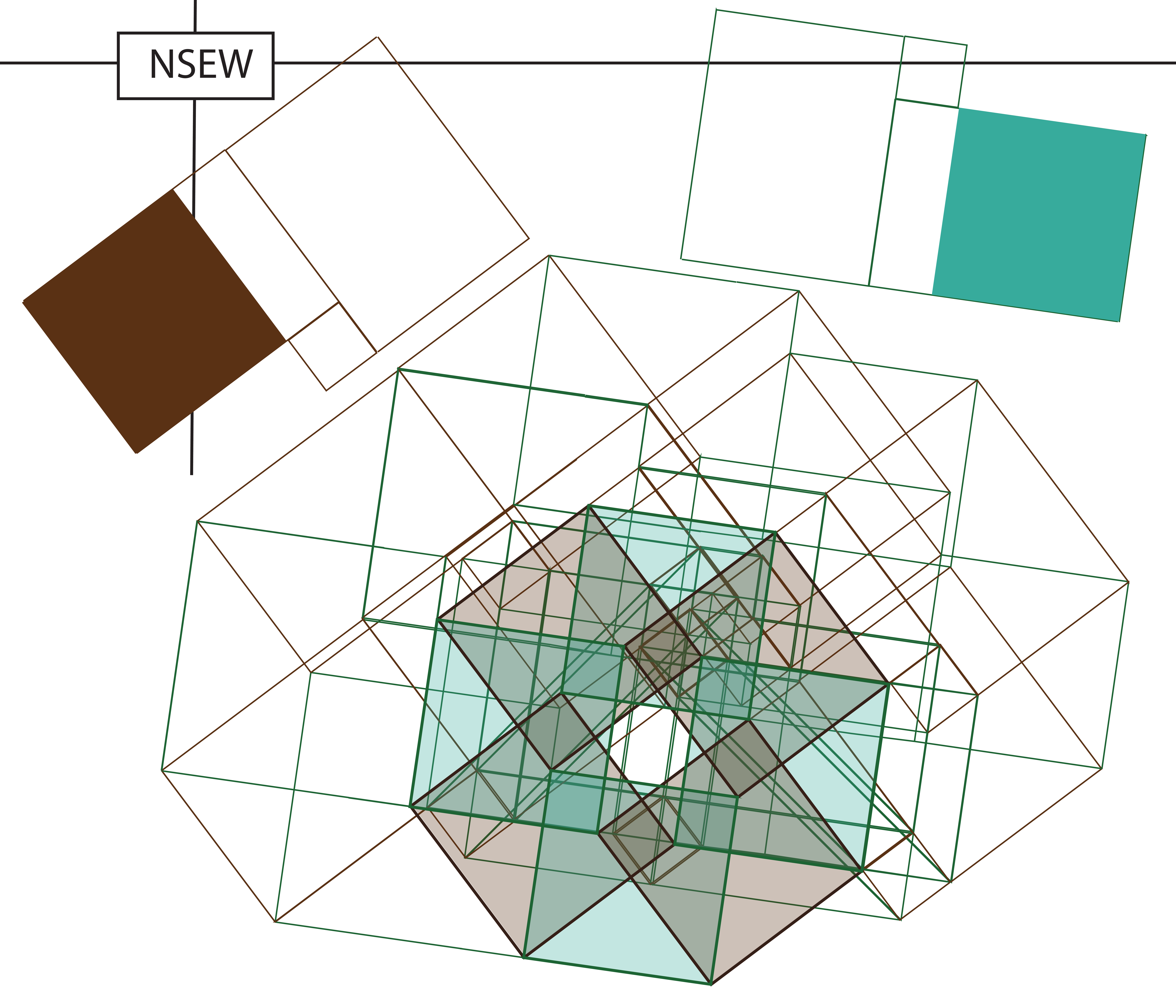} \rule{0cm}{1cm}
\end{center}
\caption{Pieces of the sum of squares by sum of squares} 
\label{4hyp}
\end{figure}


We apply these decompositions to the factors $a^4$, $b^4$, and $a^2b^2$ in the expansion of 
$$(a+b+c)(a+b-c)(a-b+c)(-a+b+c)= 2a^2b^2+2a^2c^2 +2b^2c^2-a^4-b^4-c^4.$$
Here, $a^2=r^2+s^2$, and $b^2=s^2+(p-r)^2$. In addition, the terms of the form $a^2c^2$ and $b^2 c^2$ can be reformed by decomposing $a^2$ or $b^2$ into the five pieces in the proof of the Pythagorean Theorem, and taking the Cartesian product of each with a square of area $c^2=p^2$. When the five pieces are reassembled as union of two squares, the resulting $4$-dimensional solids are also reassembled. 

\section{Products of Sums and Differences of Squares}
\label{sumdiff}

In our next two illustrations, we indicate that products of the form $(x^2-y^2)(x^2-y^2)$ and products of the form $(x^2+y^2)(z^2+w^2)$, when thought of as hyper-volumes of $4$-dimensional solids, can be reassembled into hyper-rectangles whose volumes are the expected terms:
$$(x^2-y^2)(x^2-y^2)=x^4-2x^2y^2+y^4$$ and
$$(x^2+y^2)(z^2+w^2)=x^2z^2+ x^2w^2+y^2z^2 + y^2w^2.$$
Such multiplications and factorizations were used in the last stages to illustrate that $$2a^2b^2+2a^2c^2 +2b^2c^2-a^4-b^4-c^4= 4s^2p^2.$$

\begin{figure}[htb]
 \begin{center}
\includegraphics[scale=.1]{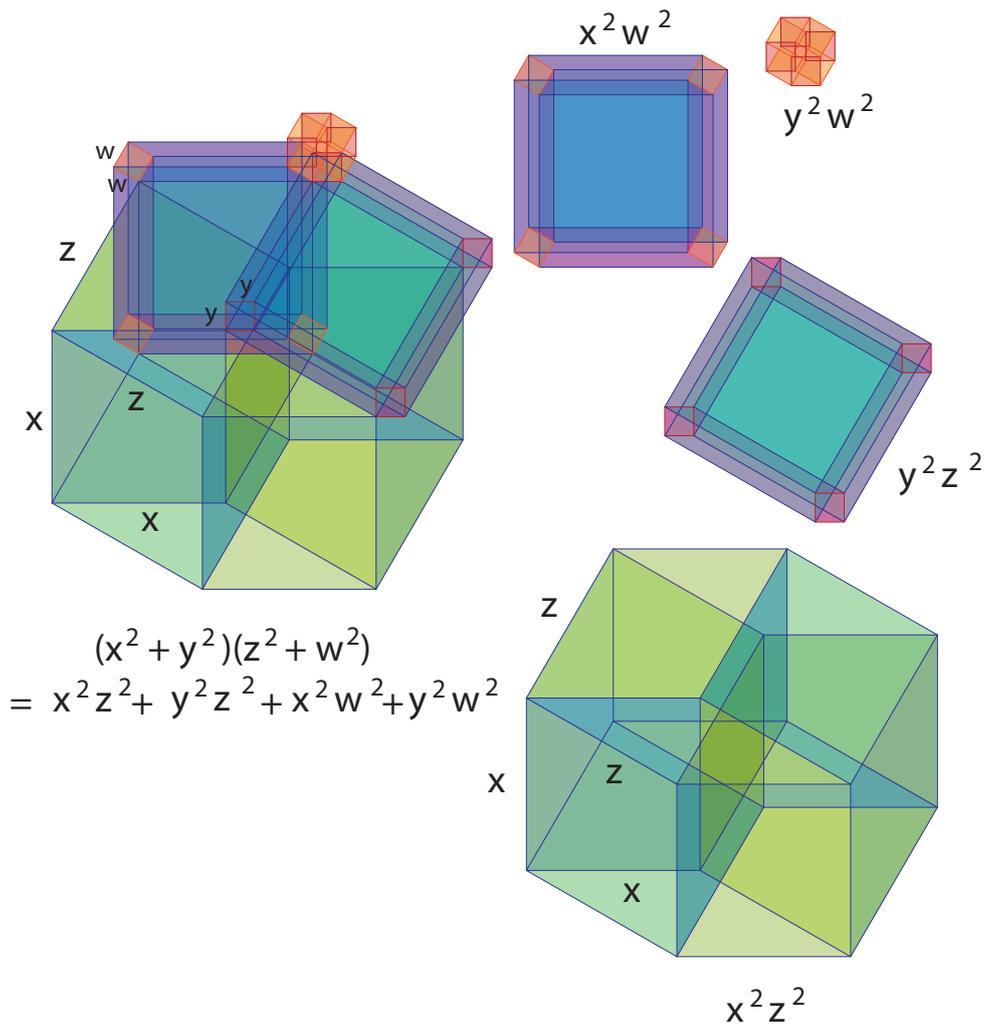}
\end{center}
\caption{Products of sums of squares} 
\label{sumofsq}
\end{figure}

\clearpage

\begin{figure}[htb]
 \begin{center}
\includegraphics[scale=.08]{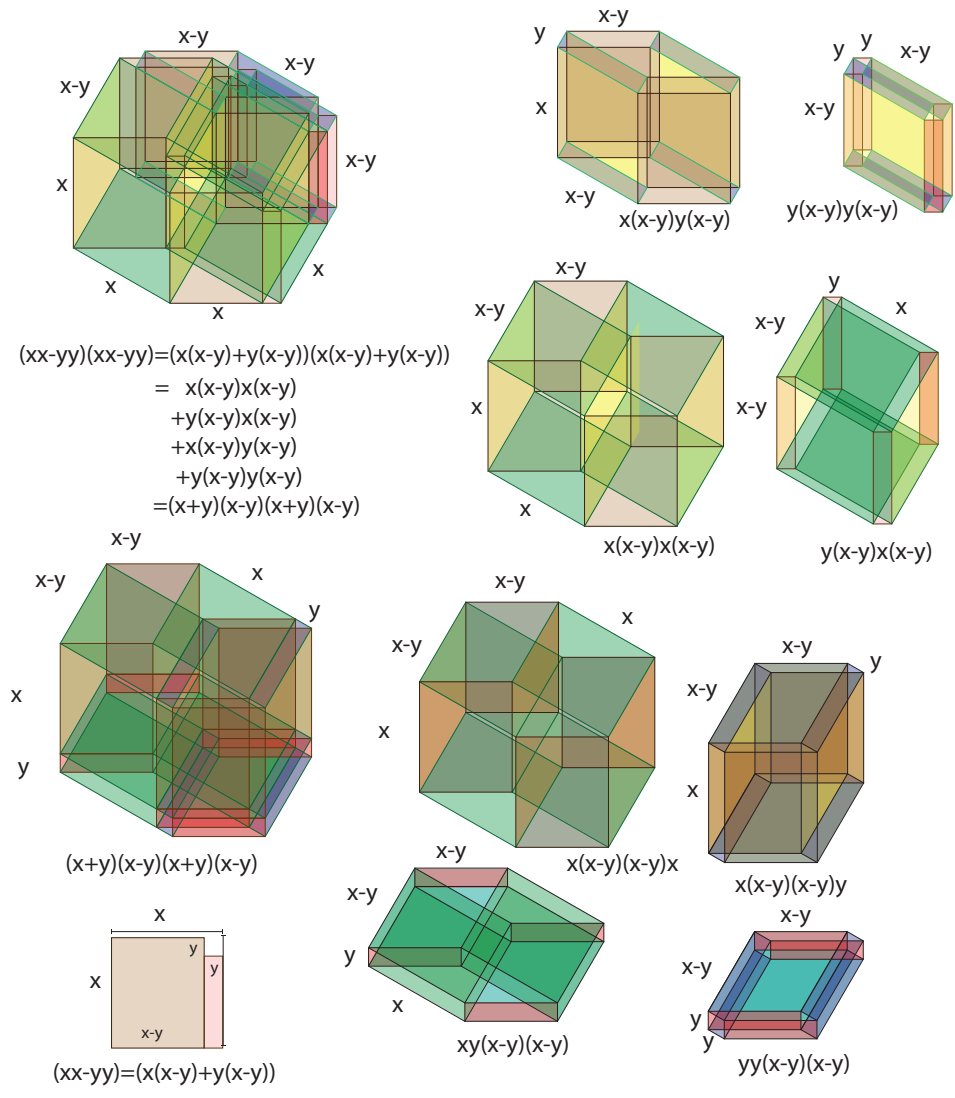}
\end{center}
\caption{Products of differences of squares} 
\label{diffofsq}
\end{figure}

\clearpage

\section{The Left-Hand-Side of the Formula}
\label{LHS}

Consider now the proof that the area of an arbitrary triangle as indicated in Fig.~\ref{righttriangle} is $p\cdot r/2$. The triangle is duplicated and reflected along a horizontal line and these two copies of the original triangle are joined to create a parallelogram. Thus twice the area of the triangle is the area of the parallelogram. Meanwhile, the parallelogram is scissors congruent to a rectangle. 

If we consider the $4$-dimensional figure that is a parallelogram-times-parallelogram, then it can be decomposed into four copies of a triangle times a triangle. The decomposition is analogous to that which is depicted in Fig.~\ref{fourTxT} so we do not indicate a new proof. In the former case, the triangles were right triangles. In general, this decomposition can be composed with a linear map from the rectangle-times-rectangle to the parallelogram-times-parallelogram. 

\begin{figure}[htb]
 \begin{center}
\includegraphics[scale=.08]{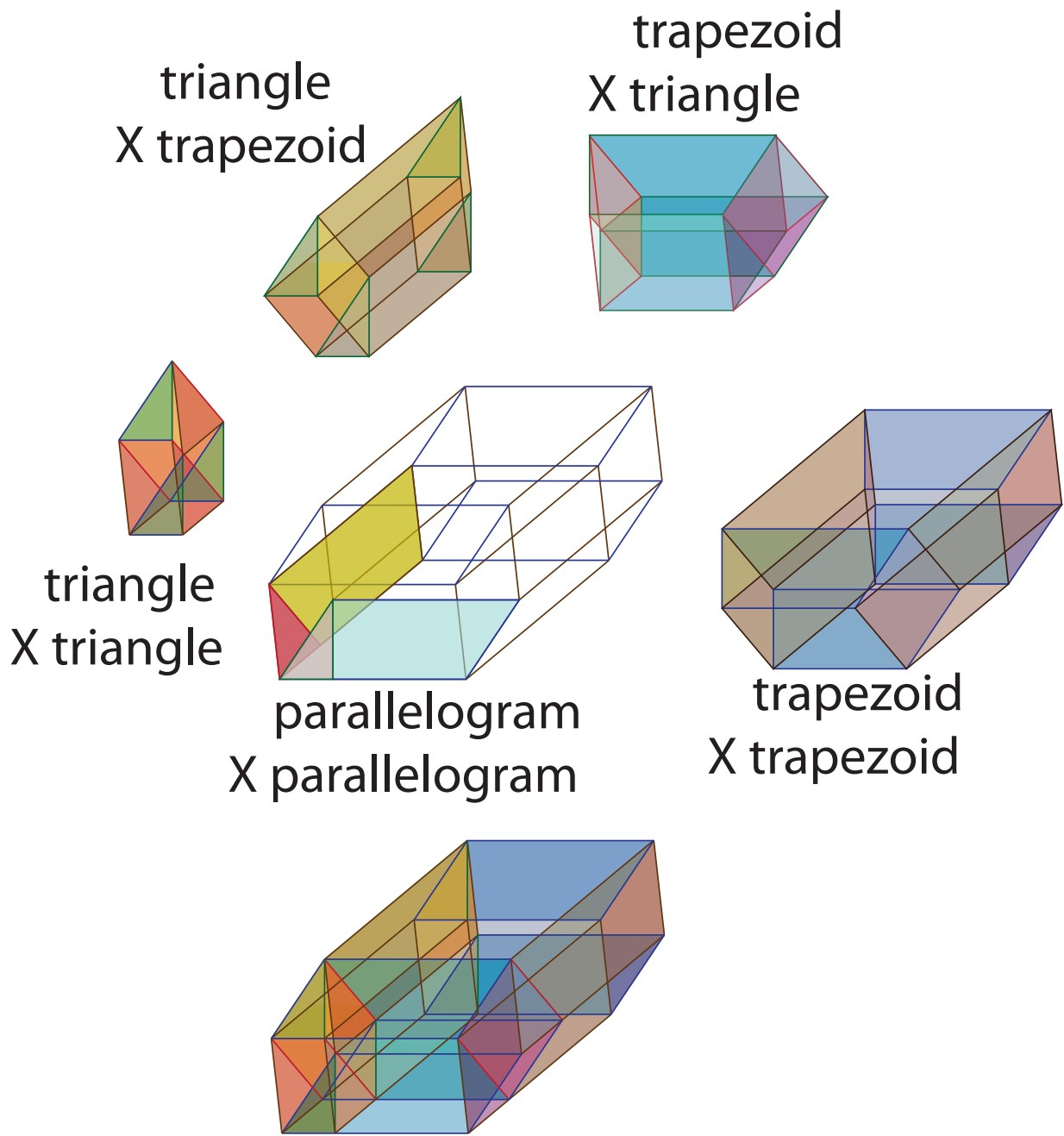}
\end{center}
\caption{The product of a parallelogram and itself} 
\label{paraxpara}
\end{figure}

\begin{figure}[htb]
 \begin{center}
\includegraphics[scale=.08]{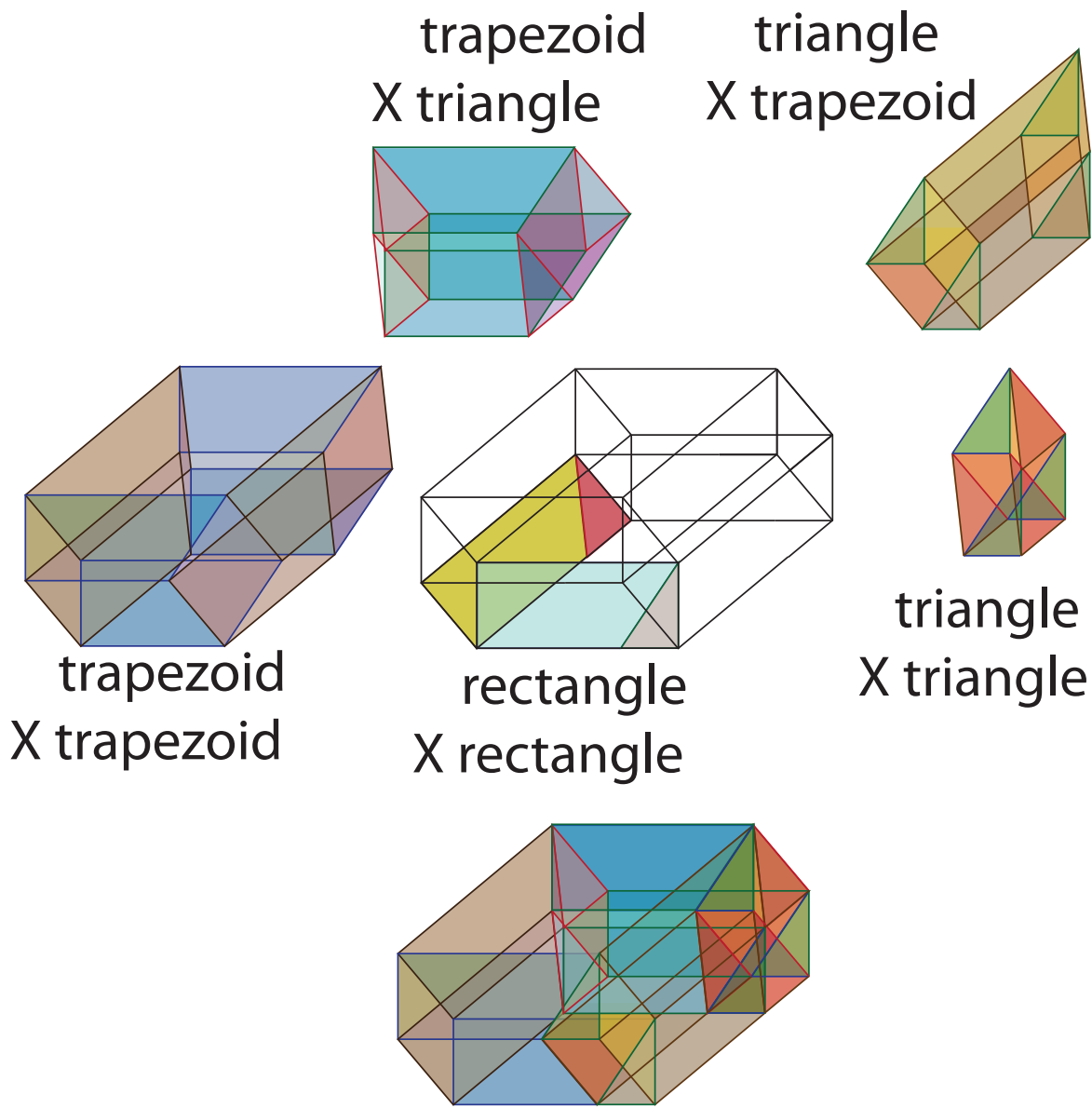}
\end{center}
\caption{Decomposing a rectangle-times=rectangle into the pieces above} 
\label{rectXrect}
\end{figure}

In Fig.~\ref{paraxpara}, the decomposition of a parallelogram into a right triangle and a trapezoid is lifted into $4$-space as the decomposition of a parallelogram-times-parallelogram into four pieces: triangle-times-triangle, triangle-times-trapezoid, trapezoid-times-triangle, and trapezoid-times-trapezoid. 
In Fig.~\ref{rectXrect}, these same four pieces are reassembled to fill a rectangle-times-rectangle.  

The left-hand side of Heron's formula is the quantity $4p^2r^2.$ Sixteen copies of the original triangle times itself are reassembled into four copies of a rectangle-times-rectangle where the dimensions of the rectangle are $p \times r$. Thus the sequence of equalities $$16 A^2$$ $$= 4 p^2r^2 $$ $$= 2a^2b^2+2a^2c^2 +2b^2c^2-a^4-b^4-c^4$$ $$ = (a+b+c)(a+b-c)(a-b+c)(-a+b+c).$$
are all realized as scissors congruences on $4$-dimensional hyper-solids. 

Our $4$-dimensional geometric proof is complete.

\section{Conclusion}
\label{con}

Even after more than two and a half centuries of study, many people --- including many outstanding mathematicians --- have not developed an intuitive grasp of elements of $4$-dimensional geometry. Of course, many of us have a more esoteric understanding of higher geometry that is facilitated by our deep understanding of linear algebra and differential forms.  Our purpose here has been to demonstrate with $2$-dimensional projections and with analytic descriptions, that many aspects of $4$-dimensional geometry can be understood. In addition, we have indicated that a number of combinatorial structures and identities are more deeply understood in terms of aspects of higher dimensional geometry. 

So while an algebraic proof of Heron's formula is tedious, yet conceptually simple, a more geometric proof that depends on the geometry of $4$-spacial dimensions removes some aspects of the tedium. One can envision the pieces fitting tightly together.

Still, there is a cautionary statement to be made. We only perceive $2$-dimensions. Those of us who are sighted see the world projected onto our retinas. As I type these lines, my hands feel the surface of a desk that I presume to be solid. A $4$-dimensional hyper-solid has a boundary that is $3$-dimensional. I can only imagine the surfaces of these pieces. As the $4$-dimensional hyper-solid is projected to $3$-space, its boundary overlaps in the same way that a ball or orange projects to my eye and appears to be a disk. The front and the back of the ball veil the inner structure.

When we use a temporal metaphor to describe the space-time of the room in which we sit, the authors move within the confines of this $10 \times 11 \times 10$ cubic foot room for a period of, say, $1$ hour. The walls, floor and ceiling of the room form a substantial portion of the boundary as they evolve in time. The``us-then" and the ``us-now" --- that is the interior of the room before the paragraph was written and the interior of the room as the paragraph ends forms the rest of the boundary. Our discussion and our gesticulations in completing this last paragraph and all the space around us (as small part of space-time) forms the interior of this hyper-cube of size $1100$ ${\mbox{\rm ft}}^3$-hours. 

Neither the space-time metaphor nor the projection method produce the full richness of a $4$-dimensional entity. But together, they can be used to develop our understanding.

\clearpage

\begin{center}\rule{4.5in}{0.001in} \end{center}

\acknowledgments{Acknowledgments}

The initial work for this project were conducted when the second author was a UCUR (University Committee for Undergraduate Scholarship) student at the University of South Alabama. It appeared in \cite{CM}, but we felt that the pictures and our techniques deserved a wider audience.
We would like to thank Dror Bar-Natan, Abhijit Champanerkar, Dylan Thurston, and the faculty at the University of South Alabama Department of Mathematics and Statistics for valuable conversations.  







\bibliographystyle{mdpi}
\makeatletter
\renewcommand\@biblabel[1]{#1. }
\makeatother

\end{document}